\documentclass[american,english]{article}
\usepackage[T1]{fontenc}
\usepackage[latin9]{inputenc}
\usepackage{geometry}
\usepackage{color}
\usepackage{array}
\usepackage{float}
\usepackage{multirow}
\usepackage{amsmath}
\usepackage{amssymb}
\usepackage{graphicx}
\usepackage{wasysym}

\makeatletter

\newcommand{\noun}[1]{\textsc{#1}}
\providecommand{\tabularnewline}{\\}
\floatstyle{ruled}
\newfloat{algorithm}{tbp}{loa}
\providecommand{\algorithmname}{Algorithm}
\floatname{algorithm}{\protect\algorithmname}

\newcommand{\lyxaddress}[1]{
\par {\raggedright #1
\vspace{1.4em}
\noindent\par}
}

\usepackage{algorithmic}

\def\vec#1{\boldsymbol{#1}}

\usepackage{moreverb}
\usepackage[colorlinks,bookmarksopen,bookmarksnumbered,citecolor=red,urlcolor=red]{hyperref}

\makeatother

\usepackage{babel}
\addto\captionsamerican{\renewcommand{\algorithmname}{Algorithm}}

\begin{document}

\title{A Comparison of Preconditioned Krylov Subspace \\
 Methods for Large-Scale Nonsymmetric Linear Systems}

\author{Aditi Ghai, Cao Lu and Xiangmin Jiao\thanks{E-mail: xiangmin.jiao@stonybrook.edu.}}
\maketitle

\lyxaddress{Department of Applied Mathematics \& Statistics and Institute for
Advanced Computational Science, \\
Stony Brook University, Stony Brook, NY 11794, USA \break }
\begin{abstract}
Preconditioned Krylov subspace (KSP) methods are widely used for solving
large-scale sparse linear systems arising from numerical solutions
of partial differential equations (PDEs). These linear systems are
often nonsymmetric due to the nature of the PDEs, boundary or jump
conditions, or discretization methods. While implementations of preconditioned
KSP methods are usually readily available, it is unclear to users
which methods are the best for different classes of problems. In this
work, we present a comparison of some KSP methods, including GMRES,
TFQMR, BiCGSTAB, and QMRCGSTAB, coupled with three classes of preconditioners,
namely Gauss-Seidel, incomplete LU factorization (including ILUT,
ILUTP, and multilevel ILU), and algebraic multigrid (including BoomerAMG
and ML). Theoretically, we compare the mathematical formulations and
operation counts of these methods. Empirically, we compare the convergence
and serial performance for a range of benchmark problems from numerical
PDEs in 2D and 3D with up to millions of unknowns and also assess
the asymptotic complexity of the methods as the number of unknowns
increases. Our results show that GMRES tends to deliver better performance
when coupled with an effective multigrid preconditioner, but it is
less competitive with an ineffective preconditioner due to restarts.
BoomerAMG with proper choice of coarsening and interpolation techniques
typically converges faster than ML, but both may fail for ill-conditioned
or saddle-point problems while multilevel ILU tends to succeed. We
also show that right preconditioning is more desirable. This study
helps establish some practical guidelines for choosing preconditioned
KSP methods and motivates the development of more effective preconditioners.
\end{abstract}

\section{Introduction}

Preconditioned Krylov subspace (KSP) methods are widely used for solving
large-scale sparse linear systems, especially those arising from numerical
methods for partial differential equations (PDEs). For most modern
applications, these linear systems are nonsymmetric due to various
reasons, such as the multiphysics nature of the PDEs, sophisticated
boundary or jump conditions, or the discretization methods themselves.
For symmetric systems, conjugate gradient (CG) \cite{Hestenes52CG}
and MINRES \cite{Paige75MINRES} are widely recognized as the best
KSP methods \cite{FS12CG}. However, the situation is far less clear
for nonsymmetric systems. Various KSP methods have been developed,
such as GMRES \cite{Saad86GMRES}, CGS \cite{Sonneveld89CGS}, QMR
\cite{FN91QMR}, TFQMR \cite{Freund93TFQMR}, BiCGSTAB \cite{vanderVorst92BiCGSTAB},
QMRCGSTAB \cite{CGS94QMRCGS}, etc. Most of these methods are described
in detail in textbooks such as \cite{BBC94Templates,Saad03IMS,Van-der-Vorst:2003aa},
and their implementations are readily available in software packages,
such as PETSc \cite{petsc-user-ref} and MATLAB \cite{MATLAB}. However,
each of these methods has its advantages and disadvantages. Therefore,
it is difficult for practitioners to choose the proper methods for
their specific applications. Moreover, a KSP method may perform well
with one preconditioner but poorly with another. As a result, users
often spend a significant amount of time to find a reasonable combination
of the KSP methods and preconditioners through trial and error, and
yet the final choice may still be far from optimal. Therefore, a systematic
comparison of the preconditioned KSP methods is an important subject.

In the literature, various comparisons of KSP methods have been reported
previously. In \cite{comparison_trefethen}, Nachtigal, Reddy, and
Trefethen presented some theoretical analysis and comparison of the
convergence properties of CGN (CG on Normal equations), GMRES, and
CGS, which were the leading methods for nonsymmetric systems in the
early 1990s. They showed that the convergence of CGN is governed by
singular values, whereas that of GMRES and CGS by eigenvalues and
pseudo-eigenvalues, and each of these methods may significantly outperform
the others for different matrices. Their work did not consider preconditioners.
The work is also outdated because newer methods have been introduced
since then, which are superior to CGN and CGS. In Saad's textbook
\cite{Saad03IMS}, some comparisons of various KSP methods, including
GMRES, BiCGSTAB, QMR, and TFQMR, were given in terms of computational
cost and storage requirements. The importance of preconditioners was
emphasized, but no detailed comparison for the different combinations
of the KSP methods and preconditioners was given. The same is also
true for other textbooks, such as \cite{Van-der-Vorst:2003aa}. In
terms of empirical comparison, Meister reported a comparison of a
few preconditioned KSP methods for several inviscid and viscous flow
problems \cite{MEISTER1998311}. His study focused on incomplete LU
factorization as the preconditioner. Benzi and coworkers \cite{Benzi99CSS,BENZI02PTL}
also compared a few preconditioners, also with a focus on incomplete
factorization and their block variants. What were notably missing
in these previous studies include the more advanced ILU preconditioners
(such as multilevel ILU \cite{Boll06MPC,lishao10}) and multigrid
preconditioners, which have advanced significantly in recent years.

The goal of this work is to perform a systematic comparison and in
turn establish some practical guidelines in choosing the best preconditioned
KSP solvers. Our study is similar to the recent work of Feng and Saunders
in \cite{FS12CG}, which compared CG and MINRES for symmetric systems.
However, we focus on nonsymmetric systems with a heavier emphasis
on preconditioners. We consider four KSP solvers, GMRES, TFQMR, BiCGSTAB
and QMRCGSTAB. Among these, the latter three enjoy three-term recurrences.
We also consider three classes of general-purpose preconditioners,
namely Gauss-Seidel, incomplete LU factorization (including ILUT,
ILUTP, and multilevel ILU), and algebraic multigrid (including variants
of classical AMG and smoothed aggregation). Each of these KSP methods
and preconditioners has its advantages and disadvantages. At the theoretical
level, we compare the mathematical formulations, operation counts,
and storage requirements of these methods. However, theoretical analysis
alone is insufficient in establishing their suitability for different
types of problems. The primary focus of this work is to compare the
methods empirically in terms of convergence and serial performance
for large linear systems. Our choice of comparing only serial performance
is partially for keeping the focus on the mathematical properties
of these methods, and partially due to the lack of efficient parallel
implementation of ILU.

A systematic comparative study requires a comprehensive set of benchmark
problems. Unfortunately, existing benchmark problems for nonsymmetric
systems, such as those in the Matrix Market \cite{boisvert1997matrix}
and the UF Sparse Matrix Collection (a.k.a. the SuiteSparse Matrix
Collection) \cite{davis2011university}, are generally too small to
be representative of the large-scale problems in current engineering
practice. They often also do not have the right-hand-side vectors,
which can significantly affect the actual performance of KSP methods.
To facilitate this comparison, we constructed a collection of benchmark
systems ourselves from discretization methods for a range of PDEs
in $2$-D and $3$-D. The sizes of these systems range from $10^{5}$
to $10^{7}$ unknowns, which are typical of modern industrial applications,
and are much larger than most benchmark problems in previous studies.
We also assess the asymptotic time complexity of different preconditioned
KSP solvers with respect to the number of unknowns. To the best of
our knowledge, this is the most comprehensive comparison of the preconditioned
KSP solvers to date for large, sparse, nonsymmetric linear systems
in terms of convergence rate, serial performance, and asymptotic complexity.
Our results also show that BoomerAMG in hypre \cite{falgout2002hypre},
which is an extension of classical AMG, typically converges faster
than Trillions/ML \cite{GeeSie06ML}, which is a variant of smoothed-aggregation
AMG. However, it is important to choose the coarsening and interpolation
techniques in BoomerAMG, and we observe that HMIS+FF1 tends to outperform
the default options in both older and newer versions of hypre. However,
both AMG methods tend to fail for very ill-conditioned systems or
saddle-point-like problems, while multilevel ILU, and also ILUTP to
some extent, tend to succeed. We also observed that right preconditioning
is, in general, more reliable than left preconditioning for large-scale
systems. Our results help establish some practical guidelines for
choosing preconditioned KSP methods. They also motivate the further
development of more effective, scalable, and robust multigrid preconditioners. 

The remainder of the paper is organized as follows. In Section~\ref{sec:background},
we review some background knowledge of numerical PDEs, KSP methods,
and preconditioners. In Section~\ref{sec:Analysis-KSP}, we outline
a few KSP methods and compare their main properties in terms of asymptotic
convergence, the number of operations per iteration, and the storage
requirement. This theoretical background will help us predict the
relative performance of the various methods and interpret the numerical
results. In Section~\ref{sec:benchmark-problems}, we describe the
benchmark problems. In Section \ref{sec:Results}, we present numerical
comparisons of the preconditioned KSP methods. Finally, Section~\ref{sec:Conclusions-and-Future}
concludes the paper with some practical recommendations and a discussion
on future work.

\section{Background\label{sec:background}}

In this section, we give a general overview of Krylov subspace methods
and preconditioners for solving a nonsymmetric linear system 
\begin{equation}
\vec{A}\vec{x}=\vec{b},\label{eq:linear_system}
\end{equation}
where $\vec{A}\in\mathbb{R}^{n\times n}$ is large, sparse, nonsymmetric,
and nonsingular, and $\vec{b}\in\mathbb{R}{}^{n}$. These systems
typically arise from PDE discretizations. We consider only real matrices,
because they are more common in applications. However, all the methods
apply to complex matrices, by replacing the matrix transposes with
the conjugate transposes. We focus on the Krylov subspaces and the
procedure in constructing the basis vectors of the subspaces, which
are often the determining factors in the overall performance of different
types of KSP methods. We defer more detailed discussions and analysis
of the individual methods to Section~\ref{sec:Analysis-KSP}.

\subsection{Nonsymmetric Systems from Numerical PDEs}

This work is concerned of solving nonsymmetric systems from discretization
of PDEs. It is important to understand the origins of these systems,
especially of their nonsymmetric structures. Consider an abstract
but general linear, time-independent scalar PDE over $\Omega\subset\mathbb{R}^{d}$,
\begin{equation}
\mathcal{L}u(\vec{x})=f(\vec{x}),\label{eq:linearPDE}
\end{equation}
with Dirichlet or Neumann boundary conditions (BCs) over $\Gamma_{D}$
and $\Gamma_{N}$, respectively, where $d=2$ or $3$, and $\mathcal{L}$
is a linear differential operator, and $f$ is a known source term.
A specific and yet quite general example is the second-order boundary
value problem
\begin{align}
-\vec{\nabla}\cdot(\mu\vec{\nabla}u)+\vec{\nu}\cdot\vec{\nabla}u+\omega^{2}u & =f\quad\,\,\,\,\,\text{ in }\Omega,\label{eq:model_problem}\\
u & =u_{D}\quad\text{ on }\Gamma_{D},\\
\vec{n}\cdot\vec{\nabla}u & =g\quad\,\,\,\,\text{ on }\Gamma_{N},
\end{align}
for which $\mathcal{L}u=-\vec{\nabla}\cdot(\mu\vec{\nabla}u)+\vec{\nu}\cdot\vec{\nabla}u+\omega^{2}u$,
where $\mu$ is scalar field, $\vec{\nu}$ is a vector field, $\omega$
is a scalar, and $\vec{n}$ denotes outward normal to $\Gamma_{N}$.
If $\omega=0$, then it is a \emph{convection-diffusion} (a.k.a. \emph{advection-diffusion})
equation, where $\mu$ corresponds to a diffusion coefficient, and
$\vec{\nu}$ corresponds to a velocity field. If $\vec{\nu}=\vec{0}$,
then it is a \emph{Helmholtz} equation, and $\omega$ typically corresponds
to a wavenumber or frequency.

For ease of discussion, we express the PDE discretization methods
using a general notion of \emph{weighted residuals}. In particular,
consider a set of \emph{test} (a.k.a. \emph{weight}) \emph{functions}
$\{\psi_{i}(\vec{x})\}$. The PDE (\ref{eq:linearPDE}) is then converted
into a set of integral equations
\begin{equation}
\int_{\Omega}\mathcal{L}u(\vec{x})\,\psi_{i}\,d\vec{x}=\int_{\Omega}f(\vec{x})\,\psi_{i}\,d\vec{x}.\label{eq:weak_form}
\end{equation}
To discretize the equations fully, consider a set of \emph{basis functions}
$\{\phi_{ij}(\vec{x})\}$ corresponding to each test function, and
let $u_{i}^{h}\approx\sum_{i}u_{i}\phi_{ij}$ be a local approximation
to $u$ with respect to $\psi_{i}$. Then, we obtain a linear system
\begin{equation}
\vec{A}\vec{u}=\vec{b},\label{eq:linearsys}
\end{equation}
where 
\begin{equation}
a_{ij}=\int_{\Omega}\mathcal{L}\phi_{ij}(\vec{x})\,\psi_{i}(\vec{x})\,d\vec{x}\mbox{ \,\,\ and\,\,\ }b_{i}=\int_{\Omega}f(\vec{x})\,\psi_{i}(\vec{x})\,d\vec{x}.
\end{equation}
This system may be further modified to apply boundary conditions.
In general, the test and basis functions have local support, and therefore
$\vec{A}$ is in general sparse. 

In \emph{finite element methods} (\emph{FEM}) and their variants,
the test functions are typically piecewise linear or higher-degree
polynomials, such as hat functions, and the same set of basis functions
$\{\phi_{j}(\vec{x})\}$ is used regardless of $\psi_{i}$. Since
the test and basis functions are weakly differentiable in FEM, integration
by parts is used to reduce the elliptic operator in $\mathcal{L}$
to a symmetric first-order differential operator in the interior of
the domain; see e.g. \cite{ern2013theory} for details of FEM. If
$\{\phi_{i}\}=\{\psi_{i}\}$, the FEM is a \emph{Galerkin} method,
and $\vec{A}$ is nonsymmetric if $\vec{\nu}\neq0$ in (\ref{eq:model_problem}).
If $\{\phi_{i}\}\neq\{\psi_{i}\}$, the FEM is a \emph{Petrov-Galerkin}
method, and $\vec{A}$ is always nonsymmetric.

In \emph{finite difference methods} (\emph{FDM}), the test functions
are Dirac delta functions at the nodes, and at each node a different
set of polynomial basis functions is constructed from the interpolation
over its stencil; see e.g. \cite{Strikwerda:2004:FDS} for details
of FDM. On uniform structured grids, FDM with centered differences
leads to symmetric matrices for Helmholtz equations with Dirichlet
BCs. However, the matrices are in general nonsymmetric for PDEs with
Neumann BCs, FDM on nonuniform or curvilinear grids, or higher-order
FDM.

While finite differences were traditionally limited to structured
or curvilinear meshes, they are generalized to unstructured meshes
or point clouds in the so-called \emph{generalized finite difference
methods} (\emph{GFDM}); see e.g. \cite{Benito08GFDM}. Like FDM, at
each node GFDM has a set of polynomial basis functions over its stencil,
which are constructed from least squares fittings instead of interpolation.
The matrices from GFDM are always nonsymmetric. A closely method is
AES-FEM \cite{CDJ16OEQ}, of which the test functions are similar
to those of FEM but the basis functions are similar to those of GFDM.
The linear systems from AES-FEM are also nonsymmetric. 

The model problem (\ref{eq:linearPDE}) is a scalar BVP, but it can
be generalized to vector-valued systems of PDEs. In this setting,
$u$ is replaced by a vector field, and $\mu$, $\vec{\nu}$ and $\omega^{2}$
are replaced by tensors, which may be determined by additional unknowns.
Systems of PDEs, such as Stokes equations, lead to saddle-point-like
problems, of which the linear systems have large diagonal blocks.
For time-dependent problems, such as \emph{hyperbolic} or \emph{advection-diffusion}
problems, nonsymmetric linear systems may arise due to finite-element
spatial discretization or implicit time stepping. In particular, hyperbolic
problems are often solved using the \emph{discontinuous Galerkin }(\emph{DG})
or \emph{finite volume} (\emph{FVM}) or methods, of which the test
functions are local polynomials over each element (or cell). A different
set of polynomial basis functions are used per element (or cell),
with flux reconstruction and flux limiters along element (or cell)
boundaries. These methods lead to nonsymmetric systems if implicit
time stepping is used.

\subsection{Krylov Subspaces}

For large-scale sparse linear systems, Krylov-subspace methods are
among the most powerful techniques. Given a matrix $\vec{A}\in\mathbb{R}^{n\times n}$
and a vector $\vec{v}\in\mathbb{R}^{n}$, the $k$th \emph{Krylov
subspace} generated by them, denoted by $\mathcal{K}_{k}(\vec{A},\vec{v})$,
is given by
\begin{equation}
\mathcal{K}_{k}(\vec{A},\vec{v})=\mbox{span}\{\vec{v},\vec{A}\vec{v},\vec{A}^{2}\vec{v},\dots,\vec{A}^{k-1}\vec{v}\}.\label{eq:Krylov}
\end{equation}
To solve (\ref{eq:linear_system}), let $\vec{x}_{0}$ be some initial
guess to the solution, and $\vec{r}_{0}=\vec{b}-\vec{A}\vec{x}_{0}$
is the initial residual vector. A Krylov subspace method incrementally
finds approximate solutions within $\mathcal{K}_{k}(\vec{A},\vec{v})$,
sometimes through the aid of another Krylov subspace $\mathcal{K}_{k}(\vec{A}^{T},\vec{w})$,
where $\vec{v}$ and $\vec{w}$ typically depend on $\vec{r}_{0}$.
To construct the basis of the subspace $\mathcal{K}(\vec{A},\vec{v})$,
two procedures are commonly used: the (restarted) \emph{Arnoldi iteration}
\cite{Arnoldi51PMI}, and the \emph{bi-Lanczos iteration} \cite{Lan50,Van-der-Vorst:2003aa}
(a.k.a. Lanczos biorthogonalization \cite{Saad03IMS} or tridiagonal
biorthogonalization \cite{TB97NLA}).

\subsubsection{The Arnoldi Iteration.}

The Arnoldi iteration is a procedure for constructing orthogonal basis
of the Krylov subspace $\mathcal{K}(\vec{A},\vec{v})$. Starting from
a unit vector $\vec{q}_{1}=\vec{v}/\Vert\vec{v}\Vert$, it iteratively
constructs 
\begin{equation}
\vec{Q}_{k+1}=[\vec{q}_{1}\mid\vec{q}_{2}\mid\dots\mid\vec{q}_{k}\mid\vec{q}_{k+1}]\label{eq:Arnoldi_basis}
\end{equation}
with orthonormal columns by solving 
\begin{equation}
h_{k+1,k}\vec{q}_{k+1}=\vec{A}\vec{q}_{k}-h_{1k}\vec{q}_{1}-\cdots-h_{kk}\vec{q}_{k},\label{eq:Arnoldi_core}
\end{equation}
where $h_{ij}=\vec{q}_{i}^{T}\vec{A}\vec{q}_{j}$ for $j\leq i$,
and $h_{k+1,k}=\Vert\vec{A}\vec{q}_{k}-h_{1k}\vec{q}_{1}-\cdots-h_{kk}\vec{q}_{k}\Vert$,
i.e., the norm of the right-hand side of (\ref{eq:Arnoldi_core}).
This is analogous to Gram-Schmidt orthogonalization. If $\mathcal{K}_{k}\neq\mathcal{K}_{k-1}$,
then the columns of $\vec{Q}_{k}$ form an orthonormal basis of $\mathcal{K}_{k}(\vec{A},\vec{v})$,
and 
\begin{equation}
\vec{A}\vec{Q}_{k}=\vec{Q}_{k+1}\tilde{\vec{H}}_{_{k}},
\end{equation}
where $\tilde{\vec{H}}_{_{k}}$ is a $(k+1)\times k$ upper Hessenberg
matrix, whose entries $h_{ij}$ are those in (\ref{eq:Arnoldi_core})
for $i\leq j+1$, and $h_{ij}=0$ for $i>j+1$.

The KSP method GMRES \cite{Saad86GMRES} is based on the Arnoldi iteration,
with $\vec{v}=\vec{r}_{0}$. If $\vec{A}$ is symmetric, the Hessenberg
matrix $\tilde{\vec{H}}_{_{k}}$ reduces to a tridiagonal matrix,
and the Arnoldi iteration reduces to the \emph{Lanczos iteration}.
The Lanczos iteration enjoys a three-term recurrence. In contrast,
the Arnoldi iteration has a $k$-term recurrence, so its computational
cost increases as $k$ increases. For this reason, one typically needs\textcolor{blue}{{}
}to restart the Arnoldi iteration for large systems (e.g., after every
30 iterations) to build a new Krylov subspace from $\vec{v}=\vec{r}_{k}$
at restart. Unfortunately, the restart may undermine the convergence
of the KSP methods \cite{Saad86GMRES}.

\subsubsection{The Bi-Lanczos Iteration.}

The bi-Lanczos iteration, also known as \emph{Lanczos biorthogonalization}
or \emph{tridiagonal biorthogonalization}, offers an alternative for
constructing the basis of the Krylov subspaces of $\mathcal{K}(\vec{A},\vec{v})$.
Unlike Arnoldi iterations, the bi-Lanczos iterations enjoy a three-term
recurrence. However, the basis will no longer be orthogonal, and we
need to use two matrix-vector multiplications per iteration, instead
of just one.

The bi-Lanczos iterations can be described as follows. Starting from
the vector $\vec{v}_{1}=\vec{v}/\Vert\vec{v}\Vert$, we iteratively
construct
\begin{equation}
\vec{V}_{k+1}=[\vec{v}_{1}\mid\vec{v}_{2}\mid\dots\mid\vec{v}_{k}\mid\vec{v}_{k+1}],\label{eq:nonorth_basis}
\end{equation}
by solving 
\begin{equation}
\beta_{k}\vec{v}_{k+1}=\vec{A}\vec{v}_{k}-\gamma_{k-1}\vec{v}_{k-1}-\alpha_{k}\vec{v}_{k},\label{eq:beta_k}
\end{equation}
analogous to (\ref{eq:Arnoldi_core}). If $\mathcal{K}_{k}\neq\mathcal{K}_{k-1}$,
then the columns of $\vec{V}_{k}$ form a basis of $\mathcal{K}_{k}(\vec{A},\vec{v})$,
and 
\begin{equation}
\vec{A}\vec{V}_{k}=\vec{V}_{k+1}\tilde{\vec{T}}_{_{k}},\label{eq:biorth_A}
\end{equation}
where 
\begin{equation}
\tilde{\vec{T}}_{_{k}}=\begin{bmatrix}\alpha_{1} & \gamma_{1}\\
\beta_{1} & \alpha_{2} & \gamma_{2}\\
 & \beta_{2} & \alpha_{3} & \ddots\\
 &  & \ddots & \ddots & \gamma_{k-1}\\
 &  &  & \beta_{k-1} & \alpha_{k}\\
 &  &  &  & \beta_{k}
\end{bmatrix}
\end{equation}
is a $(k+1)\times k$ tridiagonal matrix. To determine the $\alpha_{i}$
and $\gamma_{i}$, we construct another Krylov subspace $\mathcal{K}(\vec{A}^{T},\vec{w})$,
whose basis is given by the column vectors of 
\begin{equation}
\vec{W}_{k+1}=[\vec{w}_{1}\mid\vec{w}_{2}\mid\dots\mid\vec{w}_{k}\mid\vec{w}_{k+1}],\label{eq:biorth_basis_W}
\end{equation}
subject to the biorthogonality condition 
\begin{equation}
\vec{W}_{k+1}^{T}\vec{V}_{k+1}=\vec{V}_{k+1}^{T}\vec{W}_{k+1}=\vec{I}_{k+1}.\label{eq:biorthogonal}
\end{equation}
Since 
\begin{equation}
\vec{W}_{k+1}^{T}\vec{A}\vec{V}_{k}=\vec{W}_{k+1}^{T}\vec{V}_{k+1}\tilde{\vec{T}}_{_{k}}=\tilde{\vec{T}}_{_{k}},
\end{equation}
it then follows that 
\begin{equation}
\alpha_{k}=\vec{w}_{k}^{T}\vec{A}\vec{v}_{k}.\label{eq:alpha_k}
\end{equation}
Suppose $\vec{V}=\vec{V}_{n}$ and $\vec{W}=\vec{W}_{n}=\vec{V}^{-T}$
form complete basis vectors of $\mathcal{K}_{n}(\vec{A},\vec{v})$
and $\mathcal{K}_{n}(\vec{A}^{T},\vec{w})$, respectively. Let $\vec{T}=\vec{V}^{-1}\vec{A}\vec{V}$
and $\vec{S}=\vec{T}^{T}$. Then, 
\begin{equation}
\vec{W}^{-1}\vec{A}^{T}\vec{W}=\vec{V}^{T}\vec{A}^{T}\vec{V}^{-T}=\vec{T}^{T}=\vec{S},
\end{equation}
and 
\begin{equation}
\vec{A}^{T}\vec{W}_{k}=\vec{W}_{k+1}\tilde{\vec{S}}_{_{k}},\label{eq:biorth_At}
\end{equation}
where $\tilde{\vec{S}}_{k}$ is the leading $(k+1)\times k$ submatrix
of $\vec{S}$. Therefore, 
\begin{equation}
\gamma_{k}\vec{w}_{k+1}=\vec{A}^{T}\vec{w}_{k}-\beta_{k-1}\vec{w}_{k-1}-\alpha_{k}\vec{w}_{k}.\label{eq:gamma_k}
\end{equation}
Starting from $\vec{v}_{1}$ and $\vec{w}_{1}$ with $\vec{v}_{1}^{T}\vec{w}_{1}=1$,
and let $\beta_{0}=\gamma_{0}=1$ and $\vec{v}_{0}=\vec{w}_{0}=\vec{0}$.
Then, $\alpha_{k}$ is uniquely determined by (\ref{eq:alpha_k}),
and $\beta_{k}$ and $\gamma_{k}$ are determined by (\ref{eq:beta_k})
and (\ref{eq:gamma_k}) by up to scalar factors, subject to $\vec{v}_{k+1}^{T}\vec{w}_{k+1}=1$.
A typical choice is to scale the right-hand sides of (\ref{eq:beta_k})
and (\ref{eq:gamma_k}) by scalars of the same modulus \cite[p. 230]{Saad03IMS}.

If $\vec{A}$ is symmetric and $\vec{v}_{1}=\vec{w}_{1}=\vec{v}/\Vert\vec{v}\Vert$,
then the bi-Lanczos iteration reduces to the classical Lanczos iteration
for symmetric matrices. Therefore, it can be viewed as a different
generalization of the Lanczos iteration to nonsymmetric matrices.
Unlike the Arnoldi iteration, the cost of bi-Lanczos iteration is
fixed per iteration, so no restart is ever needed. Some KSP methods,
in particular BiCG \cite{fletcher1976conjugate} and QMR \cite{FN91QMR},
are based on bi-Lanczos iterations. A potential issue of bi-Lanczos
iteration is that it may suffer from \emph{breakdown} if $\vec{v}_{k+1}^{T}\vec{w}_{k+1}=0$
or \emph{near breakdown} if $\vec{v}_{k+1}^{T}\vec{w}_{k+1}\approx0$.
These can be resolved by a \emph{look-ahead} strategy to build a block-tridiagonal
matrix $\vec{T}$. Fortunately, breakdowns are rare, so look-ahead
is rarely implemented. 

A disadvantage of the bi-Lanczos iteration is that it requires the
multiplication with $\vec{A}^{T}$. Although $\vec{A}^{T}$ is in
principle available in most applications, multiplication with $\vec{A}^{T}$
leads to additional difficulties in performance optimization and preconditioning.
Fortunately, in bi-Lanczos iteration, $\vec{V}_{k}$ can be computed
without forming $\vec{W}_{k}$ and vice versa. This observation leads
to the transpose-free variants of the KSP methods, such as TFQMR \cite{Freund93TFQMR},
which is a transpose-free variant of QMR, and CGS \cite{Sonneveld89CGS},
which is a transpose-free variant of BiCG. Two other examples include
BiCGSTAB \cite{vanderVorst92BiCGSTAB}, which is more stable than
CGS, and QMRCGSTAB \cite{CGS94QMRCGS}, which is a hybrid of QMR and
BiCGSTAB, with smoother convergence than BiCGSTAB. These transpose-free
methods enjoy three-term recurrences and require two multiplications
with $\vec{A}$ per iteration. Note that there is not a unique transpose-free
bi-Lanczos iteration. There are primarily two types, used by CGS and
QMR, and by BiCGSTAB and QMRCGSTAB, respectively. We will address
them in more detail in Section~\ref{sec:Analysis-KSP}. 

\subsubsection{Comparison of the Iteration Procedures.}

\begin{table}[tb]
\caption{\label{tab:KrylovSubspaces}Comparisons of KSP methods based on Krylov
subspaces and iteration procedures.}

\selectlanguage{american}%
\centering{}{\small{}}%
\begin{tabular}{c|>{\raggedright}p{2cm}|c|c|c}
\hline 
\multirow{2}{*}{{\small{}Method}} & \multirow{2}{2cm}{{\small{}Iteration}} & \multicolumn{2}{c|}{{\small{}Matrix-Vector Prod.}} & \multirow{2}{*}{{\small{}Recurrence}}\tabularnewline
\cline{3-4} 
 &  & {\small{}A}\textbf{\small{}$^{T}$} & {\small{}A} & \tabularnewline
\hline 
\hline 
\selectlanguage{english}%
{\small{}GMRES \cite{Saad86GMRES}}\selectlanguage{american}%
 & {\small{}Arnoldi} & {\small{}0} & {\small{}1} & {\small{}$k$}\tabularnewline
\hline 
\selectlanguage{english}%
BiCG \cite{fletcher1976conjugate}\selectlanguage{american}%
 & \multirow{2}{2cm}{\foreignlanguage{english}{{\small{}bi-Lanczos}}} & \multirow{2}{*}{{\small{}1}} & \multirow{2}{*}{{\small{}1}} & \multirow{6}{*}{\selectlanguage{english}%
3\selectlanguage{american}%
}\tabularnewline
\cline{1-1} 
\selectlanguage{english}%
{\small{}QMR \cite{FN91QMR}}\selectlanguage{american}%
 &  &  &  & \tabularnewline
\cline{1-4} 
\selectlanguage{english}%
CGS \cite{Sonneveld89CGS}\selectlanguage{american}%
 & \multirow{2}{2cm}{{\small{}transpose-free bi-Lanczos 1}} & \multirow{4}{*}{{\small{}~~~~~0~~~~~}} & \multirow{4}{*}{{\small{}2}} & \tabularnewline
\cline{1-1} 
\selectlanguage{english}%
{\small{}TFQMR \cite{Freund93TFQMR}}\selectlanguage{american}%
 &  &  &  & \tabularnewline
\cline{1-2} 
\selectlanguage{english}%
{\small{}BiCGSTAB \cite{vanderVorst92BiCGSTAB}}\selectlanguage{american}%
 & \multirow{2}{2cm}{{\small{}transpose-free bi-Lanczos 2}} &  &  & \tabularnewline
\cline{1-1} 
\selectlanguage{english}%
{\small{}QMRCGSTAB \cite{CGS94QMRCGS}}\selectlanguage{american}%
 &  &  &  & \tabularnewline
\hline 
\end{tabular}{\small\par}\selectlanguage{english}%
\end{table}

Both the Arnoldi iteration and the bi-Lanczos iteration are based
on the Krylov subspace $\mathcal{K}(\vec{A},\vec{r}_{0})$. However,
these iteration procedures have very different properties, which are
inherited by their corresponding KSP methods, as summarized in Table~\ref{tab:KrylovSubspaces}.
These properties, for the most part, determine the cost per iteration
of the KSP methods. For KSP methods based on the Arnoldi iteration,
at the $k$th iteration the residual $\vec{r}_{k}=\mathcal{P}_{k}(\vec{A})\vec{r}_{0}$
for some degree-$k$ polynomial $\mathcal{P}_{k}$, so the asymptotic
convergence rates primarily depend on the eigenvalues and the generalized
eigenvectors in the Jordan form of $\vec{A}$ \cite{comparison_trefethen,Saad03IMS}.
For methods based on transpose-free bi-Lanczos, in general $\vec{r}_{k}=\hat{\mathcal{P}}_{k}(\vec{A})\vec{r}_{0}$,
where $\hat{\mathcal{P}}_{k}$ is a polynomial of degree $2k$. Therefore,
the convergence of these methods also depends on the eigenvalues and
generalized eigenvectors of $\vec{A}$, but at different asymptotic
rates. Typically, the reduction of error in one iteration of a bi-Lanczos-based
KSP method is approximately equal to that of two iterations in an
Arnoldi-based KSP method. Since the Arnoldi iteration requires only
one matrix-vector multiplication per iteration, compared to two per
iteration for the bi-Lanczos iteration, the costs of different KSP
methods are comparable in terms of the number of matrix-vector multiplications. 

Theoretically, the Arnoldi iteration is more robust because of its
use of orthogonal basis, whereas the bi-Lanczos iteration may breakdown
if $\vec{v}_{k+1}^{T}\vec{w}_{k+1}=0$. However, the Arnoldi iteration
typically requires restarts, which can undermine convergence. In general,
if the iteration count is small compared to the average number of
nonzeros per row, the methods based on the Arnoldi iteration may be
more efficient; if the iteration count is large, the cost of orthogonalization
in Arnoldi iteration may become higher than that of bi-Lanczos iteration.
For these reasons, conflicting results are often reported in the literature.
However, the apparent disadvantages of each KSP method may be overcome
by effective preconditioners: For Arnoldi iterations, if the KSP method
converges before restart is needed, then it may be the most effective
method; for bi-Lanczos iterations, if the KSP method converges before
any breakdown, it is typically more robust than the methods based
on restarted Arnoldi iterations. We will review the preconditioners
in the next subsection.

Note that some KSP methods use a Krylov subspace other than $\mathcal{K}(\vec{A},\vec{r}_{0})$.
The most notable examples are LSQR \cite{Paige92LSQR} and LSMR \cite{Fong11LSMR},
which use the Krylov subspace $\mathcal{K}(\vec{A}^{T}\vec{A},\vec{A}^{T}\vec{r}_{0})$.
These methods are mathematically equivalent to applying CG or MINRES
to the normal equation, respectively, but with better numerical properties
than CGN. An advantage of these methods is that they apply to least
squares systems without modification. However, they are not transpose-free,
they tend to converge slowly for square linear systems, and they require
special preconditioners. For these reasons, we do not include them
in this study.

\subsection{Preconditioners}

The convergence of KSP methods can be improved significantly by the
use of preconditioners. Various preconditioners have been proposed
for Krylov subspace methods over the past few decades. It is virtually
impossible to consider all of them. For this comparative study, we
focus on three classes of preconditioners, which are representative
for the state-of-the-art ``black-box'' preconditioners: Gauss-Seidel,
incomplete LU factorization, and algebraic multigrid.

\subsubsection{Left versus Right Preconditioning.}

Roughly speaking, a preconditioner is a matrix or transformation $\vec{M}$,
whose inverse $\vec{M}^{-1}$ approximates $\vec{A}^{-1}$, and $\vec{M}^{-1}\vec{v}$
can be computed efficiently. For nonsymmetric linear systems, a preconditioner
may be applied either to the left or the right of $\vec{A}$. With
a left preconditioner, instead of solving (\ref{eq:linear_system}),
one solves the linear system 
\begin{equation}
\vec{M}^{-1}\vec{A}\vec{x}=\vec{M}^{-1}\vec{b}
\end{equation}
by utilizing the Krylov subspace $\mathcal{K}(\vec{M}^{-1}\vec{A},\vec{M}^{-1}\vec{b})$
instead of $\mathcal{K}(\vec{A},\vec{b})$. For a right preconditioner,
one solves the linear system 
\begin{equation}
\vec{A}\vec{M}^{-1}\vec{y}=\vec{b}
\end{equation}
by utilizing the Krylov subspace $\mathcal{K}(\vec{A}\vec{M}^{-1},\vec{b})$,
and then $\vec{x}=\vec{M}^{-1}\vec{y}$. The convergence of a preconditioned
KSP method is then determined by the eigenvalues of $\vec{M}^{-1}\vec{A}$,
which are the same as those of $\vec{A}\vec{M}^{-1}$, as well as
their generalized eigenvectors. Qualitatively, $\vec{M}$ is a good
preconditioner if $\vec{M}^{-1}\vec{A}$ (or $\vec{A}\vec{M}^{-1}$
for right preconditioning) is not too far from normal and its eigenvalues
are more clustered than those of $\vec{A}$ \cite{TB97NLA}. However,
this is more useful as a guideline for developers of preconditioners,
rather than for users. 

\begin{figure}[h]
\begin{minipage}[t]{0.45\textwidth}%
\begin{center}
\includegraphics[width=1\textwidth]{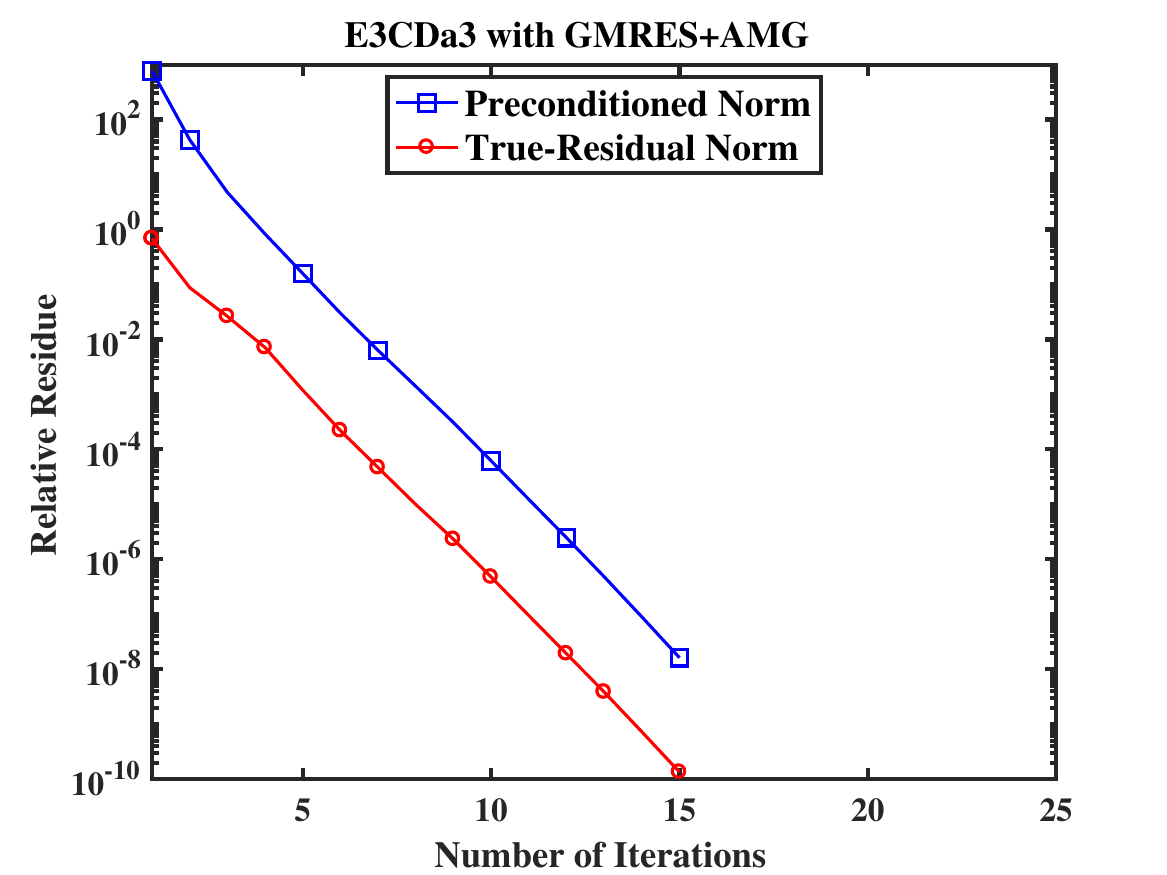}
\par\end{center}%
\end{minipage}\hfill{} %
\begin{minipage}[t]{0.45\textwidth}%
\begin{center}
\includegraphics[width=1\textwidth]{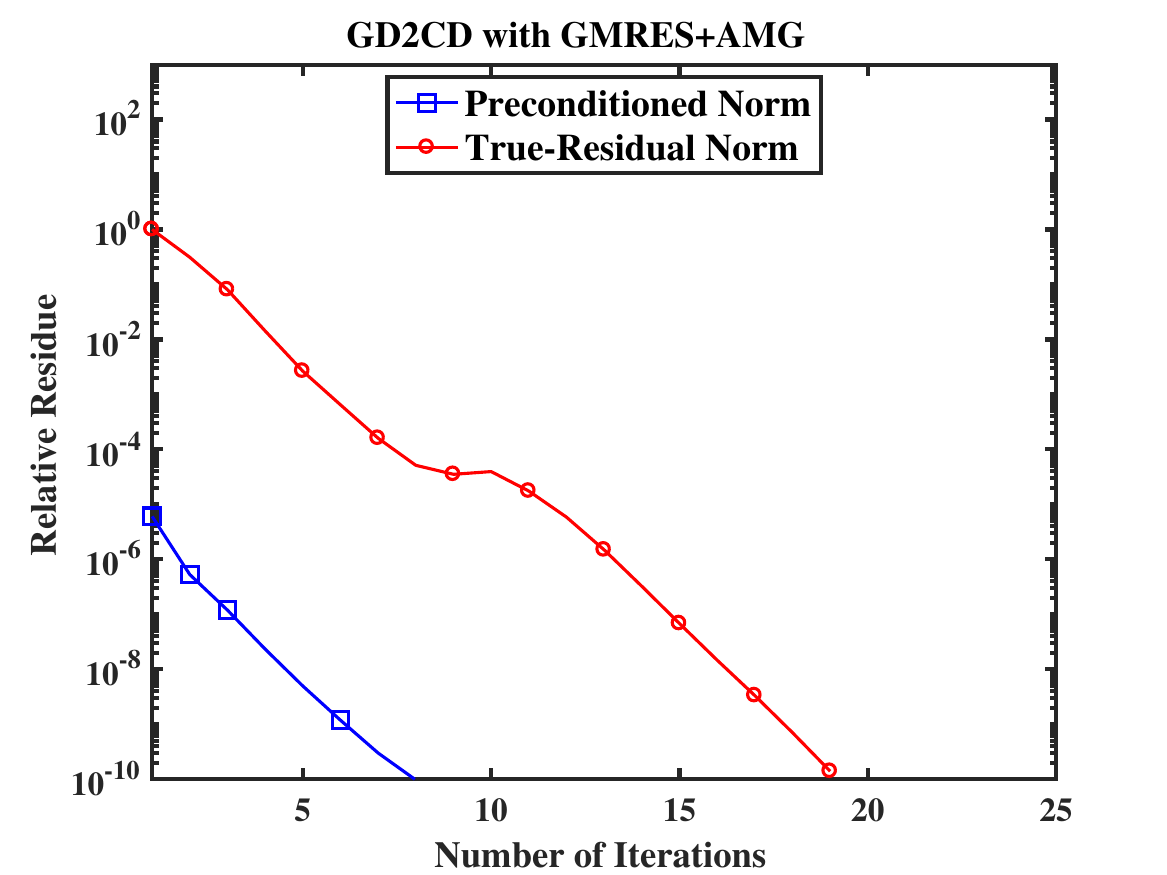}
\par\end{center}%
\end{minipage}

\caption{\label{fig:preconditioned-vs-true-residuals}\textcolor{black}{Examples
where left preconditioners lead to large discrepancies of preconditioned
and true residuals by orders of magnitude and in turn delayed (left)
or premature termination (right).}}
\end{figure}

Although the left and right-preconditioners have similar asymptotic
behavior, they can behave drastically differently in practice. It
is advisable to use right, instead of left, preconditioners for two
reasons. First, the termination criterion of a Krylov subspace method
is typically based on the norm of the residual of the preconditioned
system, which may differ significantly from the true residual. Figure~\ref{fig:preconditioned-vs-true-residuals}
shows two examples, where the norms of the preconditioned residuals
are significantly larger and smaller than the true residuals, respectively;
we will explain these matrices in Section~\ref{sec:benchmark-problems}.
Second, if the iteration terminates with a relatively large residual
$\vec{r}$, then the error of the solution is bounded by 
\begin{equation}
\left\Vert \delta\vec{x}\right\Vert \leq\left\Vert \vec{M}\vec{A}^{-1}\right\Vert \left\Vert \vec{M}^{-1}\vec{r}\right\Vert \leq\kappa(\vec{M})\left\Vert \vec{A}^{-1}\right\Vert \left\Vert \vec{r}\right\Vert ,\label{eq:error-bound}
\end{equation}
which may differ significantly from $\left\Vert \vec{A}^{-1}\right\Vert \left\Vert \vec{r}\right\Vert $
if $\kappa(\vec{M})\gg1$. The stability analysis of PDE discretization
typically depends on the boundedness of $\left\Vert \vec{A}^{-1}\right\Vert $
in (\ref{eq:linearsys}), so one should not change the residual unless
the preconditioner is derived based on \textit{a priori} knowledge
of the PDE discretization. One could overcome these issues by computing
the true residual $\left\Vert \vec{r}\right\Vert $ at each step,
but it would incur additional costs with a left preconditioner. Since
the preconditioners that we consider are algebraic in nature, we use
only right preconditioners in this study.

Note that in the so-called \emph{symmetric preconditioners}, the Cholesky
factorization of $\vec{M}^{-1}$ is applied symmetrically to both
the left and right of $\vec{A}$. A well-known example is the symmetric
successive over-relaxation (SSOR) \cite{habetler1961symmetric}. Such
preconditioners preserve symmetry for symmetric matrices. However,
since they also alter the norm of the residual and our focus is on
nonsymmetric matrices, we do not consider symmetric preconditioners
in this study.

\subsubsection{Gauss-Seidel and SOR.}

Gauss-Seidel and its generalization SOR (successive over-relaxation)
are some of the simplest preconditioners. Based on stationary iterative
methods, Gauss-Seidel and SOR are relatively easy to implement, require
virtually no setup time (at least in serial), and are sometimes fairly
effective. Therefore, they are often good choices if one needs to
implement a preconditioner from scratch.

Consider the partitioning $\vec{A}=\vec{D}+\vec{L}+\vec{U}$, where
$\vec{D}$ is the diagonal of $\vec{A}$, $\vec{L}$ is the strictly
lower triangular part, and $\vec{U}$ is the strictly upper triangular
part. Given $\vec{x}_{k}$ and $\vec{b}$, the Gauss-Seidel method
computes a new approximation to $\vec{x}_{k+1}$ as 
\begin{equation}
\vec{x}_{k+1}=(\vec{D}+\vec{L})^{-1}(\vec{b}-\vec{U}\vec{x}_{k}).
\end{equation}
SOR generalizes Gauss-Seidel by introducing a relaxation parameter
$\omega$. It computes $\vec{x}_{k+1}$ as 
\begin{equation}
\vec{x}_{k+1}=(\vec{D}+\omega\vec{L})^{-1}\left(\omega\left(\vec{b}-\vec{U}\vec{x}_{k}\right)+(1-\omega)\vec{D}\vec{x}_{k}\right).
\end{equation}
When $\omega=1$, SOR reduces to Gauss-Seidel; when $\omega>1$ and
$\omega<1$, it corresponds to over-relaxation and under-relaxation,
respectively. We choose to include Gauss-Seidel instead of SOR in
our comparison, because it is parameter free, and an optimal choice
of $\omega$ in SOR is problem dependent. Another related preconditioner
is the Jacobi or diagonal preconditioner, which is less effective
than Gauss-Seidel. A limitation of Gauss-Seidel, also shared by Jacobi
and SOR, is that the diagonal entries of $\vec{A}$ must be nonzero.

\subsubsection{\label{subsec:Incomplete-LU-Factorization.}Incomplete LU Factorization.}

Incomplete LU factorization (ILU) performs an approximate factorization
\begin{equation}
\vec{A}\approx\tilde{\vec{L}}\tilde{\vec{U}},
\end{equation}
where $\tilde{\vec{L}}$ and $\tilde{\vec{U}}$ are far sparser than
those in the LU factorization of $\vec{A}$. This approximate factorization
is typically computed in a preprocessing step. In the preconditioned
Krylov solver, $\vec{M}^{-1}\vec{y}$ is computed by forward solve
$\vec{z}=\tilde{\vec{L}}^{-1}\vec{y}$ and then back substitution
$\tilde{\vec{U}}^{-1}\vec{z}$.

There are several variants of ILU factorization. In its simplest form,
ILU does not introduce any fill, so that $\tilde{\vec{L}}$ and $\tilde{\vec{U}}$
preserve the sparsity patterns of the lower and upper triangular parts
of $\vec{A}$, respectively. This approach is often referred to as
\emph{ILU0} or ILU(0). ILU0 may be extended to preserve some of the
fills based on their \emph{levels} in the elimination tree. This is
often referred to as \emph{ILU}$(k)$, which zeros out all the fills
of level $k+1$ or higher. In addition, \emph{ILU}$(k)$ may be further
combined with thresholding on the numerical values, resulting in \emph{ILU
with dual thresholding} (\emph{ILUT}) \cite{saad1994ilut}. Most implementations
of ILU, such as those in PETSc and hypre, use some variants of ILUT,
where PETSc also allows the user to control the number of fills.

A serious issue with ILUT is that it may breakdown if there are many
zeros in the diagonal. This issue can be mitigated by pre-permuting
the matrix, but it may still fail in practice for saddle-point-like
problems. A more effective approach is to use partial pivoting, which
results in the so-called \emph{ILUTP} \cite{Saad03IMS}. The ILU implementations
in MATLAB \cite{MATLAB}, SPARSKIT \cite{Saad1994Sparsekit}, and
SuperLU \cite{lishao10}, for example, are based on ILUTP. 

A major drawback of ILUTP is that the number of nonzeros in the $\tilde{\vec{L}}$
and $\tilde{\vec{U}}$ factors may grow superlinearly with respect
to the original number of nonzeros if the drop tolerance is small.
This leads to superlinear growth of setup times and solve times. However,
a small drop tolerance may be needed for robustness, especially for
very large sparse systems. As a result, parameter tuning for ILUTP
may become an impossible task. This problem is mitigated by the \emph{multilevel
ILU}, or\emph{ MILU }for short. Unlike ILUTP, MILU uses diagonal pivoting
instead of partial pivoting, and it permutes rows and columns that
would have led to ill-conditioned $\tilde{\vec{L}}$ and $\tilde{\vec{U}}$
to the end and delays factorizing them \cite{Boll06MPC}. This approach
leads to much better robustness with a relatively small number of
fills. A robust, serial implementation of MILU is available in ILUPACK
\cite{ilupack}. Typically, MILU scales linearly with respect to the
original number of nonzeros, so it is effective for large sparse linear
systems, especially in serial. We will report some numerical comparisons
of different variants of ILU in Section~\ref{subsec:ILU-Comparison}. 

\subsubsection{\label{subsec:Algebraic-Multigrid}Algebraic Multigrid.}

Multigrid methods, including \emph{geometric multigrid} (\emph{GMG})
and \emph{algebraic multigrid} (\emph{AMG}), are the most sophisticated
preconditioners. These methods accelerate stationary iterative methods
(or more precisely, the so-called \emph{smoothers} in multigrid methods)
by constructing a series of coarser representations. The key difference
between GMG and AMG is the coarsening strategies and the associated
prolongation (a.k.a. interpolation) and restriction (a.k.a. projection)
operators between different levels. Compared to Gauss-Seidel and ILU,
multigrid preconditioners, especially GMG, are far more difficult
to implement. Fortunately, AMG preconditioners are more readily accessible
through software libraries. Similar to ILU, AMG typically requires
a significant amount of pre-processing time. Computationally, AMG
is more expensive than Gauss-Seidel and ILU in terms of both setup
time and cost per iteration, but they can accelerate convergence much
more significantly.

There are primarily two types of AMG methods: the \emph{classical
AMG} \cite{ruge1987algebraic}, and \emph{smoothed aggregation} \cite{vanvek1996algebraic}.
The main difference between the two is in their coarsening strategies.
The classical AMG uses some variants of the so-called ``classical
coarsening'' \cite{ruge1987algebraic}. It separates all the points
into \emph{coarse points}, which form the next coarser level, and
\emph{fine points}, which are interpolated from the coarse points.
In smoothed aggregation, the coarsening is done by accumulating the
aggregates, which then form the coarse-grid points \cite{vanvek1996algebraic}.
The differences in coarsening also lead to differences in their corresponding
prolongation and restriction operators. In this work, we consider
both classical and smoothed-aggregation AMG, in particular BoomerAMG
(part of hypre \cite{falgout2002hypre}) and ML \cite{GeeSie06ML},
which are variants of the two types, respectively. Both BoomerAMG
and ML are accessible through PETSc \cite{petsc-user-ref}. 

It is worth noting that BoomerAMG has many coarsening and interpolation
strategies, which may lead to drastically different performance, especially
for 3D problems. We give a brief overview of three coarsening techniques,
namely \emph{Falgout}, \emph{PMIS}, and \emph{HMIS}. Falgout coarsening
is a variant of the standard Ruge-Stüben coarsening \cite{ruge1987algebraic},\textcolor{blue}{{}
}adapted for parallel implementations. Falgout coarsening gives good
results for 2D problems, but for 3D problems, the stencil size may
grow rapidly, resulting in poor efficiency and asymptotic complexity
as the problem size increases. To overcome this issue, PMIS and HMIS
were introduced in \cite{de2006reducing}. PMIS, or \emph{Parallel
Modified Independent Set}, is a variant of the parallel maximal independent
set algorithm \cite{luby1986simple,jones1993parallel}. HMIS, or \emph{Hybrid
Modified Independent Set}, is a hybrid of PMIS and the classical Ruge-Stüben
scheme. HMIS and PMIS can significantly improve the runtime efficiency
for large 3D problems, especially when coupled with a proper interpolation
technique \cite{hypre-user,hypre-reference}.

There are various interpolation techniques available for BoomerAMG.
We consider three of them: \textit{\textcolor{black}{classical, extended+i}},
and \textit{\textcolor{black}{FF1}}. The classical interpolation is
a distance-one interpolation formulated by Ruge and Stüben \cite{ruge1987algebraic}.
This interpolation works well when combined with Falgout coarsening
for 2D problems. Up to hypre version 2.11.1, Falgout coarsening with
classical interpolation was the default for BoomerAMG. The classical
interpolation can be extended to includes coarse-points that are distance
two from the fine point to be interpolated. The extended+i interpolation
is such an extension, with some additional sophistication in computing
the weights \cite{de2008distance}. This interpolation works well
for 3D problems with HMIS and PMIS coarsening. In hypre version 2.11.2,
PMIS coarsening with extended+i interpolation is the default for BoomerAMG
\cite{hypre-user}. The \emph{FF1} interpolation, or the \textit{\textcolor{black}{modified
FF interpolation}}\textit{\textcolor{black}{\emph{ \cite{de2004coarsening}}}},
is another extension of the classical interpolation, but it includes
only one distance-two coarse point when a strong fine-fine point connection
is encountered \cite{de2008distance}. 

An important consideration of the interpolation schemes is the \emph{stencil
size}, defined as the average number of coefficients per matrix row.
Too large stencils can significantly increase the setup time and also
runtime, especially for 3D problems. However, too small stencils may
not capture enough information and adversely affect convergence rate.
Therefore, special care must be taken in choosing the coarsening strategy,
which we will discuss further in Section \ref{subsec:HYPRE-vs-ML}.
hypre recommends truncating the interpolation operator to four to
five elements per row \cite[p. 43-46]{hypre-user}. In Section \ref{subsec:HYPRE-vs-ML},
we will present a comparison of BoomerAMG versus ML as well as a comparison
of coarsening and interpolation in BoomerAMG. Our results show that
hypre tends to outperform ML, at least in serial, and FF1 tends to
outperform extended+i for BoomerAMG.

Besides the coarsening and interpolation, another important aspect
of multigrid methods is the \emph{smoothers}, which smooth the solutions
of the residual equations at each level. For both AMG and GMG, the
smoothers are typically based on stationary iterative methods (such
as Jacobi or Gauss-Seidel) or incomplete LU (such as ILU0). In particular,
the default smoother in hypre is SOR/Jacobi. Since Gauss-Seidel and
ILU0 have difficulties with saddle-point-like problems as preconditioners,
multigrid methods with these smoothers share similar issues.

\section{Analysis of Preconditioned KSP Methods\label{sec:Analysis-KSP}}

In this section, we discuss a few Krylov subspace methods in more
detail, especially the preconditioned GMRES, TFQMR, BiCGSTAB, and
QMRCGSTAB with right preconditioners. In the literature, these methods
are typically given either without preconditioners or with left preconditioners.
We present their high-level descriptions with right preconditioners.
We also present some theoretical results in terms of operation counts
and storage, which are helpful in analyzing the numerical results.

\subsection{GMRES}

Developed by Saad and Schultz \cite{Saad86GMRES}, \emph{GMRES}, or
\emph{generalized minimal residual method}, is one of most well-known
iterative methods for solving large, sparse, nonsymmetric systems.
GMRES is based on the Arnoldi iteration. At the $k$th iteration,
it minimizes $\Vert\vec{r}_{k}\Vert$ in $\mathcal{K}_{k}(\vec{A},\vec{b})$.
Equivalently, it finds an optimal degree-$k$ polynomial $\mathcal{P}_{k}(\vec{A})$
such that $\vec{r}_{k}=\mathcal{P}_{k}(\vec{A})\vec{r}_{0}$ and $\Vert\vec{r}_{k}\Vert$
is minimized. Suppose the approximate solution has the form
\begin{equation}
\vec{x}_{k}=\vec{x}_{0}+\vec{Q}_{k}\vec{z},\label{eq:GMRES_sol}
\end{equation}
where $\vec{Q}_{k}$ was given in (\ref{eq:Arnoldi_basis}). Let $\beta=\Vert\vec{r}_{0}\Vert$
and $\vec{q}_{1}=\vec{r}_{0}/\Vert\vec{r}_{0}\Vert$. It then follows
that 
\begin{equation}
\vec{r}_{k}=\vec{b}-\vec{A}\vec{x}_{k}=\vec{b}-\vec{A}(\vec{x}_{0}+\vec{Q}_{k}\vec{z})=\vec{r}_{0}-\vec{A}\vec{Q}_{k}\vec{z}=\vec{Q}_{k+1}(\beta\vec{e}_{1}-\tilde{\vec{H}}_{k}\vec{z}),\label{eq:GMRES_residual}
\end{equation}
and $\Vert\vec{r}_{k}\Vert=\Vert\beta\vec{e}_{1}-\tilde{\vec{H}}_{k}\vec{z}\Vert$.
Therefore, $\vec{r}_{k}$ is minimized by solving the least squares
system $\tilde{\vec{H}}_{k}\vec{z}\approx\beta\vec{e}_{1}$ using
QR factorization. In this sense, GMRES is closely related to MINRES
for solving symmetric systems \cite{Paige75MINRES}. Algorithm~1
gives a high-level pseudocode of the preconditioned GMRES with a right
preconditioner; for a more detailed pseudocode, see e.g. \cite[p. 284]{Saad03IMS}. 

For nonsingular matrices, the convergence of GMRES depends on whether
$\vec{A}$ is close to normal, and also on the distribution of its
eigenvalues \cite{comparison_trefethen,TB97NLA}. At the $k$th iteration,
GMRES requires one matrix-vector multiplication, $k+1$ axpy operations
(i.e., $\alpha\vec{x}+\vec{y}$), and $k+1$ inner products. Let $\ell$
denote the average number of nonzeros per row. In total, GMRES requires
$2n(\ell+2k+2)$ floating-point operations per iteration and requires
storing $k+5$ vectors in addition to the matrix itself. Due to the
high cost of orthogonalization in the Arnoldi iteration, GMRES in
practice needs to be restarted periodically. This leads to GMRES with
restart, denoted by GMRES($r$), where $r$ is the iteration count
before restart. A typical value of $r$ is 30, which is the recommended
value for large systems in PETSc, ILUPACK, etc.

Note that the GMRES implementation in MATLAB (as of R2018a) supports
only left preconditioning. Its implementation in PETSc supports both
left and right preconditioning, but the default is left preconditioning.
An extension of GMRES, called \emph{Flexible GMRES} (or \emph{FGMRES})
\cite[p. 287]{Saad03IMS}, allows adapting the preconditioner from
iteration to iteration, and it only supports right preconditioners.
Therefore, if FGMRES is available, one can use it as GMRES with right
preconditioning by fixing the preconditioner across iterations.

\begin{algorithm}[h]
\begin{minipage}[t]{0.45\textwidth}%
\textbf{\noun{Algorithm}}\textbf{ 1:}

\textbf{Right-Precond'dGMRES}

\textbf{input}: $\vec{x}_{0}$: initial guess

\hspace{1.25cm}$\vec{r}_{0}$: initial residual

\textbf{output}: $\vec{x}_{*}$: final solution

\begin{algorithmic}[1]

\STATE $\vec{q}_{1}\leftarrow\vec{r}_{0}/\Vert\vec{r}_{0}\Vert$;
$\beta\leftarrow\Vert\vec{r}_{0}\Vert$

\STATE \textbf{for} $k=1,2,\dots$

\STATE ~~~~obtain $\tilde{\vec{H}}_{k}$ and $\vec{Q}_{k}$ from
Arnoldi iteration s.t. $\vec{r}_{k}=\mathcal{P}_{k}(\vec{A}\vec{M}^{-1})\vec{r}_{0}$

\STATE ~~~~solve $\tilde{\vec{H}}_{k}\vec{z}\approx\beta\vec{e}_{1}$

\STATE ~~~~$\vec{y}_{k}\leftarrow\vec{Q}_{k}\vec{z}$

\STATE ~~~~check convergence of $\Vert\vec{r}_{k}\Vert$

\STATE \textbf{end for}

\STATE $\vec{x}_{*}\leftarrow\vec{M}^{-1}\vec{y}_{k}$

\end{algorithmic}%
\end{minipage}\hfill{} %
\begin{minipage}[t]{0.45\textwidth}%
\textbf{\noun{Algorithm}}\textbf{ 2:}

\textbf{Right-Precond'd TFQMR}

\textbf{input}: $\vec{x}_{0}$: initial guess

\hspace{1.25cm}$\vec{r}_{0}$: initial residual

\textbf{output}: $\vec{x}_{*}$: final solution

\begin{algorithmic}[1]

\STATE $\vec{v}_{1}\leftarrow\vec{r}_{0}/\Vert\vec{r}_{0}\Vert$;
$\beta\leftarrow\Vert\vec{r}_{0}\Vert$

\STATE \textbf{for} $k=1,2,\dots$

\STATE ~~~~obtain $\tilde{\vec{T}}_{k}$ and $\vec{V}_{k}$ from
bi-Lanczos s.t. $\vec{r}_{k}=\mathcal{\tilde{P}}_{k}^{2}(\vec{A}\vec{M}^{-1})\vec{r}_{0}$

\STATE ~~~~solve $\tilde{\vec{T}}_{k}\vec{z}\approx\beta\vec{e}_{1}$

\STATE ~~~~$\vec{y}_{k}\leftarrow\vec{V}_{k}\vec{z}$

\STATE ~~~~check convergence of $\Vert\vec{r}_{k}\Vert$

\STATE \textbf{end for}

\STATE $\vec{x}_{*}\leftarrow\vec{M}^{-1}\vec{y}_{k}$

\end{algorithmic}%
\end{minipage}
\end{algorithm}

\subsection{QMR and TFQMR}

Proposed by Freund and Nachtigal \cite{FN91QMR}, \emph{QMR}, or \emph{quasi-minimal
residual method}, minimizes $\vec{r}_{k}$ in a pseudo-norm within
the Krylov subspace $\mathcal{K}(\vec{A},\vec{r}_{0})$. At the $k$th
step, suppose the approximate solution has the form
\begin{equation}
\vec{x}_{k}=\vec{x}_{0}+\vec{V}_{k}\vec{z},\label{eq:22}
\end{equation}
where $\vec{V}_{k}$ was the same as that in (\ref{eq:nonorth_basis}).
Let $\beta=\Vert\vec{r}_{0}\Vert$ and $\vec{v}_{1}=\vec{r}_{0}/\Vert\vec{r}_{0}\Vert$.
It then follows that 
\begin{equation}
\vec{r}_{k}=\vec{b}-\vec{A}\vec{x}_{k}=\vec{b}-\vec{A}(\vec{x}_{0}+\vec{V}_{k}\vec{z})=\vec{r}_{0}-\vec{A}\vec{V}_{k}\vec{z}=\vec{V}_{k+1}(\beta\vec{e}_{1}-\tilde{\vec{T}}_{k}\vec{z}).
\end{equation}
QMR minimize $\Vert\beta\vec{e}_{1}-\tilde{\vec{T}}_{k}\vec{z}\Vert$
by solving the least-squares problem $\tilde{\vec{T}}_{k}\vec{z}\approx\beta\vec{e}_{1}$,
which is equivalent to minimizing the pseudo-norm 
\begin{equation}
\Vert\vec{r}_{k}\Vert_{\vec{W}_{k+1}^{T}}=\Vert\vec{W}_{k+1}^{T}\vec{r}_{k}\Vert_{2},
\end{equation}
where $\vec{W}_{k+1}$ was defined in (\ref{eq:biorth_basis_W}).

QMR requires explicit constructions of $\vec{W}_{k}$. \emph{TFQMR}
\cite{Freund93TFQMR} is a transpose-free variant, which constructs
$\vec{V}_{k}$ without forming $\vec{W}_{k}$. Motivated by CGS, \cite{Sonneveld89CGS},
at the $k$th iteration, TFQMR finds a degree-$k$ polynomial $\tilde{\mathcal{P}}_{k}(\vec{A})$
such that $\vec{r}_{k}=\mathcal{\tilde{P}}_{k}^{2}(\vec{A})\vec{r}_{0}$.
This is what we refer to as ``transpose-free bi-Lanczos 1'' in Table~\ref{tab:KrylovSubspaces}.
Algorithm~2 outlines TFQMR with a right preconditioner. Its only
difference from GMRES is in lines 3\textendash 5. Detailed pseudocode
without preconditioners can be found in \cite{Freund93TFQMR} and
\cite[p. 252]{Saad03IMS}.

At the $k$th iteration, TFQMR requires two matrix-vector multiplication,
ten axpy operations (i.e., $\alpha\vec{x}+\vec{y}$), and four inner
products. In total, TFQMR requires $4n(\ell+7)$ floating-point operations
per iteration and requires storing eight vectors in addition to the
matrix itself. This is comparable to QMR, which requires 12 axpy operations
and two inner products, so QMR requires the same number of floating-point
operations. However, QMR requires storing twice as many vectors as
TFQMR. In practice, TFQMR often outperforms QMR, because the multiplication
with $\vec{A}^{T}$ is often less optimized. In addition, preconditioning
QMR is problematic, especially with multigrid preconditioners. Therefore,
TFQMR is in general preferred over QMR. Both QMR and TFQMR may suffer
from breakdowns, but they rarely happen in practice, especially with
a good preconditioner. Unlike GMRES, the TFQMR implementation in MATLAB
supports only right preconditioning.

\subsection{BiCGSTAB}

Proposed by van der Vorst \cite{vanderVorst92BiCGSTAB}, \emph{BiCGSTAB}
is a transpose-free version of BiCG, which has smoother convergence
than BiCG and CGS. Different from CGS and TFQMR, at the $k$th iteration,
BiCGSTAB constructs another degree-$k$ polynomial 
\begin{equation}
\mathcal{Q}_{k}(\vec{A})=(1-\omega_{1}\vec{A})(1-\omega_{2}\vec{A})\cdots(1-\omega_{k}\vec{A})\label{eq:bicgstab_poly}
\end{equation}
in addition to $\mathcal{\tilde{P}}_{k}(\vec{A})$ in CGS, such that
$\vec{r}_{k}=\mathcal{Q}_{k}(\vec{A})\mathcal{\tilde{P}}_{k}(\vec{A})\vec{r}_{0}$.
BiCGSTAB determines $\omega_{k}$ by minimizing $\Vert\vec{r}_{k}\Vert$
with respect to $\omega_{k}$. This is what we referred to as ``transpose-free
bi-Lanczos 2'' in Table~\ref{tab:KrylovSubspaces}. Like BiCG and
CGS, BiCGSTAB solves the linear system $\vec{T}_{k}\vec{z}=\beta\vec{e}_{1}$
using LU factorization without pivoting, which is analogous to solving
the tridiagonal system using Cholesky factorization in CG \cite{Hestenes52CG}.
Algorithm~3 outlines BiCGSTAB with a right preconditioner, of which
the only difference from GMRES is in lines 3\textendash 5. Detailed
pseudocode without preconditioners can be found in \cite[p. 136]{Van-der-Vorst:2003aa}.

At the $k$th iteration, BiCGSTAB requires two matrix-vector multiplications,
six axpy operations, and four inner products. In total, it requires
$4n(\ell+5)$ floating-point operations per iteration and requires
storing $10$ vectors in addition to the matrix itself. Like GMRES,
the convergence rate of BiCGSTAB also depends on the distribution
of the eigenvalues of $\vec{A}$. Unlike GMRES, however, BiCGSTAB
is ``parameter-free.'' Its underlying bi-Lanczos iteration may breakdown,
but it is very rare in practice with a good preconditioner. Therefore,
BiCGSTAB is often more efficient and robust than restarted GMRES.
Like TFQMR, the BiCGSTAB implementation in MATLAB (as of R2018a) supports
only right preconditioning.

\begin{algorithm}[h]
\begin{minipage}[t]{0.45\textwidth}%
\textbf{\noun{Algorithm}}\textbf{ 3:}

\textbf{Right-Precond'd BiCGSTAB}

\textbf{input}: $\vec{x}_{0}$: initial guess

\hspace{1.25cm}$\vec{r}_{0}$: initial residual

\textbf{output}: $\vec{x}_{*}$: final solution

\begin{algorithmic}[1]

\STATE $\vec{v}_{1}\leftarrow\vec{r}_{0}/\Vert\vec{r}_{0}\Vert$;
$\beta\leftarrow\Vert\vec{r}_{0}\Vert$

\STATE \textbf{for} $k=1,2,\dots$

\STATE ~~~~obtain $\vec{T}_{k}$ and $\vec{V}_{k}$ from bi-Lanczos
s.t. $\vec{r}_{k}=\mathcal{Q}_{k}\mathcal{\tilde{P}}_{k}(\vec{A}\vec{M}^{-1})\vec{r}_{0}$

\STATE ~~~~solve $\vec{T}_{k}\vec{z}=\beta\vec{e}_{1}$

\STATE ~~~~$\vec{y}_{k}\leftarrow\vec{V}_{k}\vec{z}$

\STATE ~~~~check convergence of $\Vert\vec{r}_{k}\Vert$

\STATE \textbf{end for}

\STATE $\vec{x}_{*}\leftarrow\vec{M}^{-1}\vec{y}_{k}$

\end{algorithmic}%
\end{minipage}\hfill{} %
\begin{minipage}[t]{0.45\textwidth}%
\textbf{\noun{Algorithm 4}}\textbf{:}

\textbf{Right-Precond'd QMRCGSTAB}

\textbf{input}: $\vec{x}_{0}$: initial guess

\hspace{1.25cm}$\vec{r}_{0}$: initial residual

\textbf{output}: $\vec{x}_{*}$: final solution

\begin{algorithmic}[1]

\STATE $\vec{v}_{1}\leftarrow\vec{r}_{0}/\Vert\vec{r}_{0}\Vert$;
$\beta\leftarrow\Vert\vec{r}_{0}\Vert$

\STATE \textbf{for} $k=1,2,\dots$

\STATE ~~~~obtain $\tilde{\vec{T}}_{k}$ and $\vec{V}_{k}$ from
bi-Lanczos s.t. $\vec{r}_{k}=\mathcal{Q}_{k}\mathcal{\tilde{P}}_{k}(\vec{A}\vec{M}^{-1})\vec{r}_{0}$

\STATE ~~~~solve $\tilde{\vec{T}}_{k}\vec{z}\approx\beta\vec{e}_{1}$

\STATE ~~~~$\vec{y}_{k}\leftarrow\vec{V}_{k}\vec{z}$

\STATE ~~~~check convergence of $\Vert\vec{r}_{k}\Vert$

\STATE \textbf{end for}

\STATE $\vec{x}_{*}\leftarrow\vec{M}^{-1}\vec{y}_{k}$

\end{algorithmic}%
\end{minipage}
\end{algorithm}

\subsection{QMRCGSTAB}

One drawback of BiCGSTAB is that the residual does not decrease monotonically,
and it is often quite oscillatory. Chan ${\it el\thinspace al.}$
\cite{CGS94QMRCGS} proposed \emph{QMRCGSTAB}, which is a hybrid of
QMR and BiCGSTAB, to improve the smoothness of BiCGSTAB. Like BiCGSTAB,
QMRCGSTAB constructs a polynomial $\mathcal{Q}_{k}(\vec{A})$ as defined
in (\ref{eq:bicgstab_poly}) by minimizing $\Vert\vec{r}_{k}\Vert$
with respect to $\omega_{k}$, which they refer to as ``local quasi-minimization.''
Like QMR, it then minimizes $\Vert\vec{W}_{k+1}^{T}\vec{r}_{k}\Vert_{2}$
by solving the least-squares problem $\tilde{\vec{T}}_{k}\vec{z}\approx\beta\vec{e}_{1}$,
which they refer to as ``global quasi-minimization.'' Algorithm~4
outlines the high-level algorithm, of which its main difference from
BiCGSTAB is in lines 3 and 4. Detailed pseudocode without preconditioners
can be found in \cite{CGS94QMRCGS}.

At the $k$th iteration, QMRCGSTAB requires two matrix-vector multiplications,
eight axpy operations, and six inner products. In total, it requires
$4n(\ell+7)$ floating-point operations per iteration, and it requires
storing 13 vectors in addition to the matrix itself. Like QMR and
BiCGSTAB, the underlying bi-Lanczos may breakdown, but it is very
rare in practice with a good preconditioner. There is no built-in
implementation of QMRCGSTAB in MATLAB or PETSc. In this study, we
implement the algorithm ourselves with right preconditioning.

\subsection{Comparison of Operation Counts and Storage}

We summarize the cost and storage comparison of the four KSP methods
in Table~\ref{tab:Comparison-of-operations}. Except for GMRES, the
other methods require two matrix-vector multiplications per iteration.
However, we should not expect GMRES to be twice as fast as the other
methods, because the asymptotic convergence of the KSP methods depends
on the degrees of the polynomials, which is equal to the number of
matrix-vector products instead of the number of iterations. Therefore,
the reduction of error in one iteration of the other methods is approximately
equal to that of two iterations in GMRES. However, since GMRES minimizes
the 2-norm of the residual in the Krylov subspace if no restart is
performed, it may converge faster than the other methods in terms
of the number of matrix-vector multiplications. Therefore, GMRES may
indeed by the most efficient, especially with an effective preconditioner.
However, without an effective preconditioner, the restarted GMRES
may converge slowly and even stagnate for large systems. For the three
methods based on bi-Lanczos, computing the $2$-norm of the residual
for convergence checking requires an extra inner product. Among the
three methods, BiCGSTAB is the most efficient, requiring $8n$ fewer
floating-point operations per iteration than TFQMR and QMRCGSTAB.
In Section~\ref{sec:Results}, we will present numerical comparisons
of the different methods to complement this theoretical analysis.
\begin{table}[tb]
\caption{\label{tab:Comparison-of-operations}Comparison of operations per
iteration and memory requirements of various KSP methods. $n$ denotes
the number of rows, $\ell$ the average number of nonzeros per row,
and $k$ the iteration count. }

\centering{}%
\begin{tabular}{c|c|c|c|c|c|c}
\hline 
\multirow{2}{*}{{\small{}Method}} & \multirow{2}{*}{Minimization} & {\small{}Mat-vec} & \multirow{2}{*}{{\small{}axpy}} & {\small{}Inner} & \multirow{2}{*}{{\small{}FLOPs}} & Stored\tabularnewline
 &  & {\small{} Prod.} &  & {\small{} Prod.} &  & vectors\tabularnewline
\hline 
\hline 
{\small{}GMRES} & $\Vert\vec{r}_{k}\Vert$ & {\small{}1} & {\small{}$k$+1} & {\small{}$k$+1} & {\small{}$2n(\ell+2k+2)$} & {\small{}$k+5$}\tabularnewline
\hline 
{\small{}BiCGSTAB} & $\Vert\vec{r}_{k}\Vert$ w.r.t. $\omega_{k}$ & \multirow{4}{*}{{\small{}2}} & {\small{}6} & {\small{}4} & {\small{}$4n(\ell+5)$} & {\small{}10}\tabularnewline
\cline{1-2} \cline{4-7} 
{\small{}TFQMR } & $\Vert\vec{r}_{k}\Vert_{\vec{W}_{k+1}^{T}}$ &  & {\small{}10} & {\small{}4} & \multirow{3}{*}{{\small{}$4n(\ell+7)$}} & {\small{}$8$}\tabularnewline
\cline{1-2} \cline{4-5} \cline{7-7} 
\multirow{2}{*}{{\small{}QMRCGSTAB }} & $\Vert\vec{r}_{k}\Vert$ w.r.t. $\omega_{k}$  &  & \multirow{2}{*}{{\small{}8}} & \multirow{2}{*}{{\small{}6}} &  & \multirow{2}{*}{{\small{}$13$}}\tabularnewline
 & \& $\Vert\vec{r}_{k}\Vert_{\vec{W}_{k+1}^{T}}$ &  &  &  &  & \tabularnewline
\hline 
\end{tabular}
\end{table}

In terms of storage, TFQMR requires the least amount of memory. BiCGSTAB
requires two more vectors than TFQMR, and QMRCGSTAB requires three
more vectors than BiCGSTAB. GMRES requires the most amount of memory
when $k\apprge8$. These storage requirements are typically not large
enough to be a concern in practice. 

\subsection{\label{subsec:Analysis-of-Preconditioners}Cost of Preconditioned
Krylov Subspace Methods}

Our analysis above did not consider preconditioners. The computational
cost of Gauss-Seidel and ILU0 is approximately equal to one matrix-vector
multiplication per iteration. However, the analysis of other variants
of ILU is more complicated. Allowing a moderate drop tolerance can
significantly improve the robustness of ILU, but as explained in Section~\ref{subsec:Incomplete-LU-Factorization.},
the number of fills may grow superlinearly with respect to problem
size, especially with pivoting. We will present a comparison of different
variants of ILU in Section \ref{subsec:ILU-Comparison}.

The cost analysis of the multigrid preconditioner is far more complicated.
In general, the runtime performance is dominated by smoothing on the
finest level, which is typically\textcolor{blue}{{} }a few sweeps of
matrix-vector multiplication, depending on how many times the smoother
is called per iteration. In AMG, the setup time may be significant,
which is dominated by the coarsening step and the construction of
the prolongation and restriction operators.

\section{\label{sec:benchmark-problems}Benchmark Problems}

For comparative studies, the selection of benchmark problems is important.
Unfortunately, most existing benchmark problems for nonsymmetric systems,
such as those in the Matrix Market \cite{boisvert1997matrix} and
the UF Sparse Matrix Collection \cite{davis2011university}, are generally
very small, so they are not representative of the large-scale problems
used in current engineering practice. In addition, they often do not
have the right-hand-side vectors, which can significantly affect the
performance of KSP methods.

For this study, we collected a set of large benchmark problems from
a range of PDEs in $2$-D and $3$-D. Table~\ref{tab:test_matrices}
summarizes the IDs, corresponding PDE, sizes, and estimated condition
numbers of 29 matrices. These cases were selected from a much larger
number of systems that we have collected and tested. The numbers of
unknowns range from about $10^{5}$ to $10^{7}$, which are representative
of engineering applications. The condition numbers are in 1-norm,
estimated using MATLAB's \textsf{condest} function. They were unavailable
for some largest matrices because the computations ran out of memory
on our computer cluster. 

We identify each matrix by its discretization type, spatial dimension,
type of PDE, etc. In particular, the first letter or first two letters
indicate the discretization methods, where \emph{E} stands for \emph{FEM},
\emph{AE} for \emph{AES-FEM}, \emph{GD} for \emph{GFD}, \emph{DG}
for \emph{discontinuous Galerkin}, \emph{D} for finite difference
methods (\emph{FDM}), and \emph{V} for \emph{finite volume methods}
(FVM). It is followed by the dimension (2 or 3) of the domain. The
next two or three capital letters indicate the type of the PDEs, where
\emph{CD} stands for \emph{convection-diffusion}, \emph{HM} for \emph{Helmholtz},
\emph{INS} for \emph{incompressible Navier-Stokes}, \emph{CNS} for
\emph{compressible Navier-Stokes}, \emph{MP} for \emph{mixed Poisson},
and \emph{STH} for \emph{Stokes} with \emph{Taylor-Hood elements}.
A lower-case letter may follow to indicate different types of boundary
conditions or parameters. If different mesh sizes of the same problems
were used, we append a digit to the matrix ID, where a larger digit
corresponds to a finer mesh. 

\begin{table}[h]
\caption{\label{tab:test_matrices}Summary of test matrices. }

\centering{}%
\begin{tabular}{c|c|c|c|c}
\hline 
\textbf{Matrix ID} & \textbf{PDE} & \textbf{Size } & \textbf{\#Nonzeros} & \textbf{Cond. No.}\tabularnewline
\hline 
\hline 
E2CDa & \multirow{16}{*}{conv. diff.} & 1,044,226 & 7,301,314 & 8.31e5\tabularnewline
\cline{1-1} \cline{3-5} 
E2CDb &  & 9,006,001 & 62,964,007 & 1.28e7\tabularnewline
\cline{1-1} \cline{3-5} 
E2CDc &  & 9,006,001 & 62,964,007 & 1.34e7\tabularnewline
\cline{1-1} \cline{3-5} 
E2CDd &  & 9,006,001 & 62,940,021 & 3.98e8\tabularnewline
\cline{1-1} \cline{3-5} 
E2CDe &  & 9,006,001 & 62,940,021 & 3.99e9\tabularnewline
\cline{1-1} \cline{3-5} 
E3CDa1 &  & 237,737 & 1,819,743 & 8.90e3\tabularnewline
\cline{1-1} \cline{3-5} 
E3CDa2 &  & 1,529,235 & 23,946,925  & 3.45e4\tabularnewline
\cline{1-1} \cline{3-5} 
E3CDa3 &  & 13,110,809 & 197,881,373 & $-$\tabularnewline
\cline{1-1} \cline{3-5} 
E3CDb &  & 4,173,281 & 59,843,949 & 1.72e5\tabularnewline
\cline{1-1} \cline{3-5} 
E3CDc1 &  & 4,173,281 & 59,843,949 & 1.38e6\tabularnewline
\cline{1-1} \cline{3-5} 
E3CDc2 &  & 16,974,593 & 253,036,801 & -\tabularnewline
\cline{1-1} \cline{3-5} 
AE2CD &  & 1,044,226  & 13,487,418 & 9.77e5\tabularnewline
\cline{1-1} \cline{3-5} 
AE3CD &  & 13,110,809 & 197,882,439 & $-$\tabularnewline
\cline{1-1} \cline{3-5} 
GD2CD &  & 1,044,226  & 7,476,484 & 2.38e6\tabularnewline
\cline{1-1} \cline{3-5} 
GD3CD &  & 1,529,235 & 23,948,687 & 6.56e4\tabularnewline
\cline{1-1} \cline{3-5} 
DG2CD &  & 1,033,350 & 12,385,260 & 1.80e7\tabularnewline
\hline 
D2HMa & \multirow{3}{*}{ Helmholtz} & 1,340,640 & 6,694,058 & 7.23e8\tabularnewline
\cline{1-1} \cline{3-5} 
D2HMb &  & 1,320,000 & 6,586,400 & 4.3e5\tabularnewline
\cline{1-1} \cline{3-5} 
D2HMc &  & 1,320,000 & 6,586,400 & 2.18\tabularnewline
\hline 
E2INSa1 & \multirow{5}{*}{incomp. Navier-Stokes} & 982,802  & 19,578,988 & 2.9e5 \tabularnewline
\cline{1-1} \cline{3-5} 
E2INSa2 &  & 1,232,450 & 24,562,751 & 3.7e5\tabularnewline
\cline{1-1} \cline{3-5} 
E2INSb &  & 9,971,828 & 199,293,332 & 1.7e2\tabularnewline
\cline{1-1} \cline{3-5} 
E2INSc &  & 9,963,354 & 199,124,203 & 6.9e2\tabularnewline
\cline{1-1} \cline{3-5} 
E3INS &  & 3,090,903 & 234,996,071 & $-$\tabularnewline
\hline 
V3CNS & comp. Navier-Stokes & 381,689 & 37,464,962 & 2.2e11\tabularnewline
\hline 
\multicolumn{5}{c}{\emph{Saddle-point-like problems}}\tabularnewline
\hline 
E3STH1 & \multirow{2}{*}{Stokes} & 29,114 & 2,638,666 & 5.32e6\tabularnewline
\cline{1-1} \cline{3-5} 
E3STH2 &  & 859,812 & 82,754,416 & $-$\tabularnewline
\hline 
E3MP1 & \multirow{2}{*}{mixed Poisson} & 343,200 & 7,269,600 & 1.75e5\tabularnewline
\cline{1-1} \cline{3-5} 
E3MP2 &  & 1,150,200 & 24,510,600 & $-$\tabularnewline
\hline 
\end{tabular}
\end{table}

Except for V3NS, which was from the UF Sparse Matrix Collection, we
generated the systems by ourselves using a combination of in-house
implementations (especially for (G)FDM and AES-FEM) and the FEniCS
software \cite{alnaes2015fenics} (especially for mixed elements and
DG). For completeness, we describe the origins of these matrices in
terms of their equations, parameters, and discretization methods as
follows.

\paragraph{Convection-Diffusion Equation with FEM.}

\begin{figure}[h]
\begin{minipage}[t]{0.45\textwidth}%
\begin{center}
\includegraphics[width=1\columnwidth]{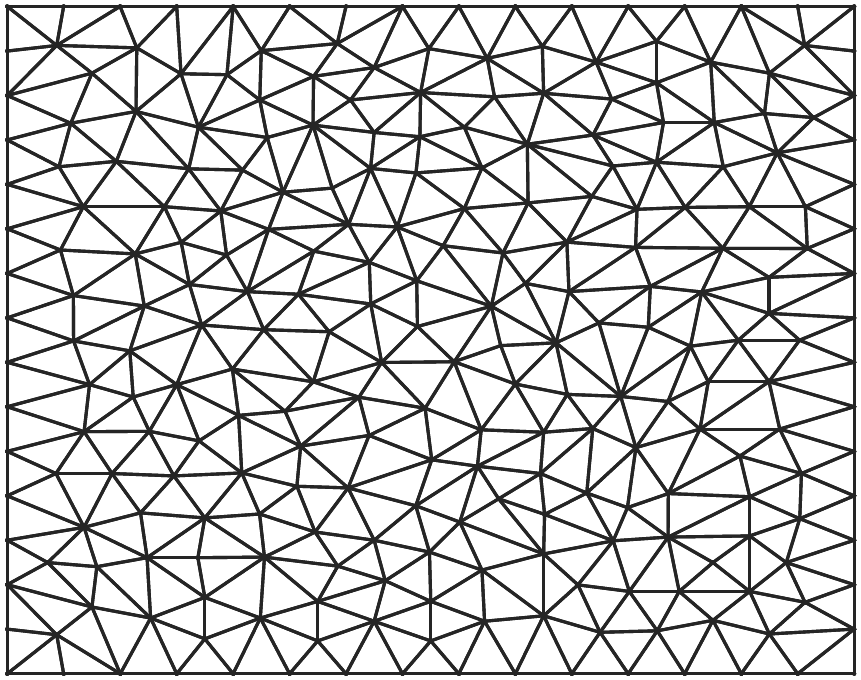}
\par\end{center}
\caption{\label{fig:Representative-example-2-D}Example unstructured 2D mesh
for convection-diffusion equation.}
\end{minipage}\hfill{} %
\begin{minipage}[t]{0.45\textwidth}%
\begin{center}
\includegraphics[width=0.98\columnwidth]{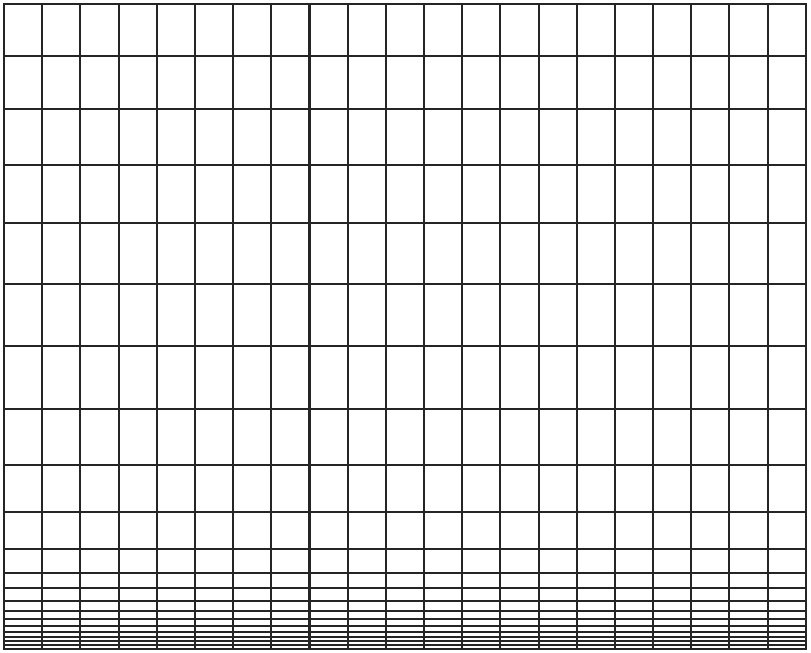}
\par\end{center}
\begin{flushleft}
\caption{\label{fig:nonuniform-mesh}Example nonuniform structured mesh for
Helmholtz equation.}
\par\end{flushleft}%
\end{minipage}
\end{figure}

Test cases E$d$CD{*} solve the convection-diffusion equation using
finite elements. Except for E$d$CDa, which we generated using our
in-house code, the other systems were generated using used FEniCS
v2017.2. E2CDa solves the problem on $\Omega=[0,1]^{2}$ with $\mu=1$,
$\vec{\nu}=\vec{1}$, $f=0$, and homogeneous Dirichlet BCs. The mesh
was generated using Triangle \cite{ShewchukTRIANGLE96}. Figure~\ref{fig:Representative-example-2-D}
shows the pattern of the mesh at a much coarser resolution. E3CDa1-3
are similar to E2CDa, except that the domain is $\Omega=[0,1]^{3}$,
triangulated using TetGen \cite{Si2006}. We generated three meshes
at different resolutions to facilitate the study of the asymptotic
complexity of preconditioned KSP methods as the number of unknowns
increases.

E2CDb is based on Example 6.1.2 in \cite{lman2014finite}. The domain
is $\Omega=[-1,1]^{2}$, with $\mu=1/200$ and $f=0$. The wind velocity
is $\vec{\nu}=\left[0,1+(x+1)^{2}/4\right]^{T}$, which is vertical
but increases in strength from left to right. Natural boundary conditions
are applied on the top wall. Dirichlet boundary conditions are imposed
elsewhere, with $u=1$ on the lower wall, and $u$ decreases cubically
and quadratically to zero on the left and right walls, respectively.
E2CDc is similar to E2CDb, except for a smaller diffusion coefficient
$\mu=1/2000$. E2CDd and E2CDe are based on Example 6.1.4 in \cite{lman2014finite},
known as the\textit{ recirculating wind problem or double-glazing
problem}, which models the temperature distribution in a cavity $\Omega=[-1,1]^{2}$
heated by an external wall. The source term is $f=0$, and the wind
velocity is $\vec{\nu}=\left[2y(1-x^{2}),-2x(1-y^{2})\right]$, which
determines a recirculating flow. Dirichlet boundary conditions are
imposed with $u=1$ on the right wall and $u=0$ elsewhere. For E2CDd
and E2CDe, $\mu$ is $1/200$ and $1/2000$, respectively. 

E3CDb-c solve the equation on $\Omega=[-1,1]^{3}$, with Dirichlet
boundary conditions 
\[
u(\vec{x})=\begin{cases}
1, & 0\leq x,y\leq0.5\text{ and }z=0,\\
0, & \text{otherwise}.
\end{cases}
\]
The source term is $f=0$, and the wind velocity is $\vec{\nu}=\left[1-z,1-z,1\right]$.
For E3CDb and E3CDc1-2, the dynamic viscosity is $\mu=10^{-12}$ and
$10^{-6}$, respectively. These problems are convection dominant,
and the FEM solutions may suffer from spurious oscillations on coarse
meshes. Hence, these problems may not be relevant physically but are
challenging algebraically. 

\paragraph{Convection-Diffusion Equation with GFDM and AES-FEM.}

GD$d$CD, for $d=2$ and $3$, were generated using our in-house implementation
of GFDM. Similarly, AE$d$CD were generated using our in-house implementation
of AES-FEM. The boundary conditions and parameters are the same as
E$d$CDa.

\paragraph{Convection-Diffusion Equation with DG.}

DG2CD solves the convection-diffusion equation using \emph{discontinuous
Galerkin} (\emph{DG}). We implemented it in FEniCS. The domain is
$\Omega=[0,1]^{2}$, with Dirichlet boundary conditions $u=\sin(5\pi y)$
on all sides. The dynamic viscosity is $\mu=2$, the wind velocity
is $\vec{\nu}=[x-0.5,0]$, and the source term is $f=3$.

\paragraph{Helmholtz Equation with FDM.}

D2HMa-c solve the Helmholtz equations using FDM on curvilinear meshes,
contributed by our collaborator Dr. Marat Khairoutdinov, who is a
climate scientist, and our colleague Oliver Yang. These equations
arise from a 3D Poisson equation for pressure $p$ in the longitude-latitude
coordinate system in a global climate model, 
\begin{equation}
\frac{\partial^{2}p}{\partial x^{2}}+\mu(y)\frac{\partial}{\partial y}\left(\mu(y)\frac{\partial p}{\partial y}\right)+\frac{\mu(y)^{2}}{\rho(z)}\frac{\partial}{\partial z}\left(\rho(z)\frac{\partial p}{\partial z}\right)=\mu(y)^{2}f,\label{eq:Poisson-longitute-latitude}
\end{equation}
with natural boundary conditions. By applying Fourier transform along
in $x$-direction to the 3D system, we obtain a collection of 2D Helmholtz
equations in the frequency domain,
\begin{equation}
\mu\left(y\right)\frac{\partial}{\partial y}\left[\mu\left(y\right)\frac{\partial p}{\partial y}\right]+\frac{\mu\left(y\right){}^{2}}{\rho\left(z\right)}\frac{\partial}{\partial z}\left[\rho\left(z\right)\frac{\partial p}{\partial z}\right]-\omega^{2}p=\tilde{f},\label{eq:FDM-Hemlholtz}
\end{equation}
where $\omega$ corresponds to the frequency. This equation is then
discretized using conservative, cell-centered finite differences.
The mesh in the $z$-direction is highly anisotropic, as illustrated
in Figure~\ref{fig:nonuniform-mesh} at a much coarser resolution.
The frequency $\omega$ is zero for one of the $yz$-planes, and in
this case the system is singular, which is a challenging problem in
its own right and is beyond the scope of this study. For the adjacent
planes, $\omega$ can be arbitrarily close to zero, depending on the
grid resolution along $x$. D2HMa-c correspond to $\omega^{2}=10^{-13}$,
$10^{-10}$, and $10^{-6}$, respectively. The strong anisotropy makes
these Helmholtz equations difficult, especially for small $\omega$. 

\paragraph{Incompressible Navier-Stokes Equations with FEM.}

E$d$INS{*} solve the incompressible Navier-Stokes equations, which
read
\begin{align}
\rho(\dot{\vec{u}}+\vec{u}\cdot\vec{\nabla}\vec{u})-\vec{\nabla}\cdot\vec{\sigma}(\vec{u},p)= & \vec{f},\label{eq:ns-momentum}\\
\vec{\nabla}\cdot\vec{u}= & 0,\label{eq:ns-continuity}
\end{align}
where $\rho$ denotes density, $\vec{u}$ the fluid velocity, $\dot{\vec{u}}$
the acceleration, $\vec{\sigma}$ the stress tensor, and $\vec{f}$
an applied body force. For a Newtonian fluid,
\begin{equation}
\vec{\sigma}(\vec{u},p)=2\mu\vec{\epsilon}(\vec{u})-p\vec{I},\label{eq:stress-tensor-Newtonian}
\end{equation}
where $\mu$ is the dynamic viscosity and $\vec{\epsilon}(\vec{u})$
is the strain-rate-tensor 
\begin{equation}
\vec{\epsilon}(\vec{u})=\frac{1}{2}\left(\vec{\nabla}\vec{u}+\left(\vec{\nabla}\vec{u}\right)^{T}\right).\label{eq:strain-rate-tensor}
\end{equation}
Eqs.~(\ref{eq:ns-momentum}) and (\ref{eq:ns-continuity}) are referred
to as \emph{momentum} and \emph{continuity} equations, respectively.
For the finite element methods, these equations are typically solved
using Chroin's method (see e.g. \cite{rannacher1992chorin}), a.k.a.
the \emph{projection method} \cite{BroCorMin01}, which first solves
the momentum equation and then solves Poisson equations to recover
the divergence-free property. We discretize these equations using
finite elements with FEniCS. We use the linear system from the first
step, which is nonsymmetric. We obtain the right-hand side vectors
after running the simulation for a few time steps, so that the flows
are better developed and the vectors are more representative.

E2INSa1 and E2INSa2 correspond to the channel-flow problem over $\Omega=[0,1]^{2}$,\textcolor{red}{{}
}with $\rho=1$ and $\mu=1$. We use Dirichlet boundary conditions
$p=8$ and $p=0$ on the left and right walls, and $\vec{u}=\vec{0}$
on the top and bottom walls. E3NS is a similar problem in 3D, with
natural boundary conditions along the $z$-direction. E2INSb and E2INSc
correspond to a flow around a cylinder, which is a classical benchmark
problem \cite{schafer1996benchmark}. For E2INSb, $\nu=0.001=\mu/\rho$,
with a corresponding to Reynolds number of $100$; for E2INSc, $\nu=0.0005$,
with a corresponding Reynolds number of $200$. The remainder of the
parameters were the same as those in \cite{schafer1996benchmark}.
E2INSb-c are very well-conditioned, with condition numbers less than
$10^{3}$, so they represent the easiest problems in this benchmark
set.

\paragraph{Compressible Navier-Stokes Equations with FVM.}

V3CNS solves the compressible Navier-Stokes equations, which is a
system of PDEs based on the conservation of mass, momentum, and energy.
 We use the matrix RM07R in Group Fluorem from the UF Sparse Matrix
Collection, which solves the equation using the \emph{finite-volume
methods} (\emph{FVM}); see \cite{pacull2011study} for more detail
of the case. This problem is the most ill-conditioned in this benchmark
set.

\paragraph{Stokes Equations with Mixed FEM.}

E3STHa1 and E3STHa2 solve the Stokes equations, which describe a steady,
incompressible Newtonian flow with low Reynolds numbers. For a domain
$\Omega\subset\mathbb{R}^{d}$, the Stokes equations read
\begin{align}
-\vec{\nabla}\cdot(\mu\vec{\nabla}\vec{u})+\vec{\nabla}p & =\vec{f},\label{eq:stokes_model_problem}\\
\vec{\nabla}\cdot\vec{u} & =0,
\end{align}
where $\vec{u}$ denotes the fluid velocity, $\mu$ the dynamic viscosity,
$p$ the pressure, and $\vec{f}$ an applied body force. This problem
is unstable with linear elements. We solve it using mixed Taylor-Hood
elements \cite{ern2013theory}, which use quadratic elements for $\vec{u}$
and linear elements for $p$. The domain is a unit cube $\Omega=[0,1]^{3}$,
with $\mu=2$ and $\vec{f}=\vec{0}$. The boundary condition is $\vec{u}=\left[-\sin(\pi y),0,0\right]$
for $x=1$ and $\vec{u}=\vec{0}$ everywhere else. The resulting linear
system may be either symmetric and indefinite or nonsymmetric and
positive definite; we used the nonsymmetric variant. A notable property
of the linear system is that there is a large zero diagonal block,
similar to saddle-point problems. We refer to these systems as \emph{saddle-point-like
problems}. 

\paragraph{Poisson Equation with Mixed FEM.}

E3MP is another example of saddle-point-like problems, from the mixed
formulation of the Poisson equation $-\Delta u=f$. This formulation
introduces a vector field, namely $\vec{v}=\vec{\nabla}u$, and then
constructs a system of first-order PDEs
\begin{align}
\vec{v}-\vec{\nabla}u & =\vec{0}\qquad\,\,\,\,\text{ in }\Omega,\label{eq:mixed-Poisson-1}\\
\vec{\nabla}\cdot\vec{v} & =-f\qquad\text{ in }\Omega.\label{eq:mixed-Poisson-2}
\end{align}
See e.g. \cite{boffi2013mixed} for details on mixed FEM. These equations
are similar to those arising from the \emph{Darcy flow} in porous
media. We solve the equations using FEniCS with linear Brezzi-Douglas-Marini
elements for $u$ and degree-$0$ discontinuous elements for $\vec{v}$.
The domain is $\Omega=[0,1]^{3}$, with homogeneous Dirichlet boundary
conditions along the left and right walls, and Neumann boundary conditions
$\vec{v}\cdot\vec{n}=\sin(5x)$ elsewhere. The source term was $f=10\exp(-((x-0.5)^{2}+(y-0.5)^{2}+(z-0.5)^{2})/0.02)$. 

\section{Numerical Comparisons\label{sec:Results}}

In this section, we present numerical comparisons of the preconditioned
Krylov methods described in Sections~\ref{sec:Analysis-KSP}. For
GMRES, TFQMR, and BiCGSTAB, we use the built-in implementations in
PETSc v3.7.1. For GMRES, we use 30 as the restart parameter, which
is the recommended value for large-scale linear systems in PETSc,
and we denote the method by GMRES(30). For QMRCGSTAB, which is unavailable
in PETSc, we use our implementations based on the low-level functions
in PETSc. In terms of preconditioners, we consider Gauss-Seidel (GS),
ILU0, ILUTP (SuperLU v5.2.1), MILU (ILUPACK v2.4), variants of classical
AMG (BoomerAMG\textcolor{blue}{{} }in hypre v2.11.0), and smoothed-aggregation
AMG (Trilinos/ML 5.0), all as right preconditioners. For a uniform
comparison, we set the relative convergence tolerance for the $2$-norm
of the residual, i.e., the 2-norm of the residual divided by the 2-norm
of the right-hand side, to $10^{-10}$ for all methods.

As an overview, Table~\ref{tab:results_outline} shows the runtimes
of all the benchmark problems described in Section~\ref{sec:benchmark-problems}
with GMRES(30), which is the most robust. All the tests were conducted
in serial on a single node of a cluster\textcolor{blue}{{} }with two
2.6 GHz Intel Xeon E5-2690v3 processors and 128 GB of memory. For
accurate timing, both turbo and power saving modes were turned off
for the processors, and each node was dedicated to run one problem
at a time. In the table, the best timing results were in bold, `-'
indicates stagnation, `{*}' indicates ``did not converge'' after
10,000 iterations, and `$\times$' indicates runtime errors in SuperLU
(in particular, nontermination of factorization after 24 hours for
E3STHa2 and malloc errors for the others). For ILUTP, we used SuperLU
with the default drop tolerance of $10^{-4}$. For hypre/BoomerAMG,
we used HMIS coarsening with FF1 interpolation for all cases except
for E2CDc, which used PMIS+FF1 due to non-convergence with HMIS+FF1.
The other parameters are all default. It is clear that hypre is the
overall winner for convection-diffusion and Helmholtz equations. However,
Gauss-Seidel is the most efficient for well-conditioned systems, such
as those from incompressible Navier-Stokes equations. We will not
consider the well-conditioned systems further. On the other end, multilevel
ILU tends to be the most robust for ill-conditioned and saddle-point-like
problems. In the following subsections, we give more detailed comparisons
in terms of convergence of different KSP methods, ILU, and multigrid
preconditioners.

\begin{table}[h]
\caption{\label{tab:results_outline}Overview of runtimes (in seconds) of GMRES(30)
with Gauss-Seidel, ILU0, ILUTP, MILU, hypre, and ML. `$-$', `{*}',
and `$\times$' indicate stagnation, ``did not converge,'' and runtime
error, respectively.}

\centering{}%
\begin{tabular}{>{\centering}p{2.7cm}|c|>{\centering}p{1.5cm}|>{\centering}p{1.5cm}|>{\centering}p{1.5cm}|>{\centering}p{1.7cm}|>{\centering}p{1.6cm}}
\hline 
\textbf{Matrix ID} & \textbf{GS} & \textbf{ILU0} & \textbf{ILUTP} & \textbf{MILU} & \textbf{HYPRE} & \textbf{ML}\tabularnewline
\hline 
\hline 
E2CDa & 1,517.3 & 1,390.5 & 85.3 & 143.0 & \textbf{7.52} & 17.88\tabularnewline
\hline 
E2CDb & 5,554.2 & 2,527.4 & 917.9 & 1,173.6 & \textbf{38.14} & 246.7\tabularnewline
\hline 
E2CDc & 3,355.2 & 1,377.8 & 574.2 & 1,107.2 & \textbf{71.1{*}} & 236.3\tabularnewline
\hline 
E2CDd & {*} & {*} & 881.4 & 1,309.6 & \textbf{44.41} & 747.8\tabularnewline
\hline 
E2CDe & {*} & {*} & 2,287.1 & 1,291.5 & \textbf{77.10} & 2,936.1\tabularnewline
\hline 
E3CDa1 & 20.41 & 12.63 & 1,158.4 & 17.35 & \textbf{5.75} & 8.48\tabularnewline
\hline 
E3CDa2 & 263.8 & 157.9 & 32,389.5 & 179.3 & \textbf{51.40} & 70.05\tabularnewline
\hline 
E3CDa3 & 7,191.2 & 5,126.1 & $\times$ & 2,430.1 & \textbf{572.93} & 853.4\tabularnewline
\hline 
E3CDb & 366.4 & 213.6 & 55,372.5 & 1,314.3 & \textbf{52.28} & 84.17\tabularnewline
\hline 
E3CDc1 & 292.13 & 172.04 & 54,793.1 & 1,313.1 & \textbf{55.0} & 77.05\tabularnewline
\hline 
E3CDc2 & 2,744.4 & 1,548.7 & $\times$ & 2,125.7 & \textbf{232.9} & 360.60\tabularnewline
\hline 
AE2CD & 416.5 & 280.5 & 128.9 & 59.2 & \textbf{18.2} & 25.3\tabularnewline
\hline 
AE3CD & 6140.3 & 5018.5 & $\times$ & 2,440.9 & \textbf{534.3} & 819.1\tabularnewline
\hline 
GD2CD & 1,543.1 & 1,230.5 & 87.01 & 142.65 & \textbf{8.45} & 23.49\tabularnewline
\hline 
GD3CD & 178.40 & 113.7 & 42,819.9 & 168.1 & \textbf{50.79} & 69.66\tabularnewline
\hline 
DG2CD & {*} & {*} & 180.79 & 145.4 & \textbf{10.60} & 21.49\tabularnewline
\hline 
D2HMa & {*} & {*} & 223.5 & 131.2 & \textbf{5.42} & 211.7\tabularnewline
\hline 
D2HMb & {*} & {*} & 30.67 & 78.1 & \textbf{5.96} & 35.93\tabularnewline
\hline 
D2HMc & 1.11 & \textbf{0.91} & 18.36 & 9.6 & \textbf{0.91} & 7.40\tabularnewline
\hline 
E2INSa1 & 846.2 & 778.9 & 273.34 & \textbf{171.2} & 938.8 & 1020.7\tabularnewline
\hline 
E2INSa2 & 1,336.2 & 1,048.9 & 192.79 & \textbf{217.5} & 1,175.6 & 1861.3\tabularnewline
\hline 
E2INSb & \textbf{58.8} & 78.16 & $\times$ & 672.7 & 148.7 & 165.8\tabularnewline
\hline 
E2INSc & \textbf{107.2} & 129.3 & $\times$ & 426.8 & 339.1 & 232.9\tabularnewline
\hline 
E3INS & \textbf{383.9} & 392.1 & $\times$ & 617.4 & 785.5 & 732.1\tabularnewline
\hline 
V3CNS & {*} & {*} & {*} & \textbf{21,039.6} & {*} & {*}\tabularnewline
\hline 
\multicolumn{7}{c}{\emph{Saddle-point-like problems}}\tabularnewline
\hline 
E3STHa1 & $-$ & $-$ & 86.7 & \textbf{7.74} & $-$ & $-$\tabularnewline
\hline 
E3STHa2 & $-$ & $-$ & $\times$ & \textbf{375.5} & $-$ & $-$\tabularnewline
\hline 
E3MPa1 & $-$ & $-$ & 4,033.2 & \textbf{129.4} & $-$ & $-$\tabularnewline
\hline 
E3MPa2 & $-$ & $-$ & 43,385.4 & \textbf{352.6} & $-$ & $-$\tabularnewline
\hline 
\end{tabular}
\end{table}

\subsection{Comparison of Krylov-Subspace Methods}

We first compare the four different Krylov-subspace methods. We consider
three preconditioners: GS (available as SOR with relaxation parameter
$\omega=1$ in PETSc), ILU0 (the default serial preconditioner in
PETSc), and BoomerAMG with HMIS coarsening and FF1 interpolation (which
differs from the default option in hypre). In Sections~\ref{subsec:ILU-Comparison}
and \ref{subsec:HYPRE-vs-ML}, we will compare the different variants
of ILU and AMG, respectively.

\subsubsection{\label{subsec:Convergence-Comparison.}Convergence Comparison.}

\begin{figure}[h]
\begin{minipage}[t]{0.45\textwidth}%
\begin{center}
\includegraphics[width=1\textwidth]{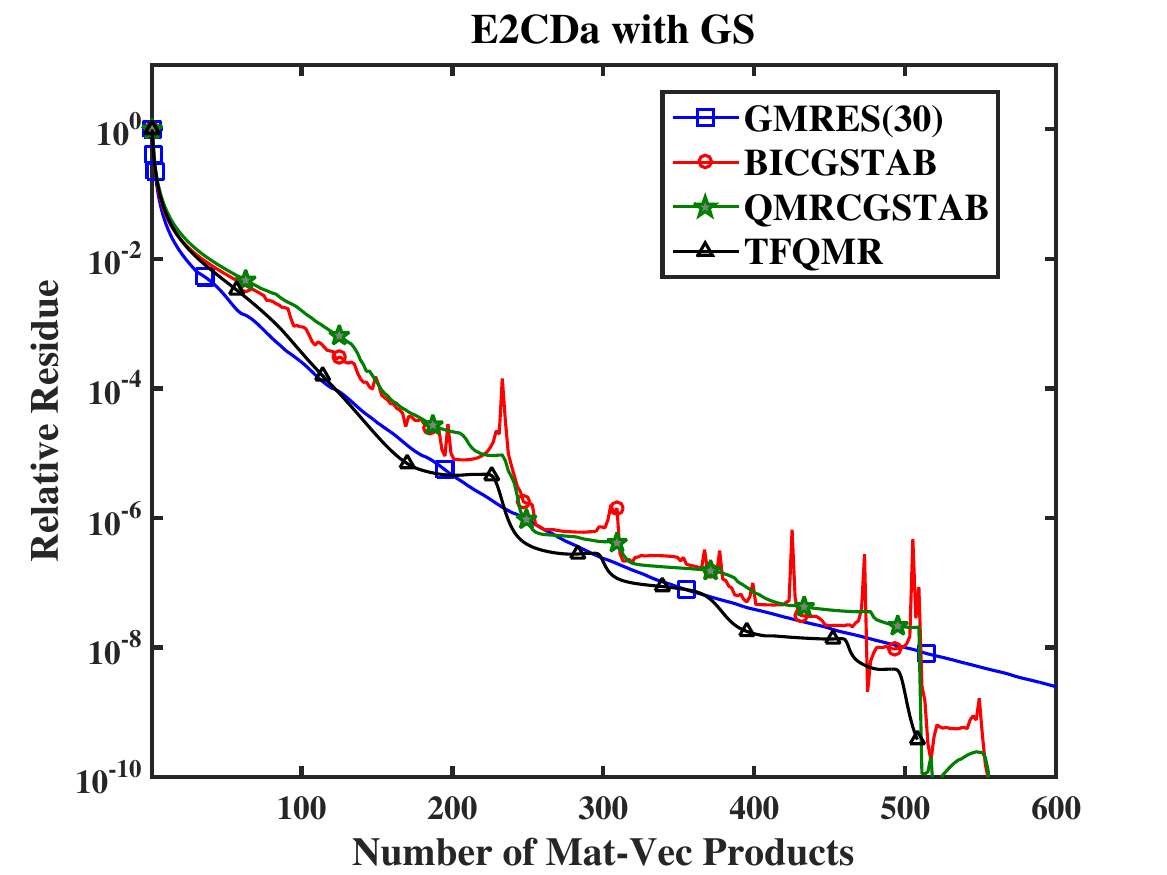}
\par\end{center}%
\end{minipage}\hfill{} %
\begin{minipage}[t]{0.45\textwidth}%
\begin{center}
\includegraphics[width=1\textwidth]{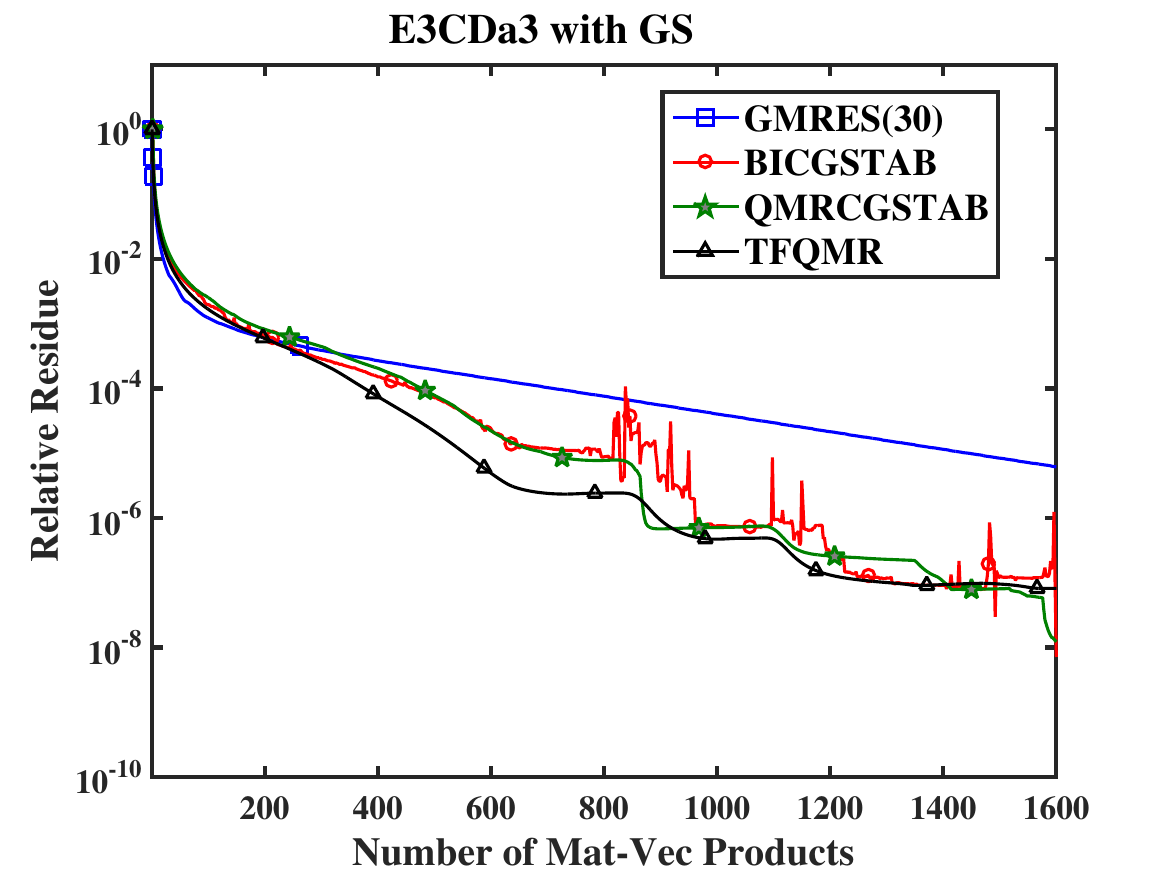}
\par\end{center}%
\end{minipage}

\begin{minipage}[t]{0.45\textwidth}%
\begin{center}
\includegraphics[width=1\textwidth]{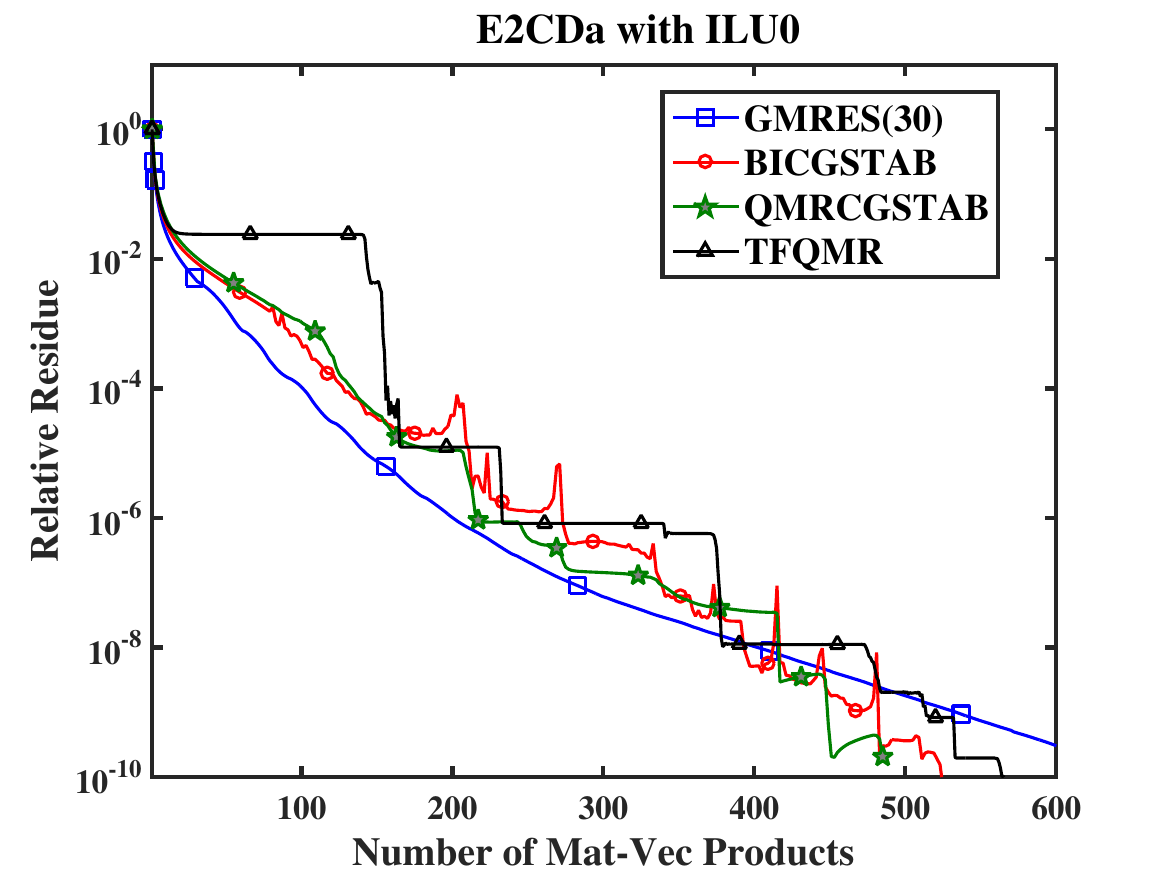}
\par\end{center}%
\end{minipage}\hfill{} %
\begin{minipage}[t]{0.45\textwidth}%
\begin{center}
\includegraphics[width=1\textwidth]{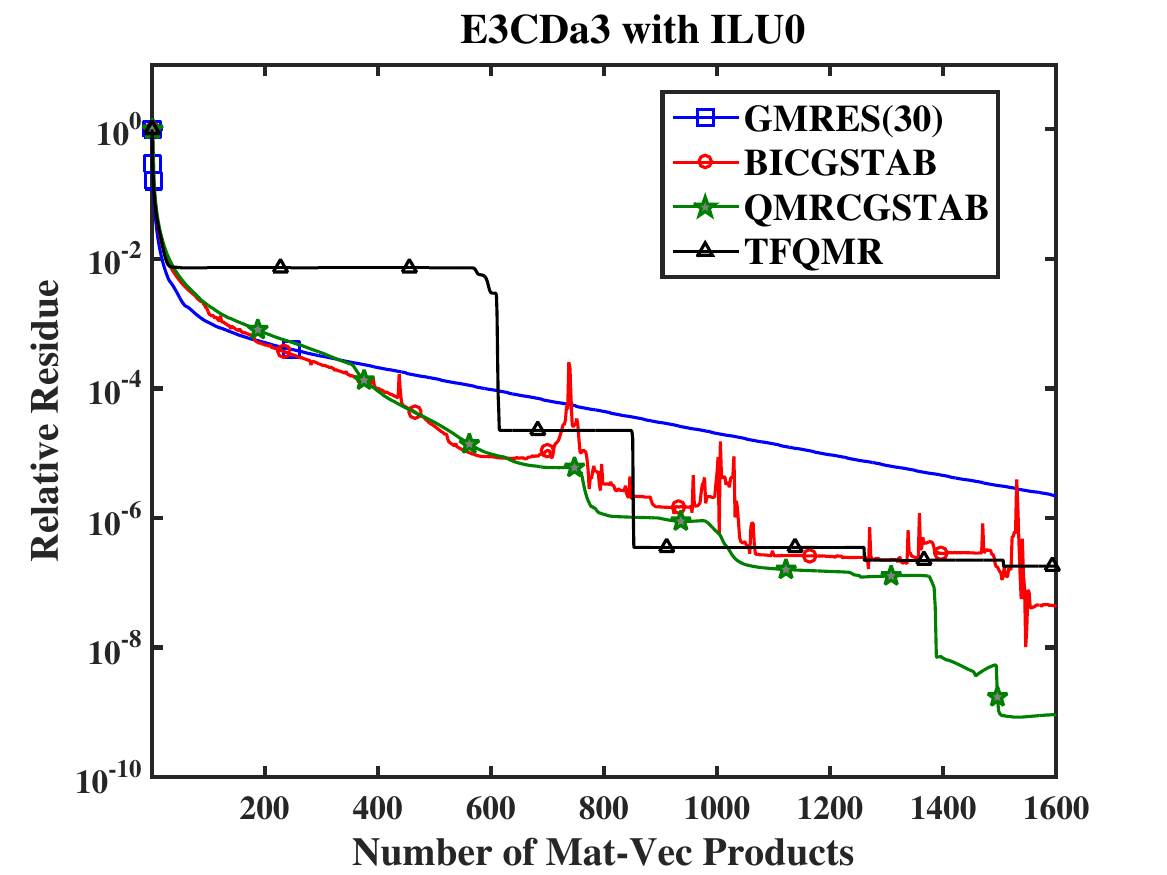}
\par\end{center}%
\end{minipage}

\begin{minipage}[t]{0.45\textwidth}%
\begin{center}
\includegraphics[width=1\textwidth]{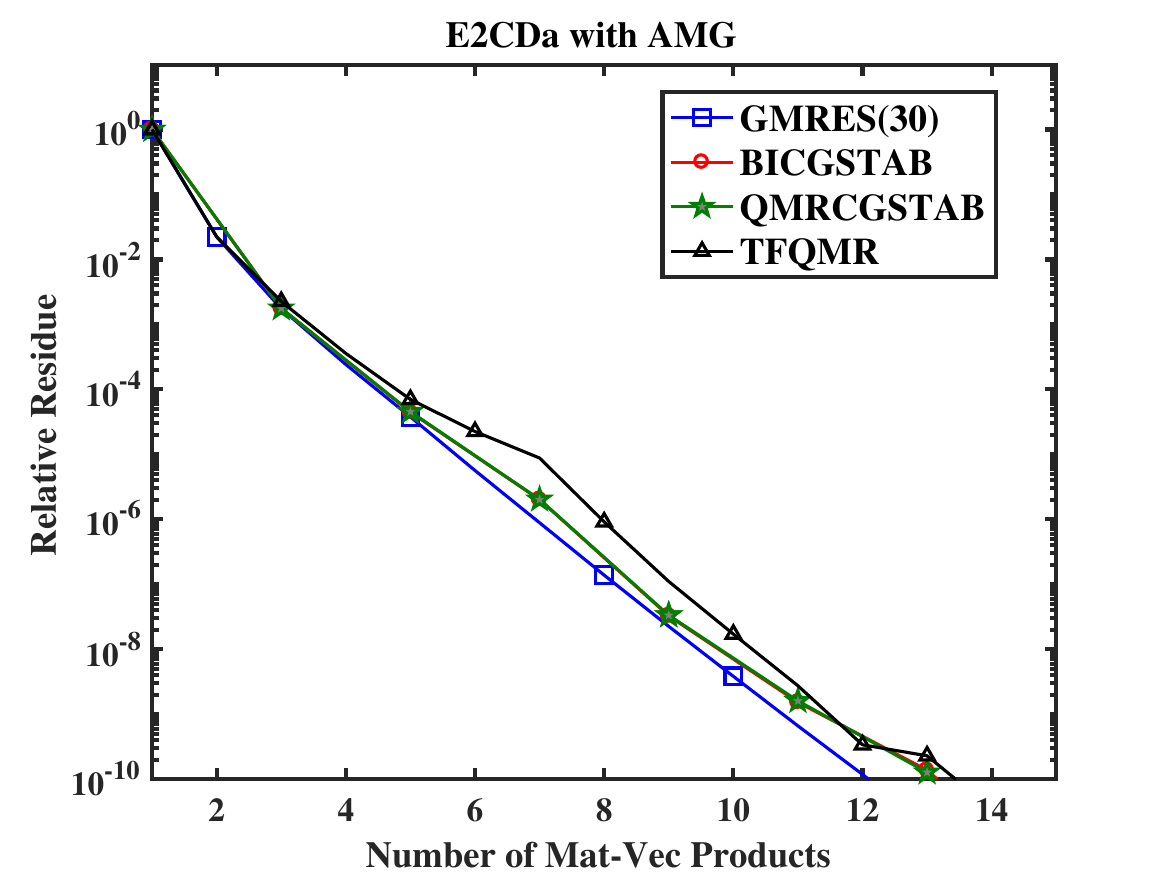}
\par\end{center}%
\end{minipage}\hfill{} %
\begin{minipage}[t]{0.45\textwidth}%
\begin{center}
\includegraphics[width=1\textwidth]{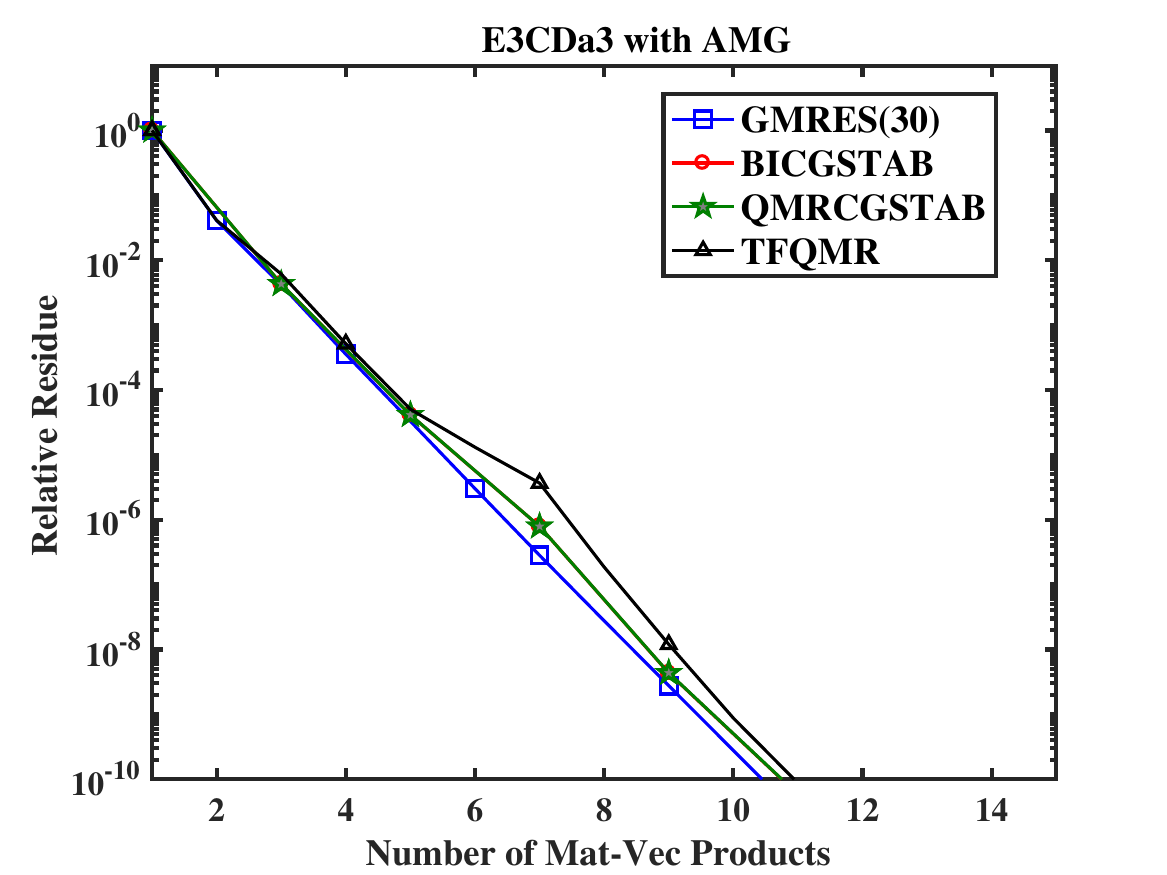}
\par\end{center}%
\end{minipage}\caption{\label{fig:FEM-residual}Residuals vs. numbers of matrix-vector products
for E2CDa (left) and E3CDa3 (right).}
\end{figure}

\begin{figure}[h]
\begin{minipage}[t]{0.45\textwidth}%
\begin{center}
\includegraphics[width=1\textwidth]{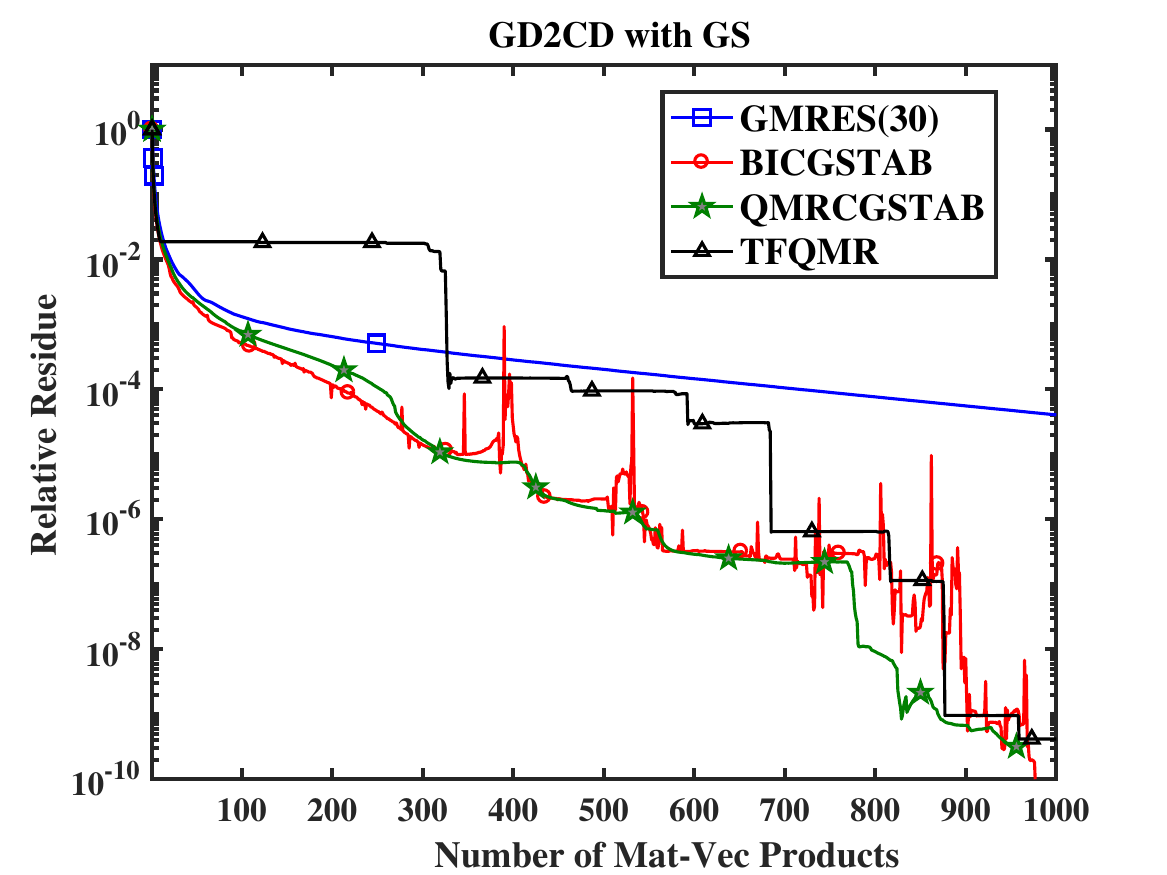}
\par\end{center}%
\end{minipage}\hfill{} %
\begin{minipage}[t]{0.45\textwidth}%
\begin{center}
\includegraphics[width=1\textwidth]{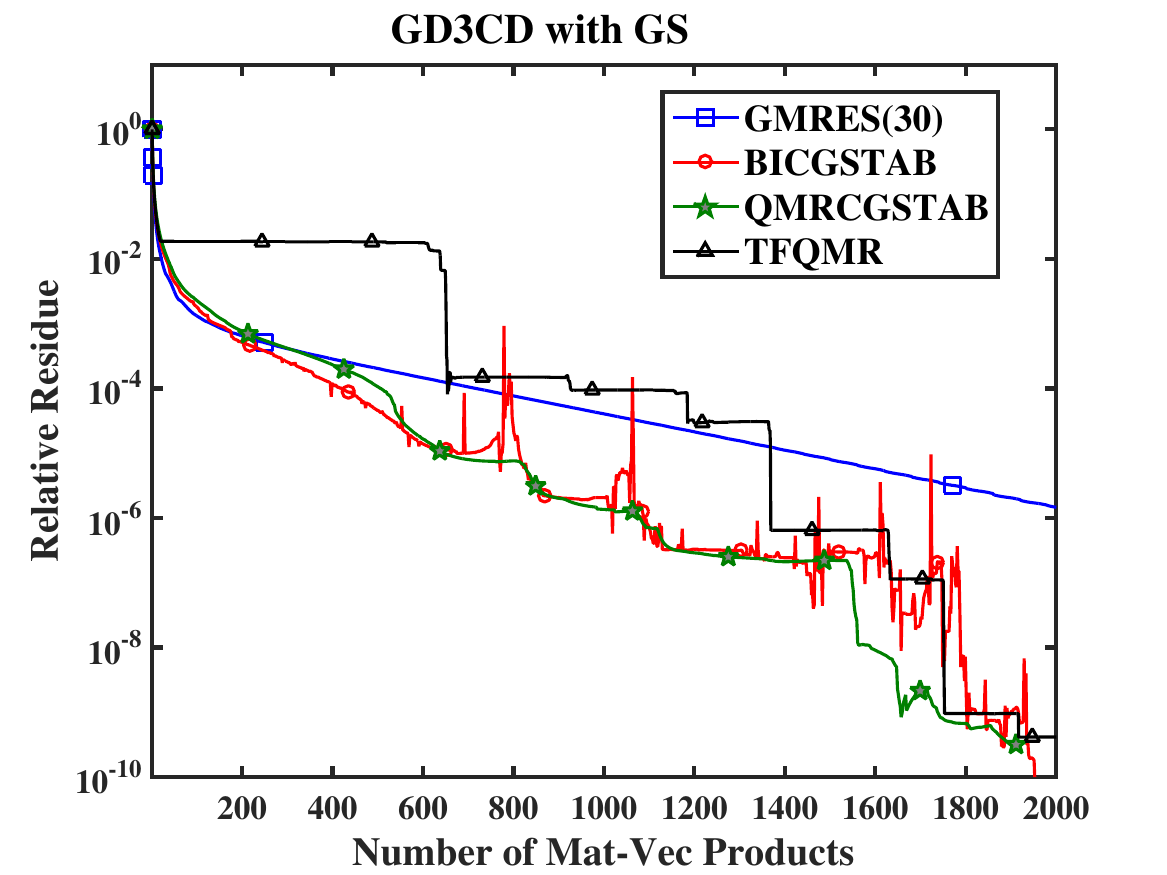}
\par\end{center}%
\end{minipage}

\begin{minipage}[t]{0.45\textwidth}%
\begin{center}
\includegraphics[width=1\textwidth]{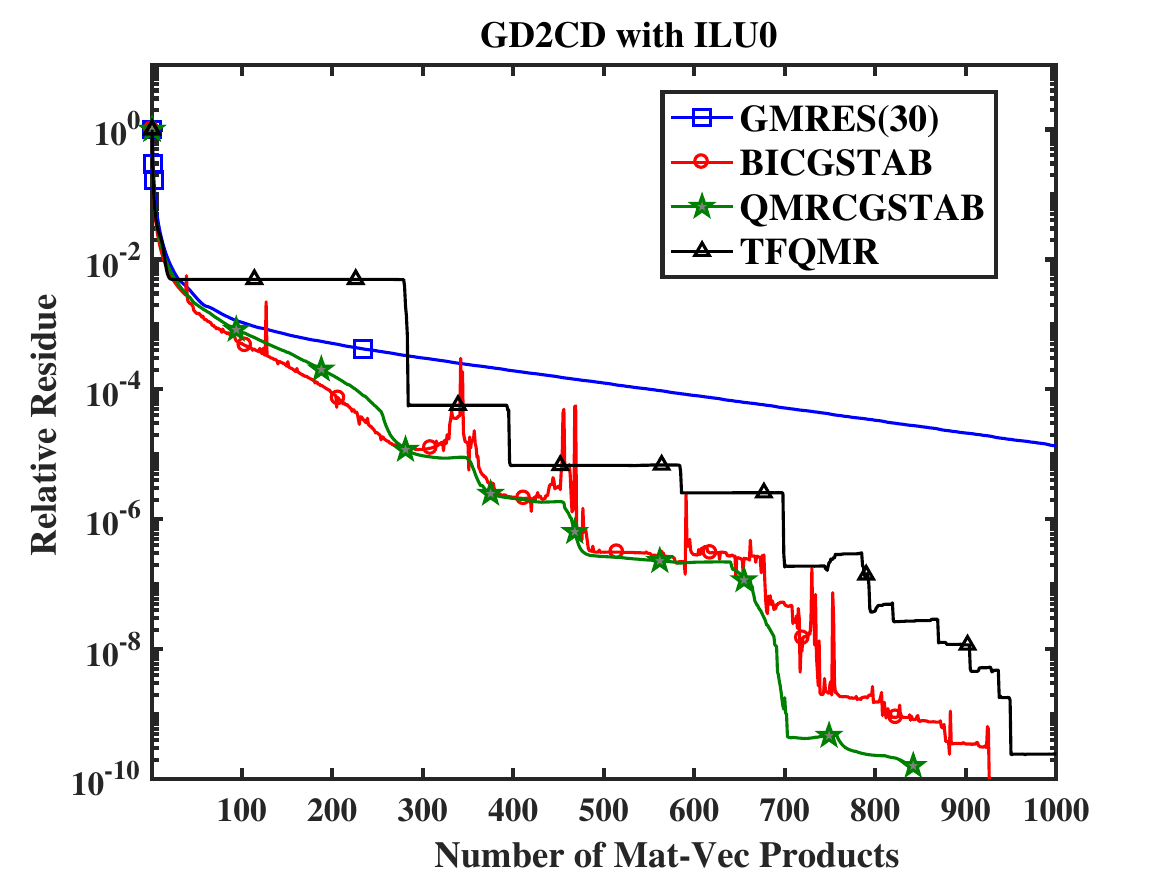}
\par\end{center}%
\end{minipage}\hfill{} %
\begin{minipage}[t]{0.45\textwidth}%
\begin{center}
\includegraphics[width=1\textwidth]{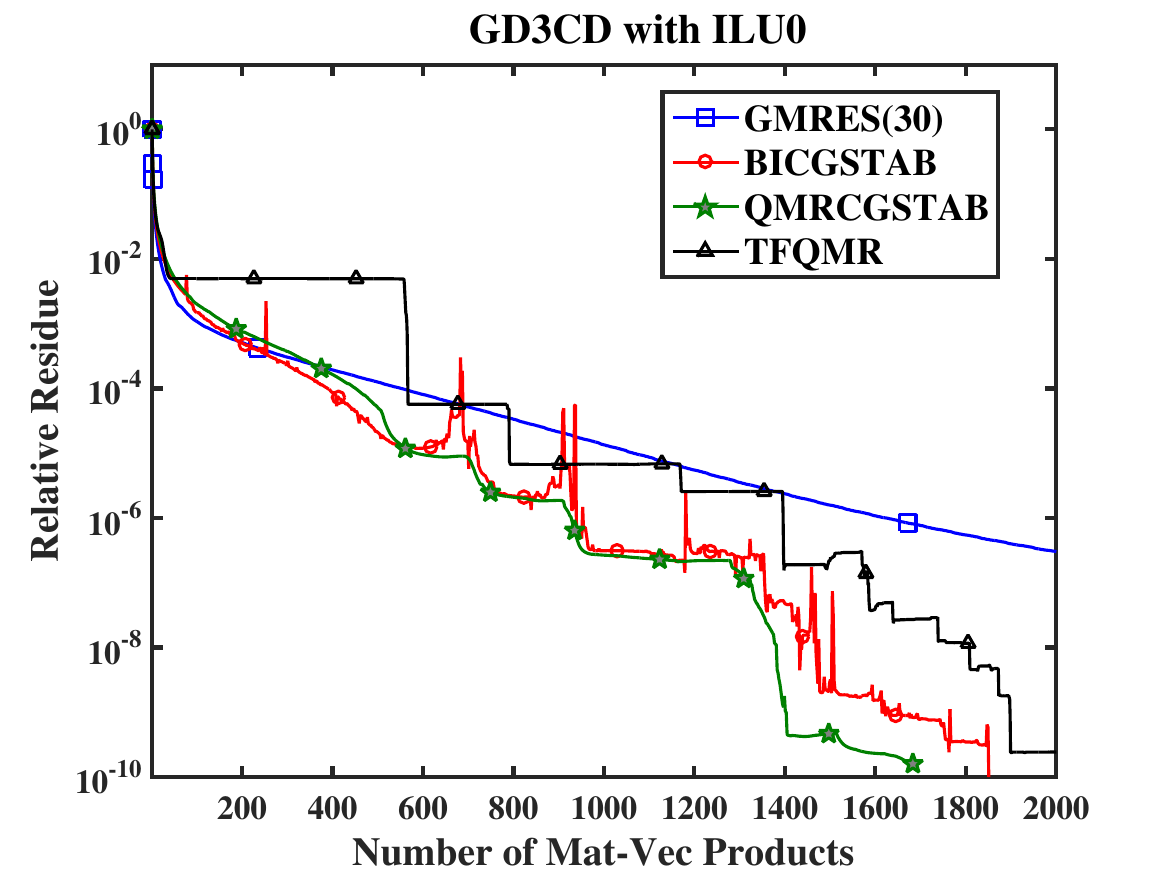}
\par\end{center}%
\end{minipage}

\begin{minipage}[t]{0.45\textwidth}%
\begin{center}
\includegraphics[width=1\textwidth]{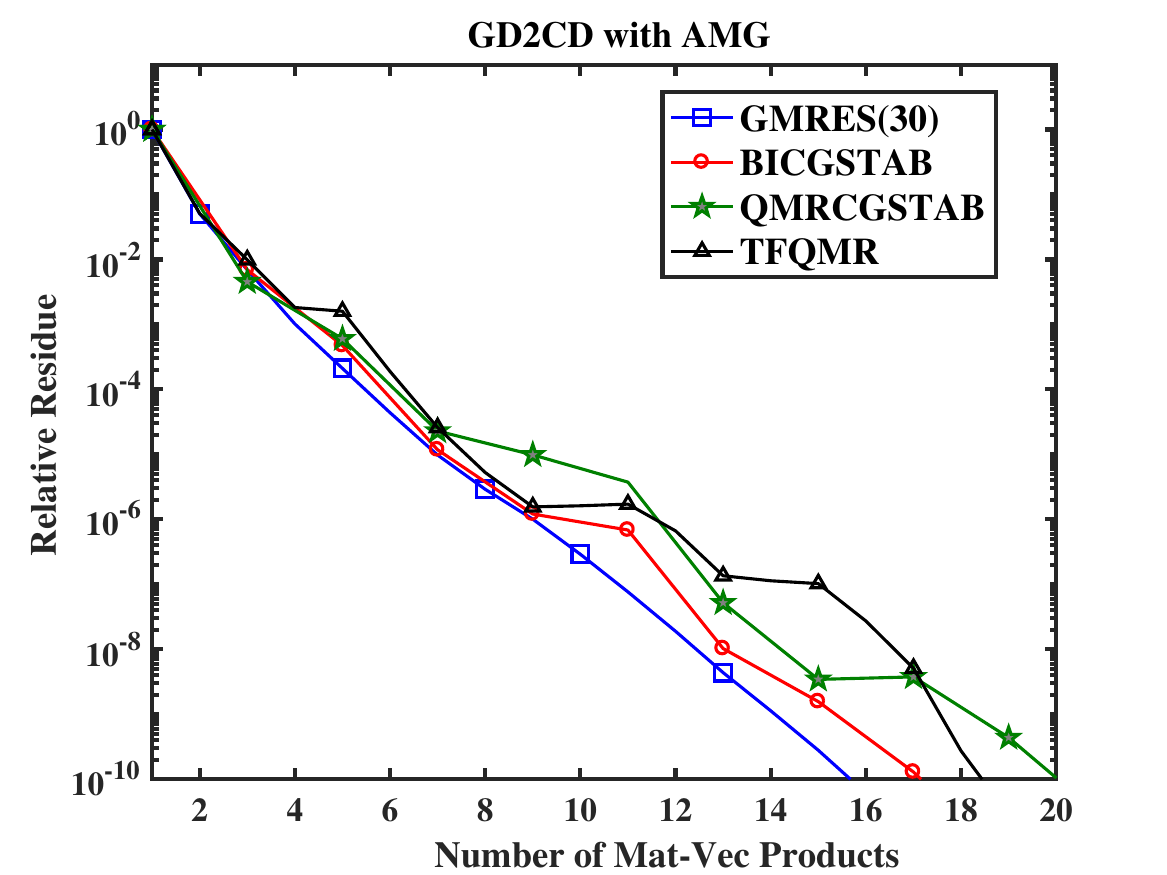}
\par\end{center}%
\end{minipage}\hfill{} %
\begin{minipage}[t]{0.45\textwidth}%
\begin{center}
\includegraphics[width=1\textwidth]{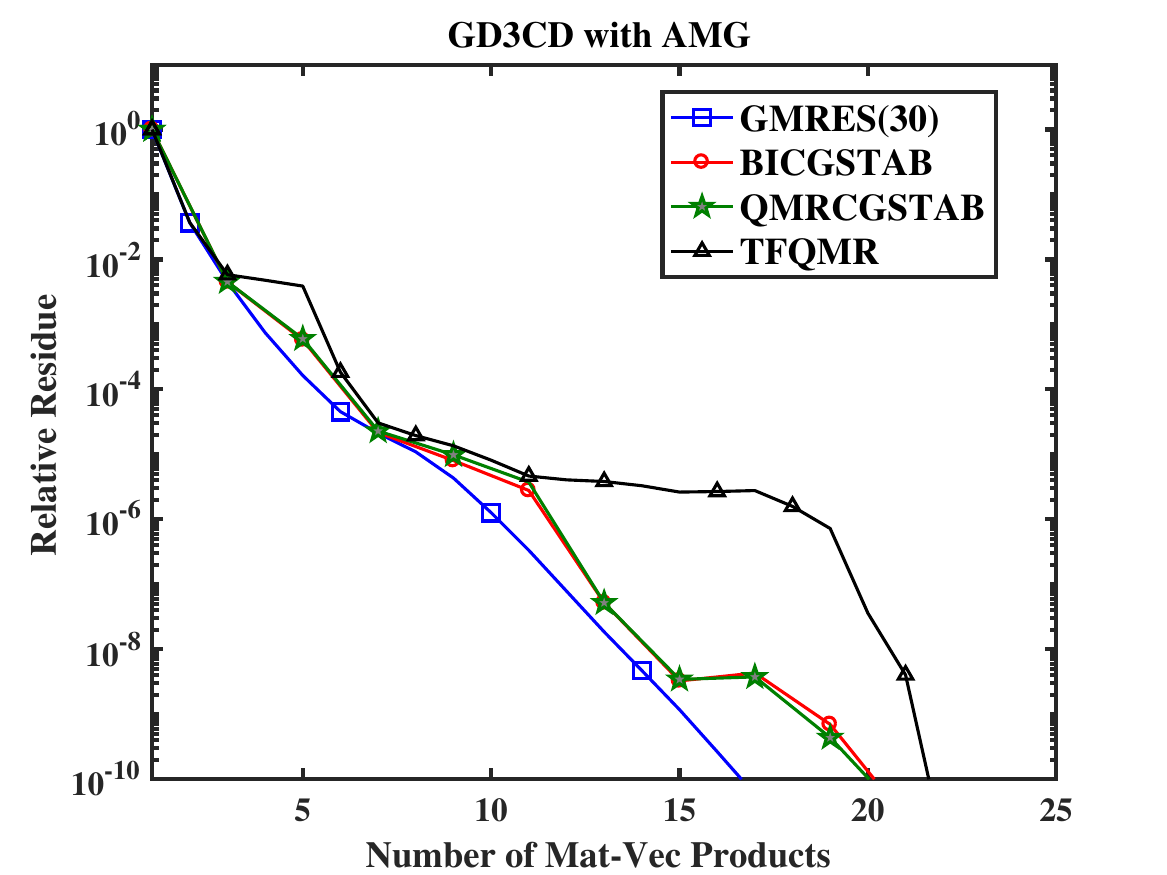}
\par\end{center}%
\end{minipage}

\caption{\label{fig:GFD-residual}Residuals vs. numbers of matrix-vector products
for GD2CD (left) and GD3CD (right).}
\end{figure}

\begin{figure}[h]
\begin{minipage}[t]{0.45\textwidth}%
\begin{center}
\includegraphics[width=1\textwidth]{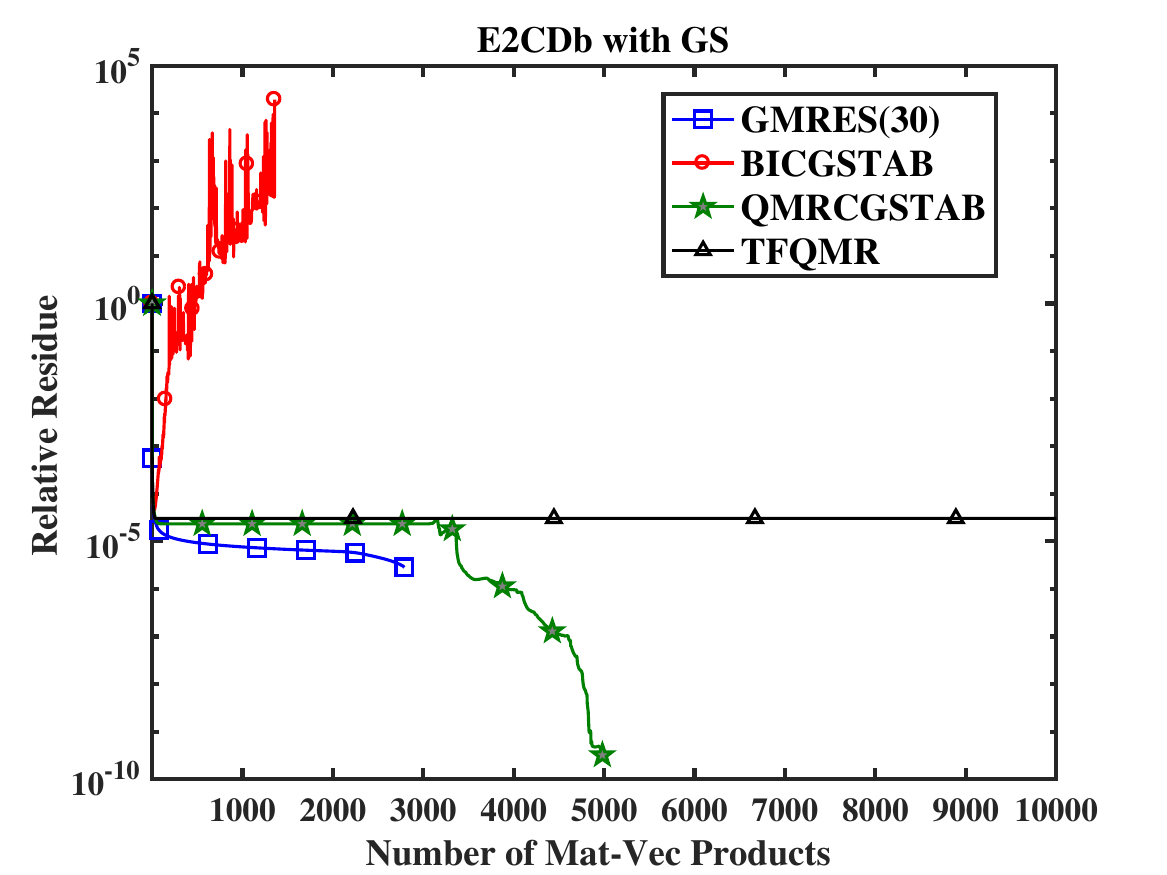}
\par\end{center}%
\end{minipage}\hfill{} %
\begin{minipage}[t]{0.45\textwidth}%
\begin{center}
\includegraphics[width=1\textwidth]{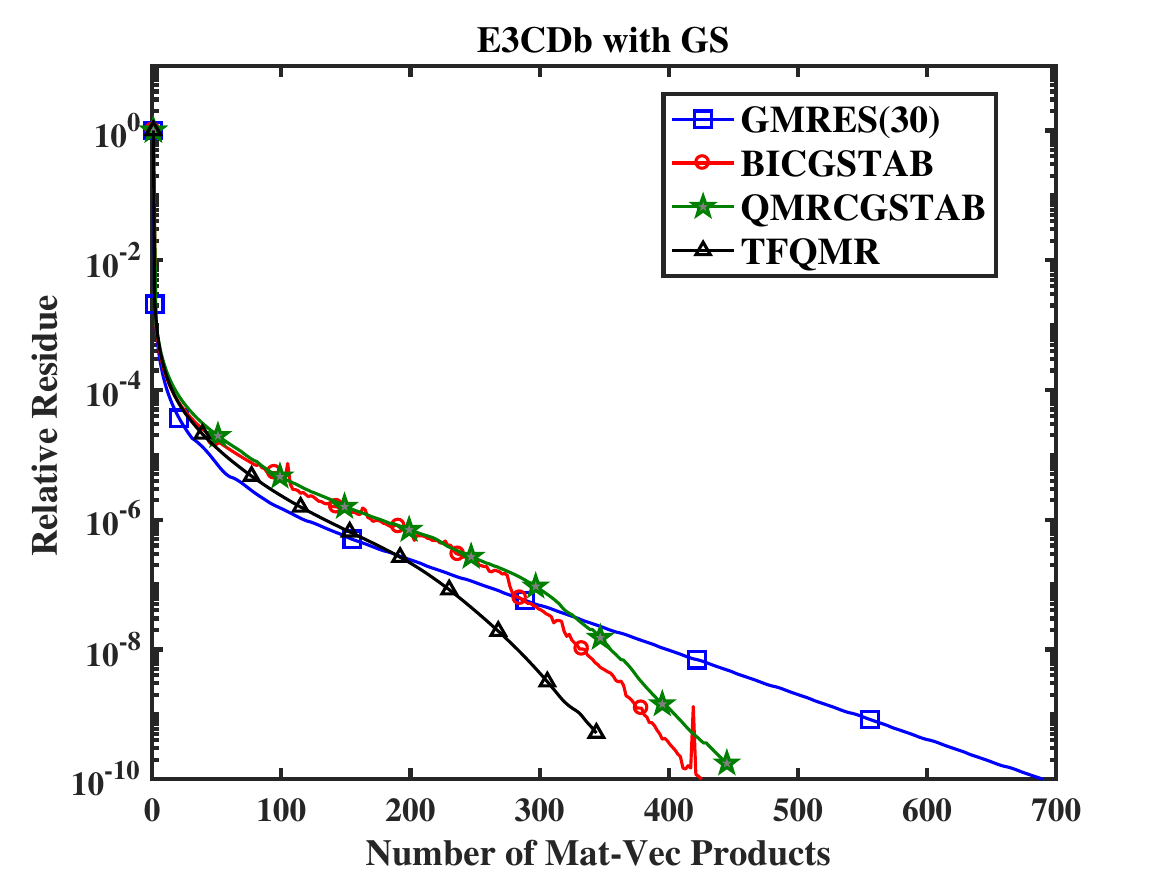}
\par\end{center}%
\end{minipage}

\begin{minipage}[t]{0.45\textwidth}%
\begin{center}
\includegraphics[width=1\textwidth]{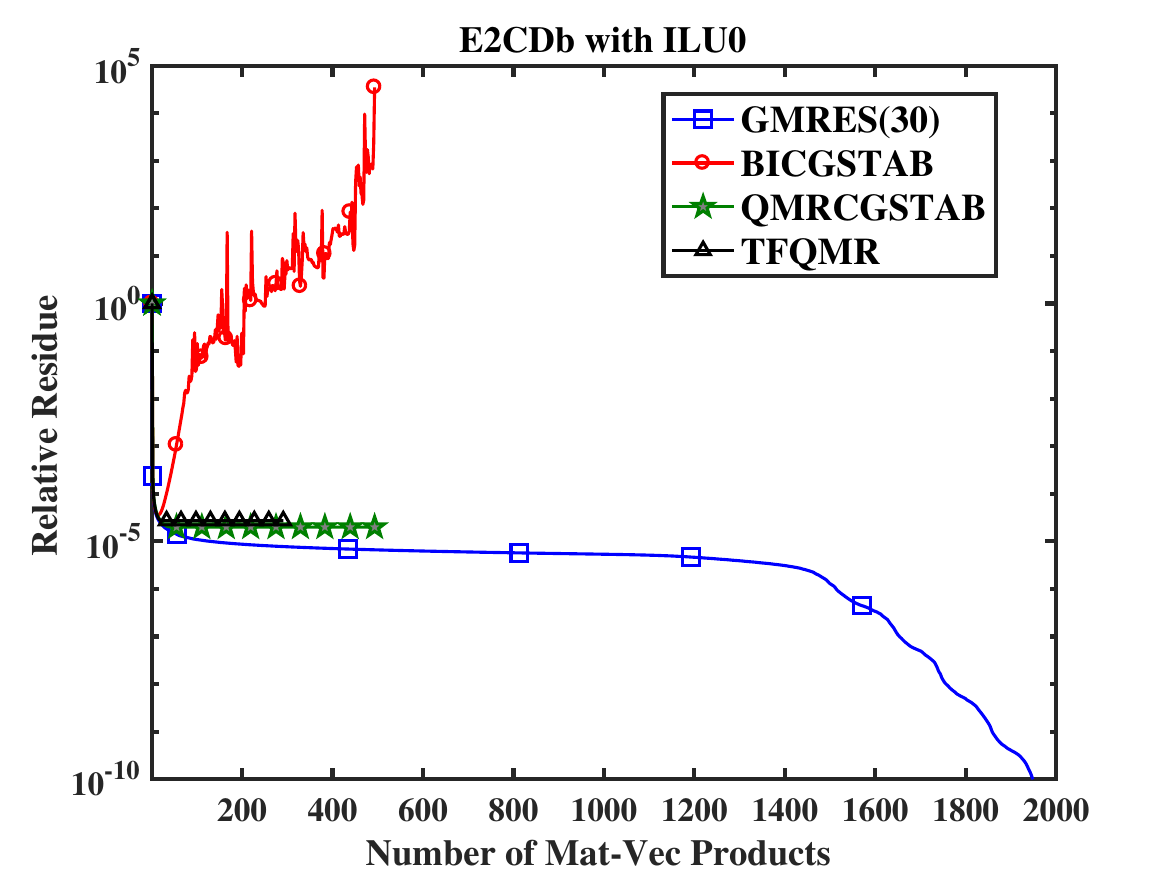}
\par\end{center}%
\end{minipage}\hfill{} %
\begin{minipage}[t]{0.45\textwidth}%
\begin{center}
\includegraphics[width=1\textwidth]{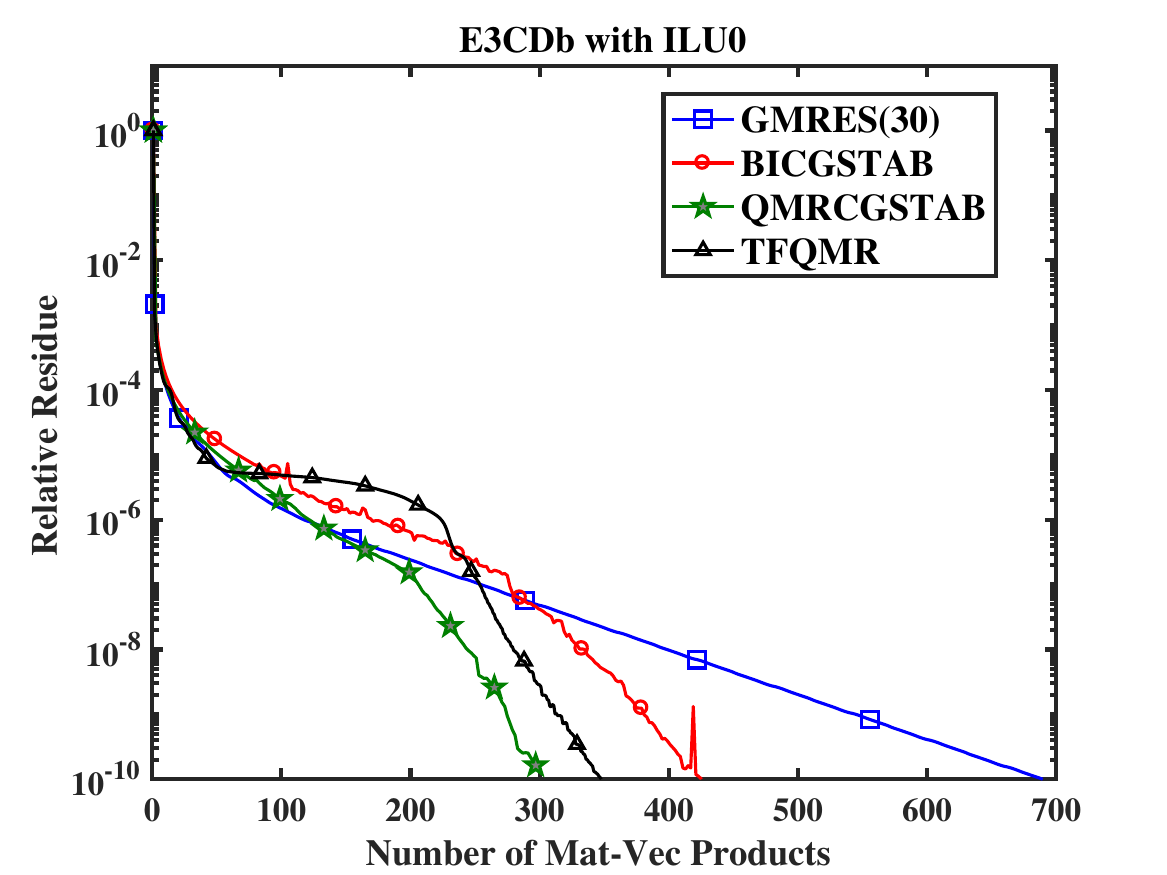}
\par\end{center}%
\end{minipage}

\begin{minipage}[t]{0.45\textwidth}%
\begin{center}
\includegraphics[width=1\textwidth]{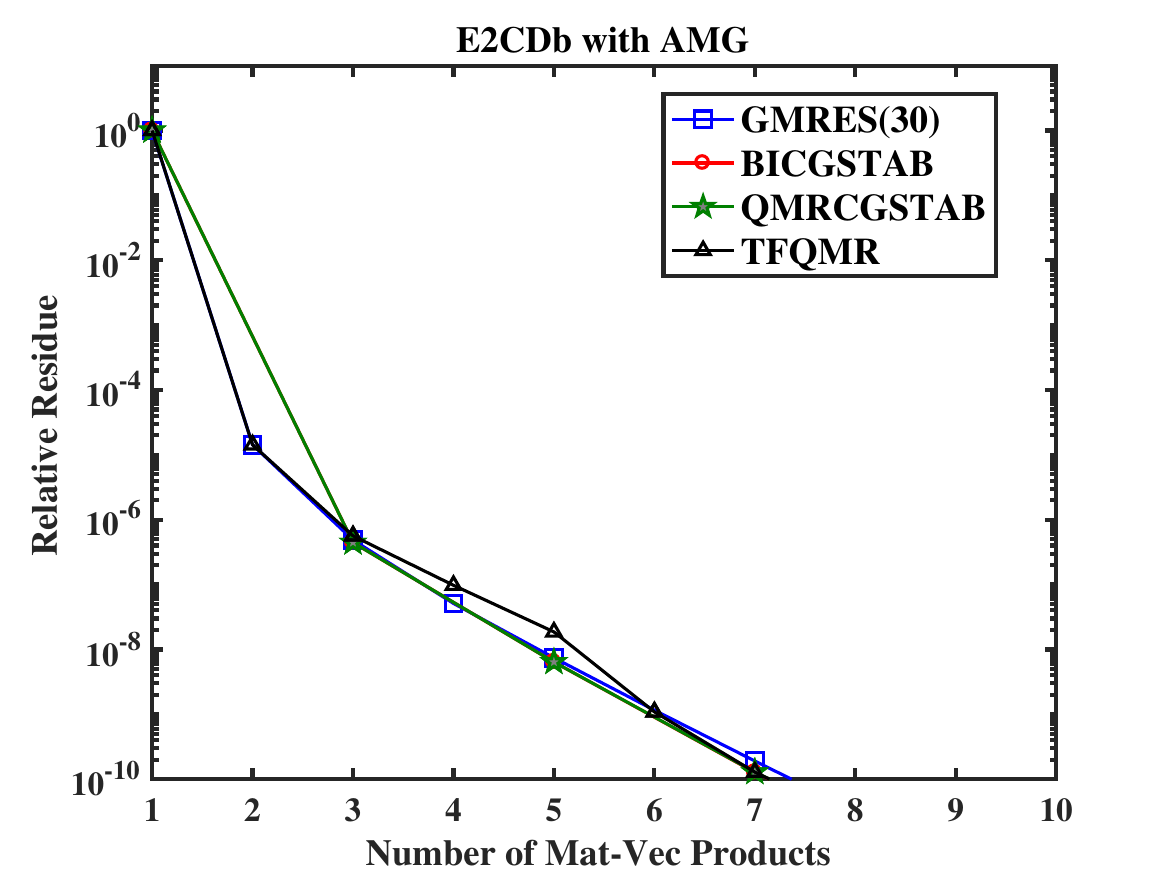}
\par\end{center}%
\end{minipage}\hfill{} %
\begin{minipage}[t]{0.45\textwidth}%
\begin{center}
\includegraphics[width=1\textwidth]{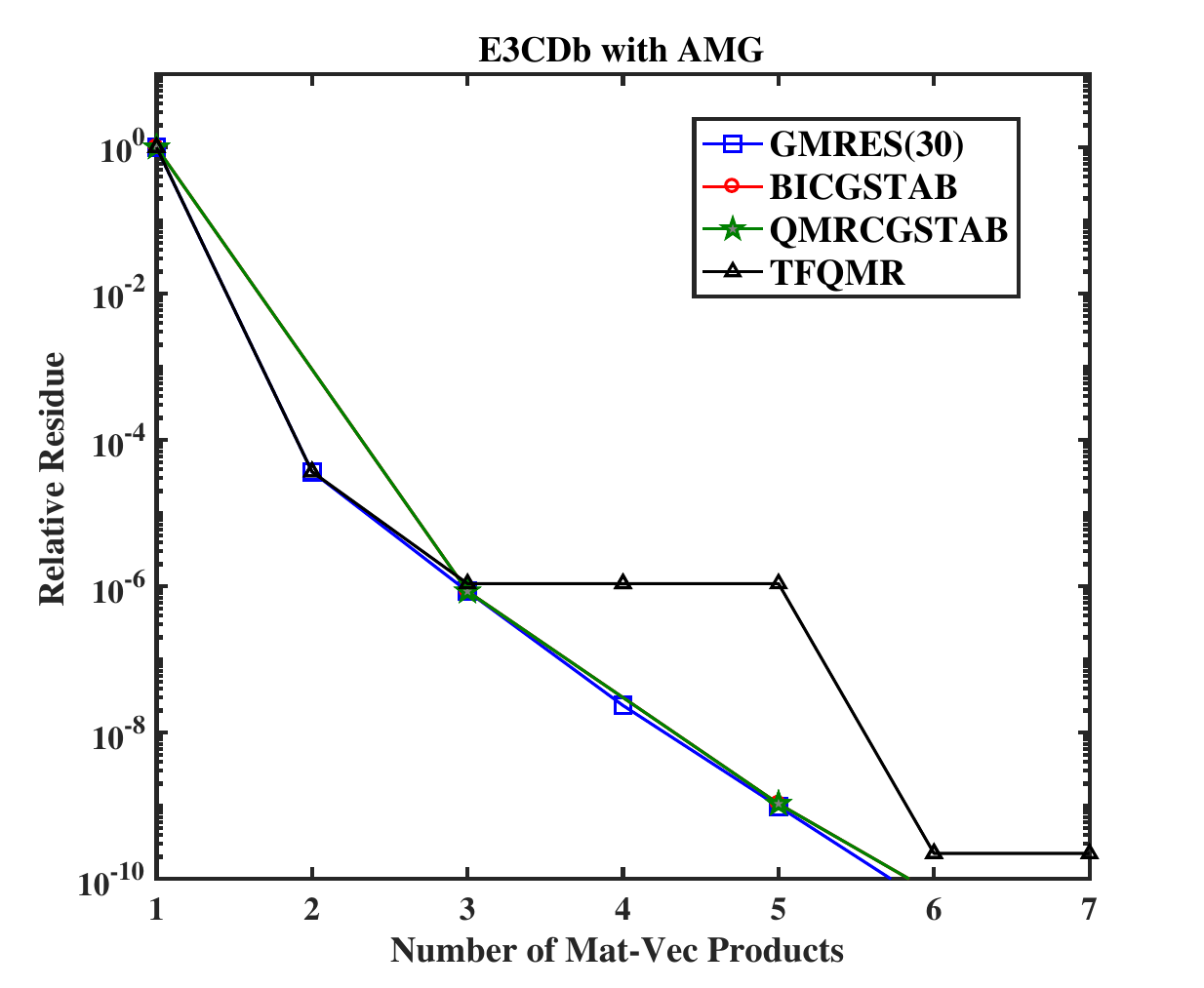}
\par\end{center}%
\end{minipage}

\caption{\textcolor{blue}{\label{FEM-residual-2}}Residuals vs. numbers of
matrix-vector products for E2CDb (left) and E3CDb (right).}
\end{figure}

\begin{figure}[h]
\begin{minipage}[t]{0.45\textwidth}%
\begin{center}
\includegraphics[width=1\textwidth]{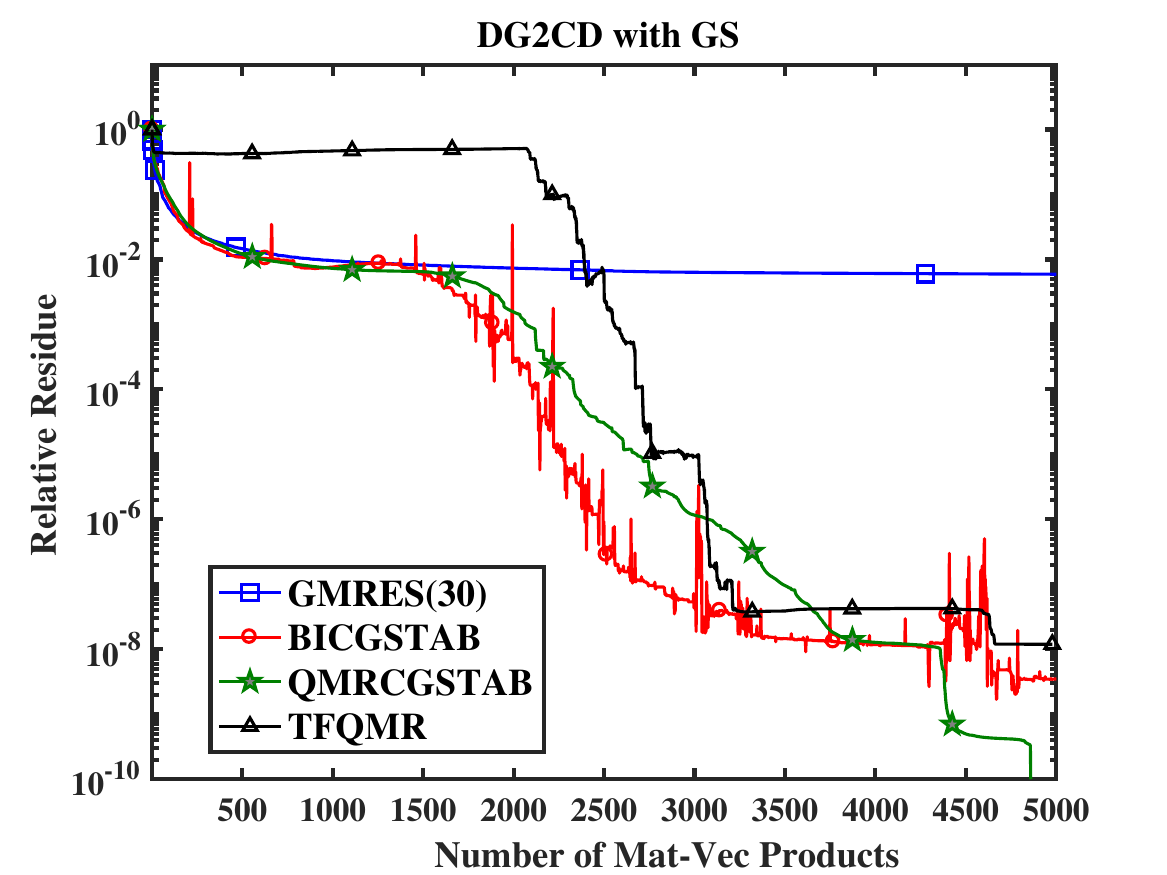}
\par\end{center}%
\end{minipage}\hfill{} %
\begin{minipage}[t]{0.45\textwidth}%
\begin{center}
\includegraphics[width=1\textwidth]{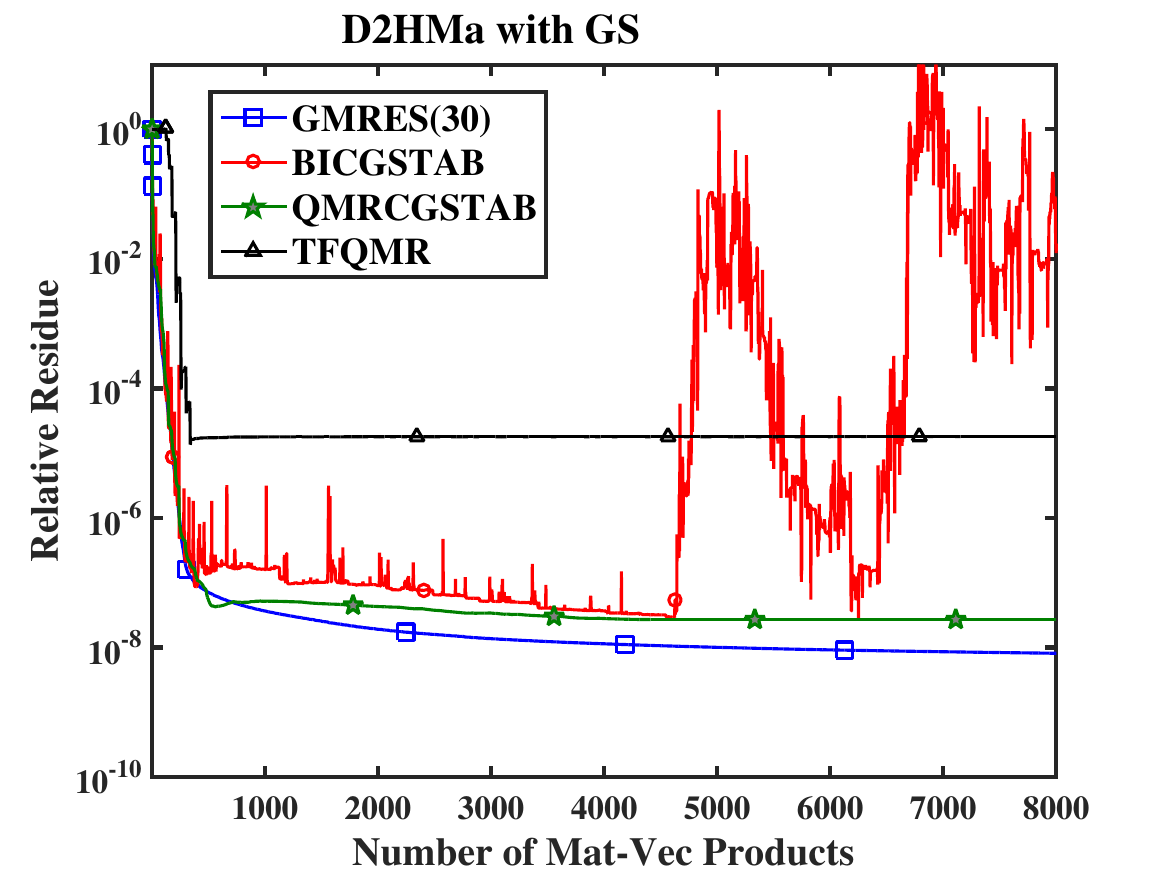}
\par\end{center}%
\end{minipage}

\begin{minipage}[t]{0.45\textwidth}%
\begin{center}
\includegraphics[width=1\textwidth]{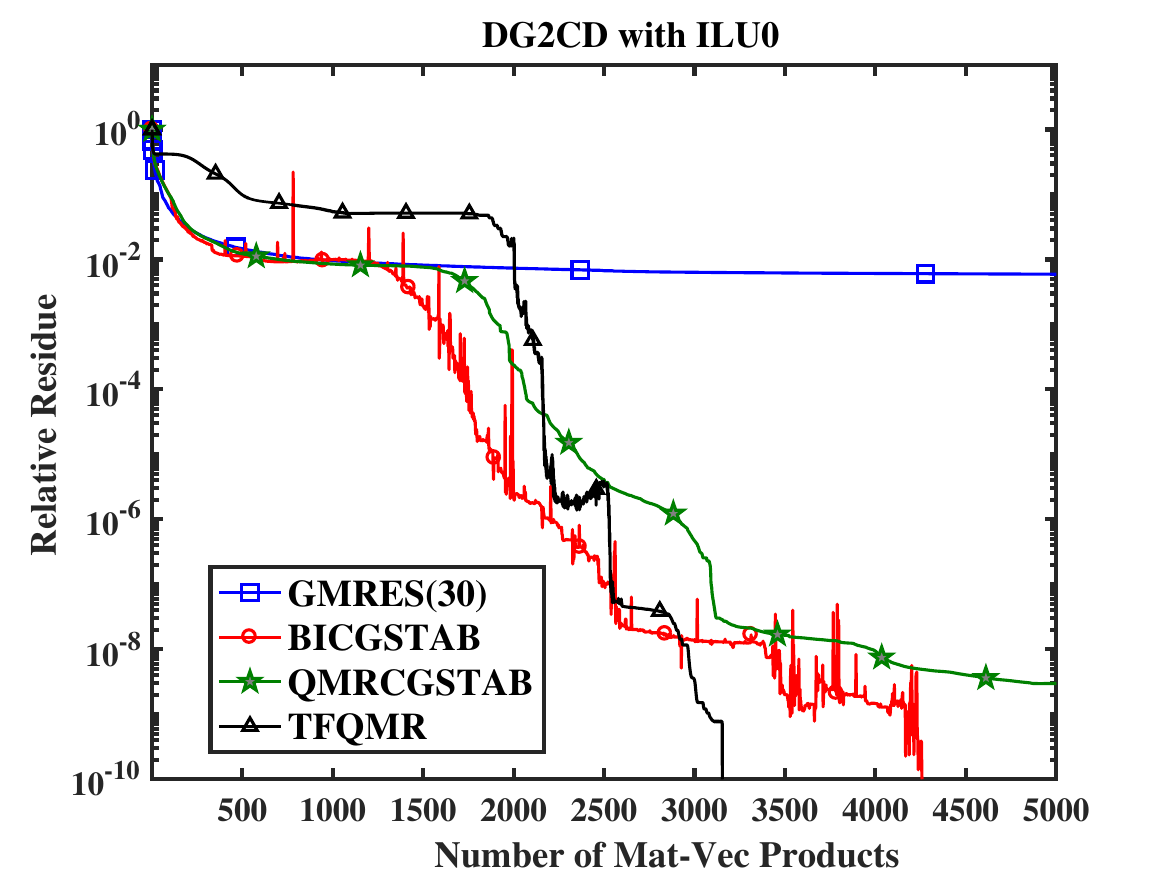}
\par\end{center}%
\end{minipage}\hfill{} %
\begin{minipage}[t]{0.45\textwidth}%
\begin{center}
\includegraphics[width=1\textwidth]{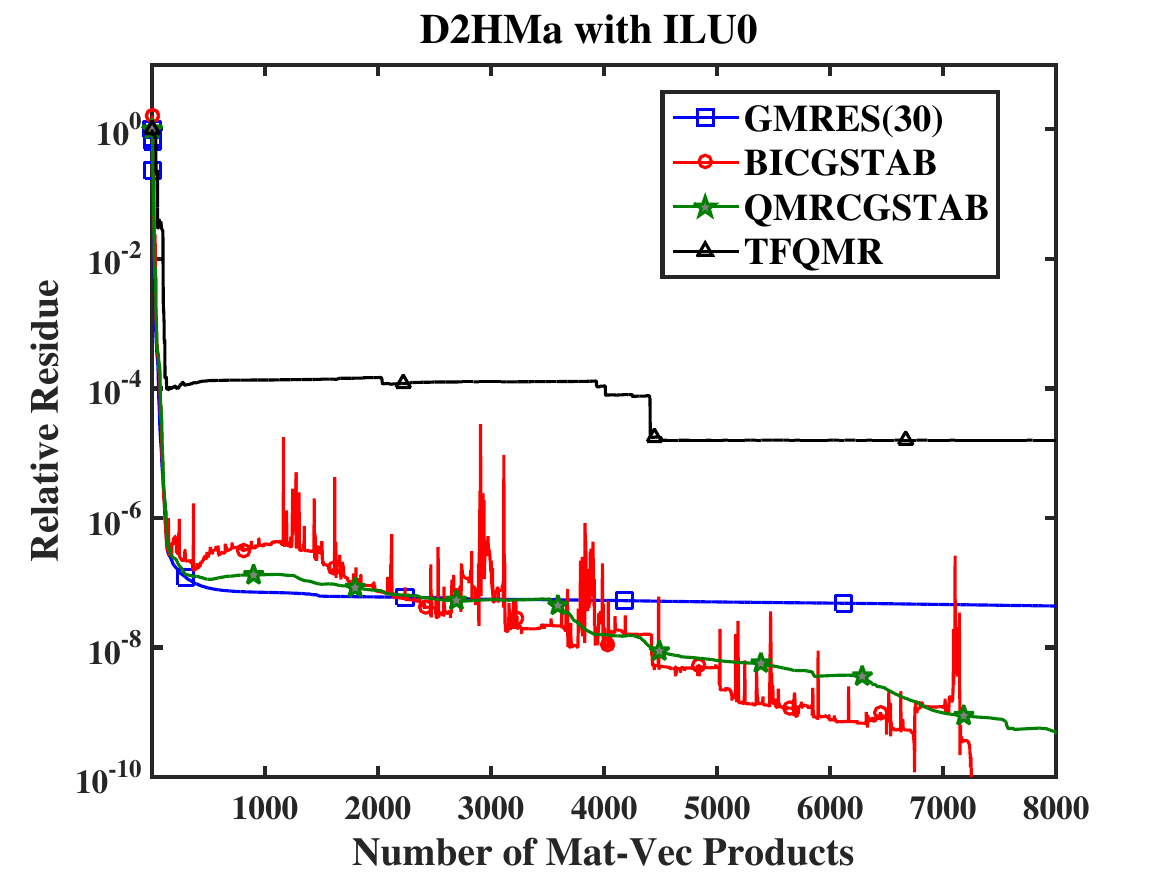}
\par\end{center}%
\end{minipage}

\begin{minipage}[t]{0.45\textwidth}%
\begin{center}
\includegraphics[width=1\textwidth]{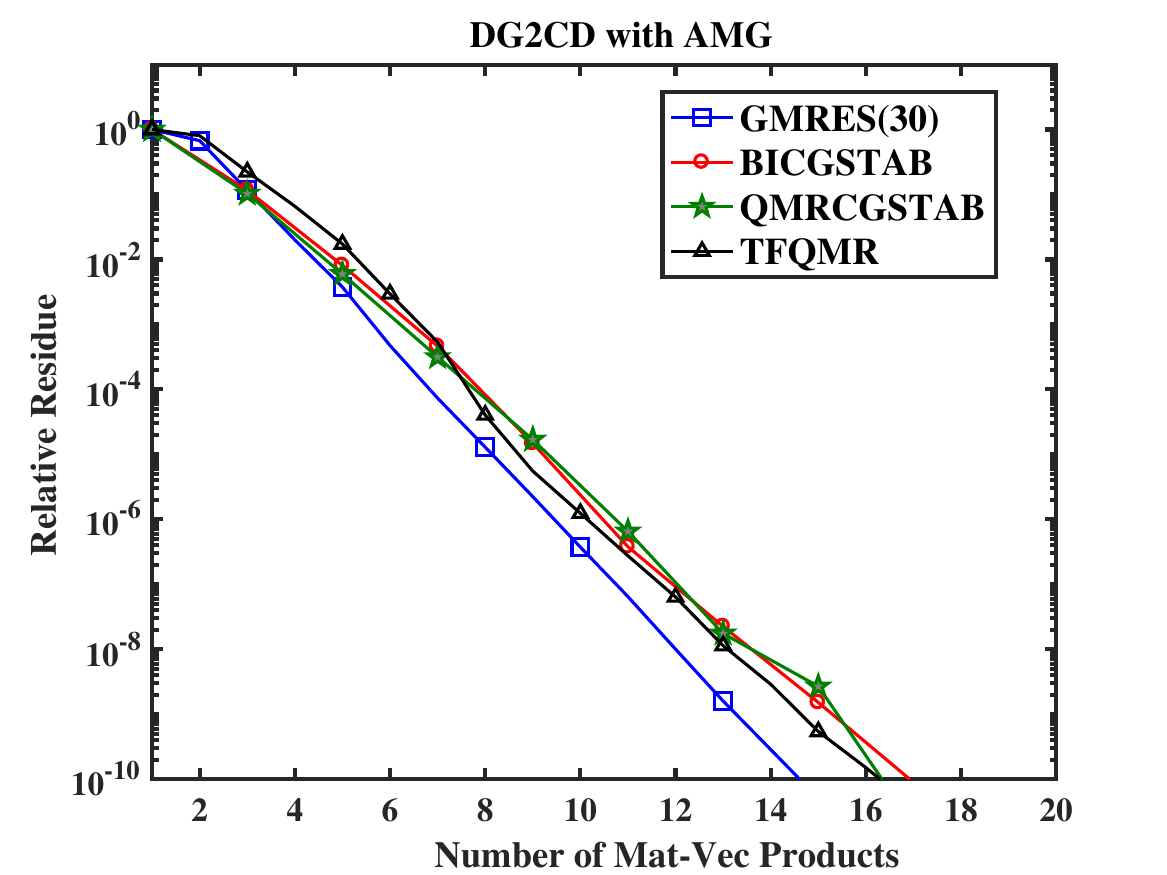}
\par\end{center}%
\end{minipage}\hfill{} %
\begin{minipage}[t]{0.45\textwidth}%
\begin{center}
\includegraphics[width=1\textwidth]{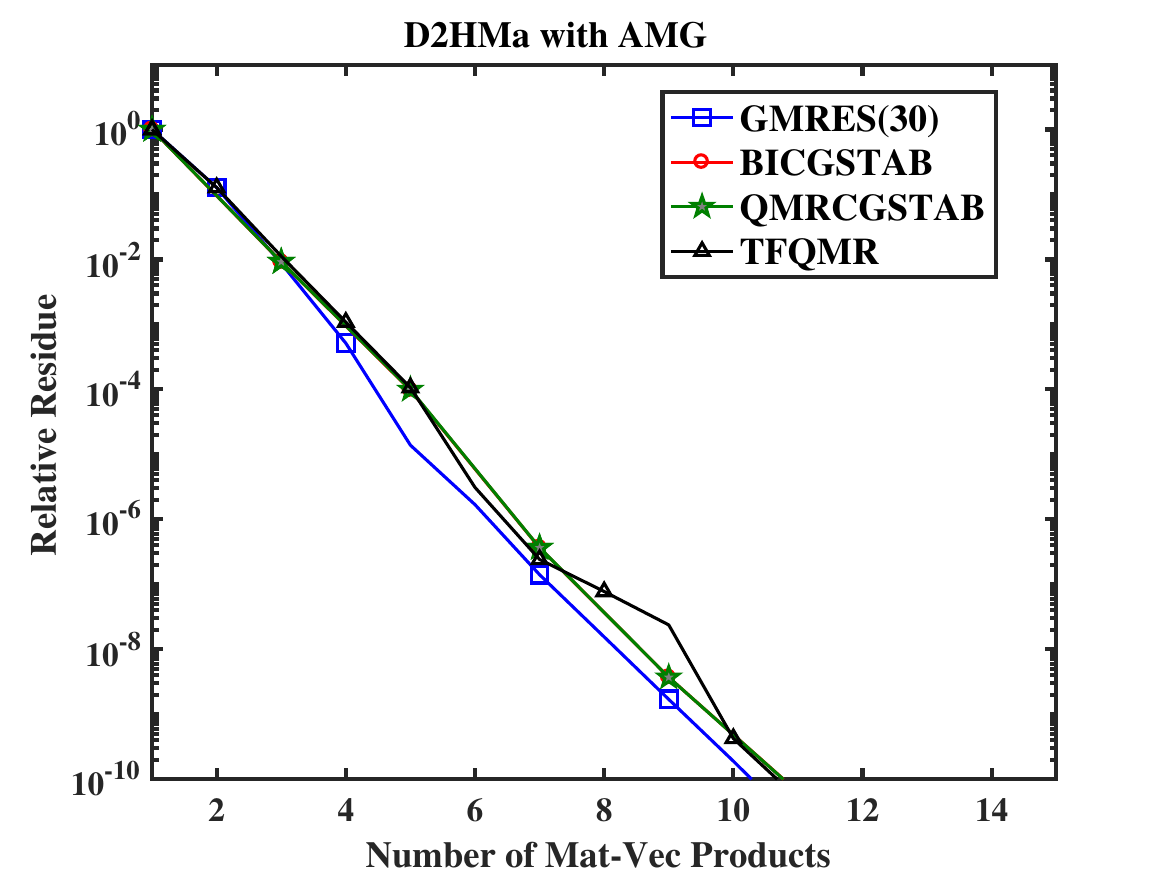}
\par\end{center}%
\end{minipage}

\caption{\textcolor{blue}{\label{fig:FD-DG-residual}}Residuals vs. numbers
of matrix-vector products for DG2CD (left) and D2HMa (right).}
\end{figure}

\begin{figure}[h]
\begin{minipage}[t]{0.45\textwidth}%
\begin{center}
\includegraphics[width=1\textwidth]{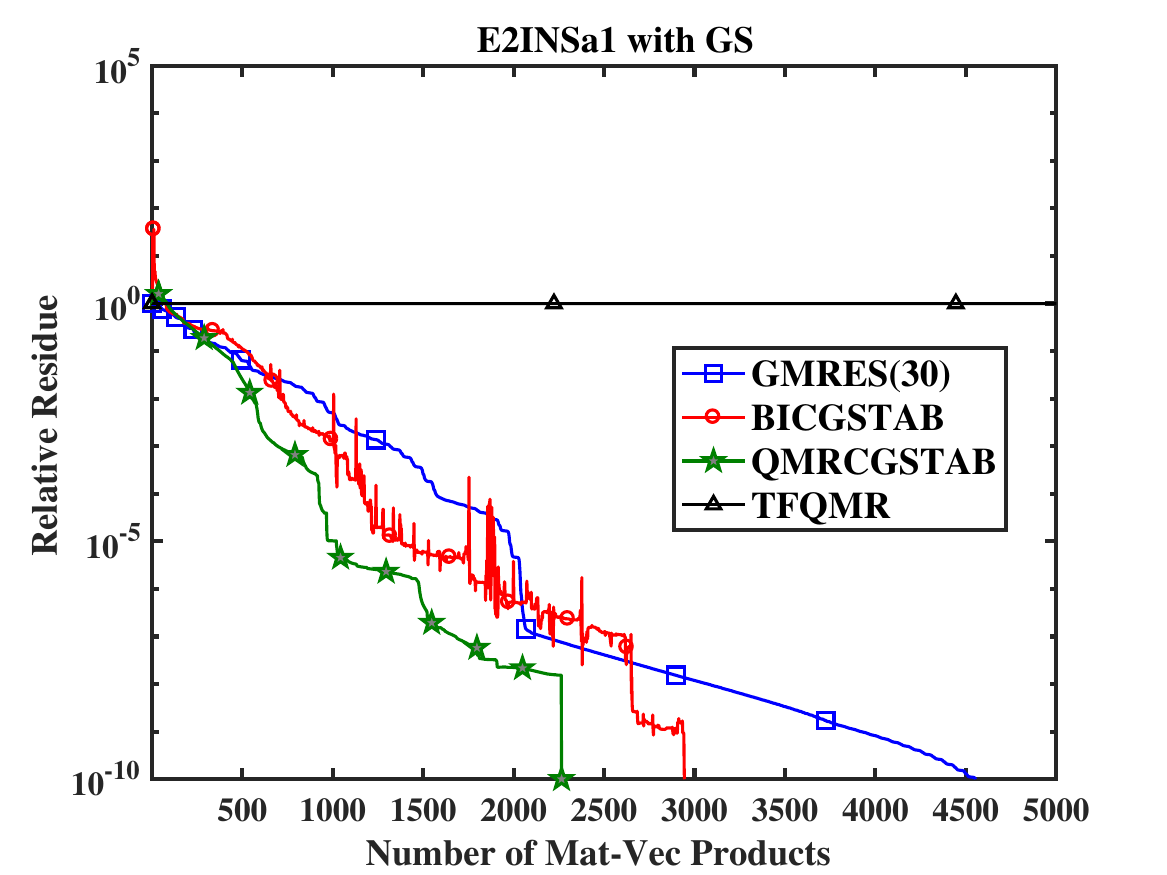}
\par\end{center}%
\end{minipage}\hfill{} %
\begin{minipage}[t]{0.45\textwidth}%
\begin{center}
\includegraphics[width=1\textwidth]{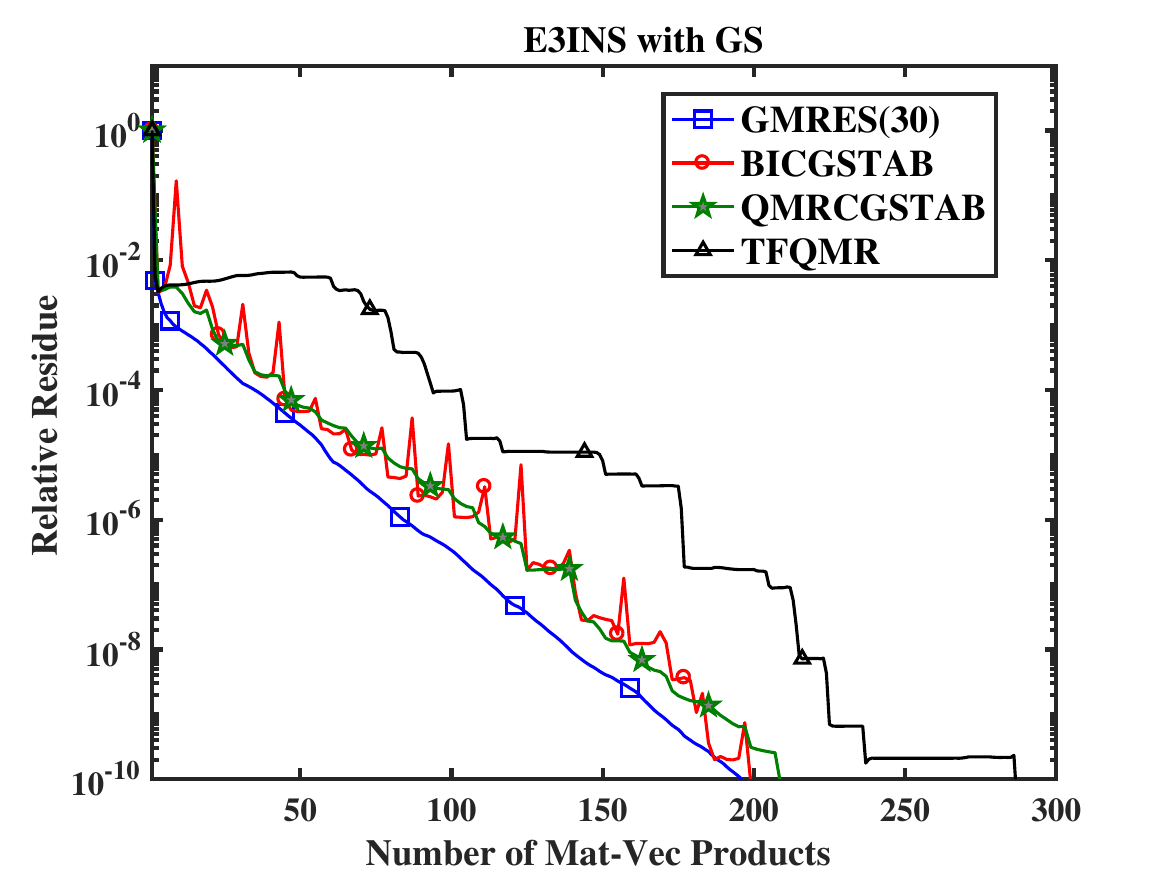}
\par\end{center}%
\end{minipage}

\begin{minipage}[t]{0.45\textwidth}%
\begin{center}
\includegraphics[width=1\textwidth]{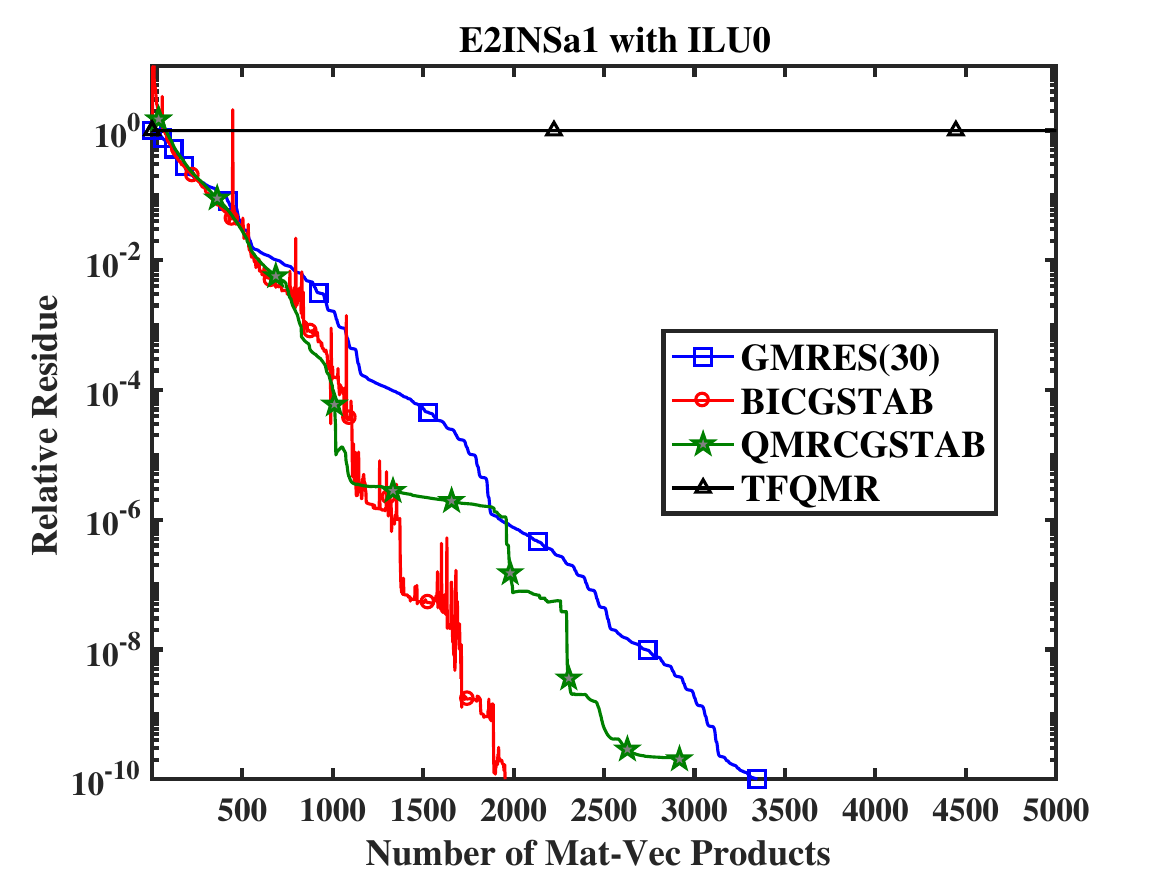}
\par\end{center}%
\end{minipage}\hfill{} %
\begin{minipage}[t]{0.45\textwidth}%
\begin{center}
\includegraphics[width=1\textwidth]{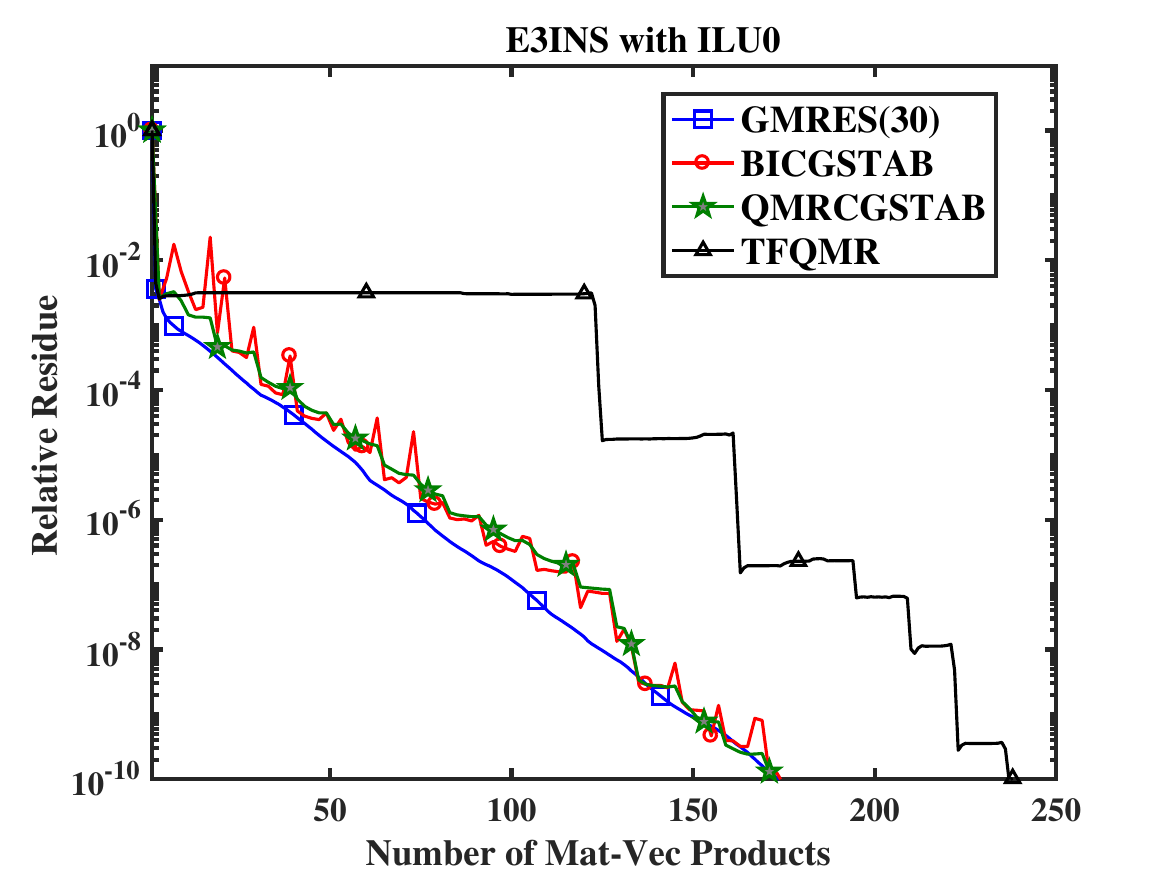}
\par\end{center}%
\end{minipage}

\begin{minipage}[t]{0.45\textwidth}%
\begin{center}
\includegraphics[width=1\textwidth]{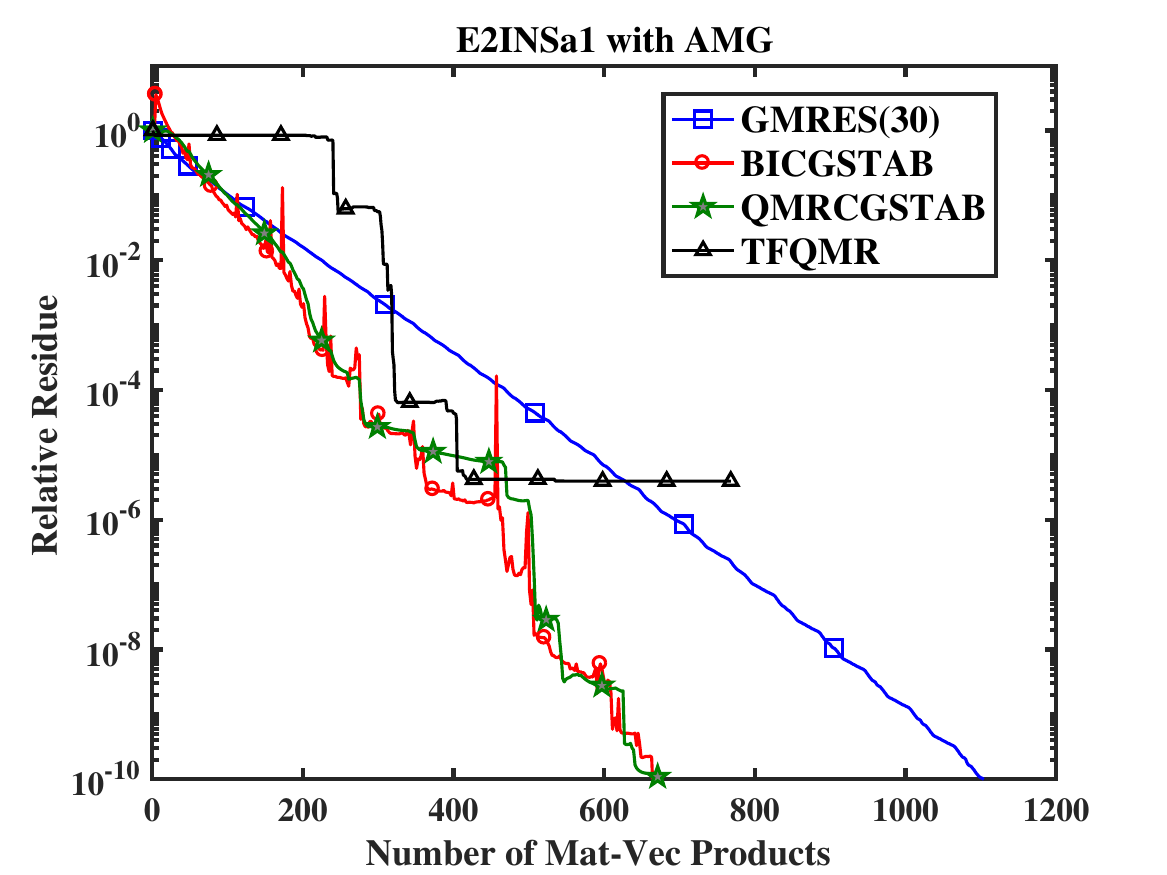}
\par\end{center}%
\end{minipage}\hfill{} %
\begin{minipage}[t]{0.45\textwidth}%
\begin{center}
\includegraphics[width=1\textwidth]{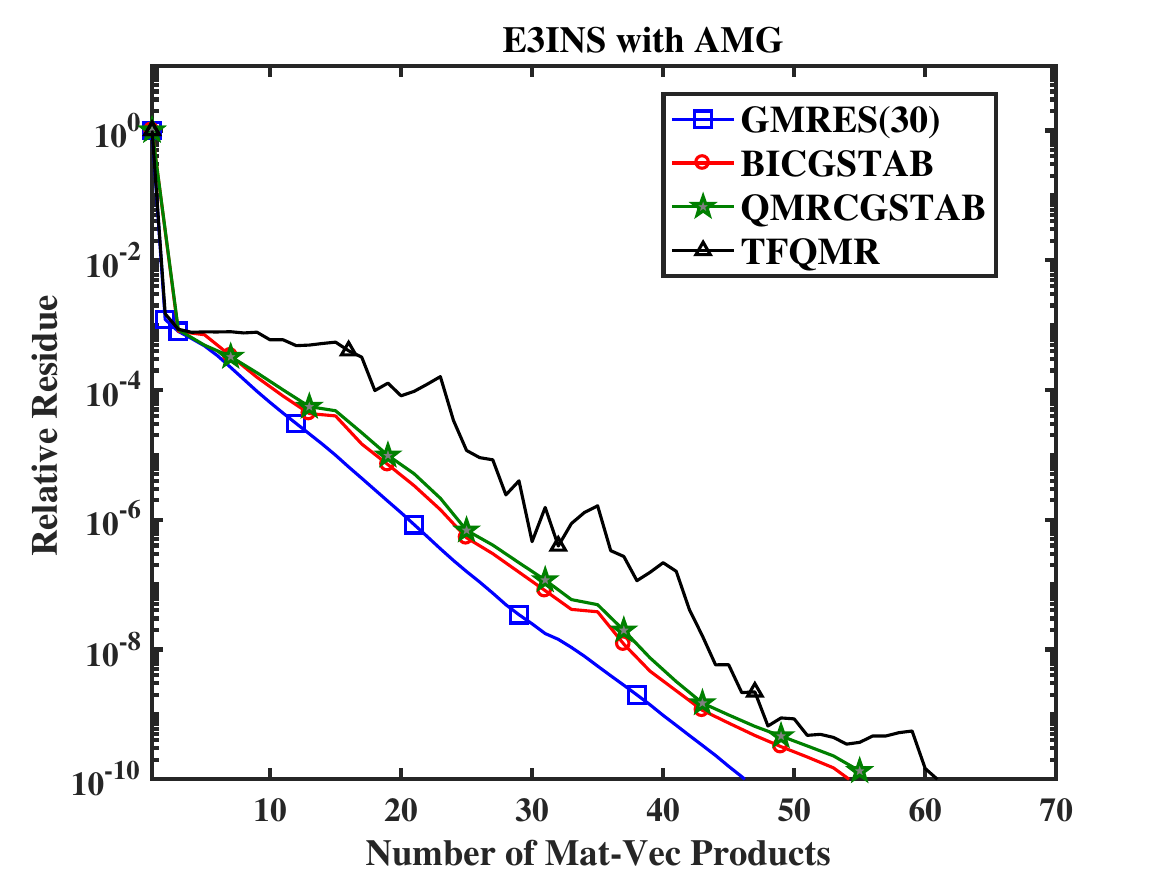}
\par\end{center}%
\end{minipage}

\caption{\textcolor{blue}{\label{fig:NS-residual}}Residuals vs. numbers of
matrix-vector products for E2INSa1 (left) and E3INS (right).}
\end{figure}

Theoretically, the asymptotic convergence of the Krylov subspace methods
may be analyzed using the distributions of eigenvalues and (generalized)
eigenvectors. In practice, however, the convergence is complicated
by the restarts in Arnoldi iterations and the nonorthogonal basis
in the bi-Lanczos iteration. To complement theoretical analysis, we
present the convergence history of ten test cases in Figures\textcolor{blue}{~}\ref{fig:FEM-residual}\textendash \ref{fig:NS-residual},
which are representative of the others. Note that we plot the relative
residuals with respect to the numbers of matrix-vector products instead
of iteration counts, because the former is a better indication of
the overall computational cost and also of the degrees of the polynomials.
For ease of cross-comparison of different preconditioners, we truncated
the $x$-axis to be the same for Gauss-Seidel and ILU0 for each matrix.

Figure~\ref{fig:FEM-residual} shows the convergence histor\textcolor{black}{y
of test cases E2CDa and E3CDa3}\textcolor{blue}{. }With Gauss-Seidel
or ILU0 preconditioners, GMRES converged fast initially but then slowed
down drastically due to restart, whereas BiCGSTAB had highly oscillatory
residuals because its formulation does not minimize the residual in
any norm or pseudo-norm. QMRCGSTAB was much smoother than BiCGSTAB,
and it sometimes converged faster than BiCGSTAB. The convergence of
TFQMR exhibited a staircase pattern, indicating frequent near stagnations
due to its sensitivity to rounding errors. With BoomerAMG, however,
the four methods had about the same convergence trajectories. GMRES
converged slightly faster than the others, but all the methods converged
quite smoothly, without any apparent oscillation or stagnation. These
results indicate that an effective multigrid preconditioner can potentially
overcome the disadvantages of each of these KSP methods, including
slow convergence with restated GMRES, oscillations with BiCGSTAB,
and near stagnation with TFQMR. Figure~\ref{fig:GFD-residual} shows
the convergence results for GD2CD and GD3CD. The results are qualitatively
similar to those of E2CDa and E3CDa3, except that the near stagnation
of TFQMR is even more apparent. 

Figure~\ref{FEM-residual-2} shows the convergence results for E2CDb
and E3CDb. The most notable new feature is that for E2CDb, BiCGSTAB
diverged with Gauss-Seidel and ILU0. QMRCGSTAB with Gauss-Seidel converged,
so did GMRES with ILU0, but all other methods with Gauss-Seidel and
ILU0 stagnated. For E3CDb, all the methods converged with Gauss-Seidel
and ILU0. Although BiCGSTAB converged relatively smoothly, the residuals
had notable jumps near the end. TFQMR had a plateau even with AMG,
again due to its sensitivity to round-off errors. Note that QMRCGSTAB
overcame the oscillations with BiCGSTAB and the near stagnations with
TFQMR, and it converged the fastest in some cases.

Figure~\ref{fig:FD-DG-residual} shows the convergence results for
DG2CD and D2HMa. The condition number of DG2CD is about $10^{7}$,
and its results were qualitatively similar as GD2CD and GD3CD. The
condition number of D2HMa is nearly $10^{9}$, and this ill-conditioning
caused difficulties for virtually all the methods with GS and ILU0
preconditioners. GMRES and TFQMR both stagnated. BiCGSTAB oscillated
wildly, while QMRCGSTAB had smoother trajectories. With BoomerAMG,
all the methods converged rapidly and smoothly, while GMRES converged
slightly faster than the others.

Figure~\ref{fig:NS-residual} shows the convergence results for E2INSa1
and E3INS. For E2INSa1, TFQMR stagnated with all the preconditioners.
The other solvers delivered similar performance to each other. GMRES
performed the best for E3INS. However, QMRCGSTAB and BiCGSTAB outperformed
GMRES for E2INSa1, which required many iterations even with AMG. 

From these results, we make the following observations. If relatively
few iterations are needed to converge, especially with an effective
multigrid preconditioner, GMRES converges the fastest, because it
minimizes the 2-norm of the residual. However, if many iterations
are needed, then the advantage of GMRES diminishes, and a method with
a three-term recurrence can be more reliable. BiCGSTAB and TFQMR suffer
from oscillatory residuals and frequent plateaus, respectively. By
design, QMRCGSTAB overcomes both of these shortcomings \cite{Chan1998}.

\subsubsection{\label{subsec:Timing-Comparison}Timing Comparison.}

\begin{figure}[h]
\begin{minipage}[t]{0.45\textwidth}%
\begin{center}
\includegraphics[width=1\textwidth]{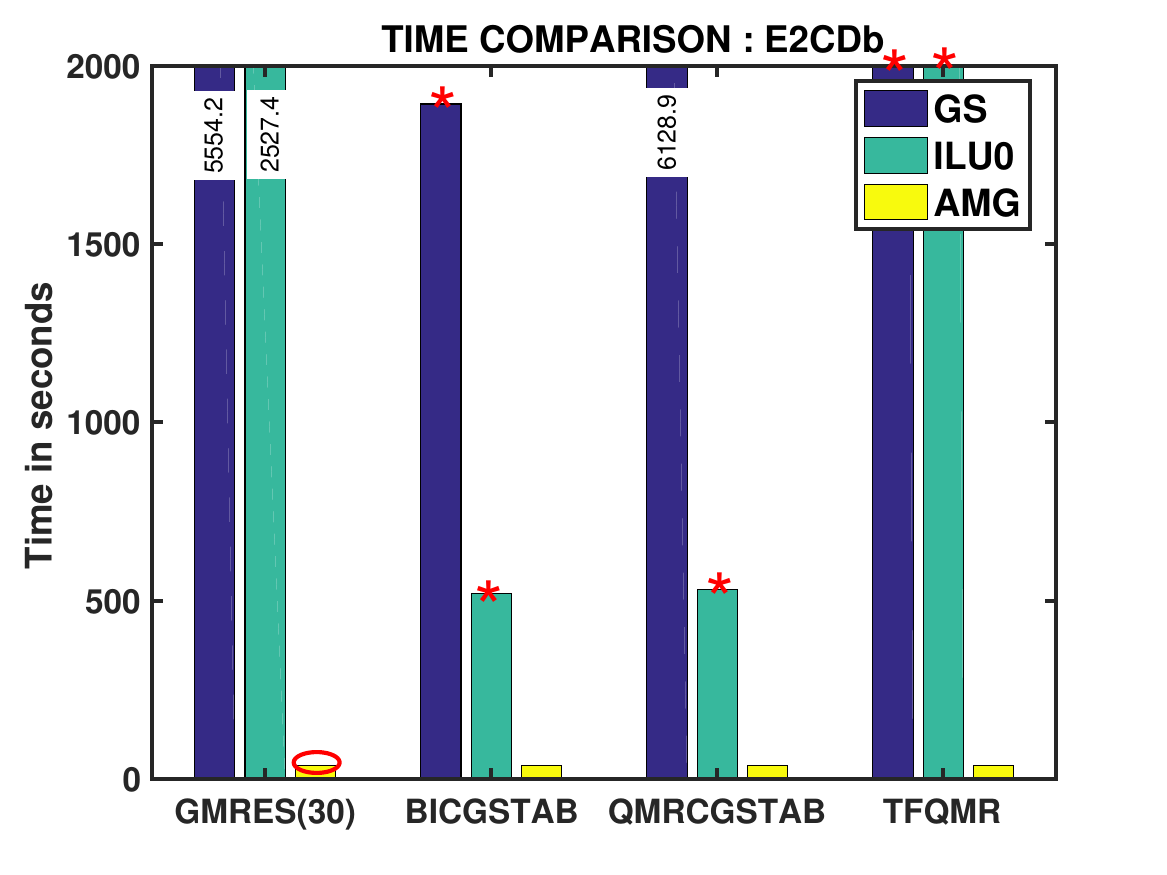}
\par\end{center}%
\end{minipage}\hfill{} %
\begin{minipage}[t]{0.45\textwidth}%
\begin{center}
\includegraphics[width=1\textwidth]{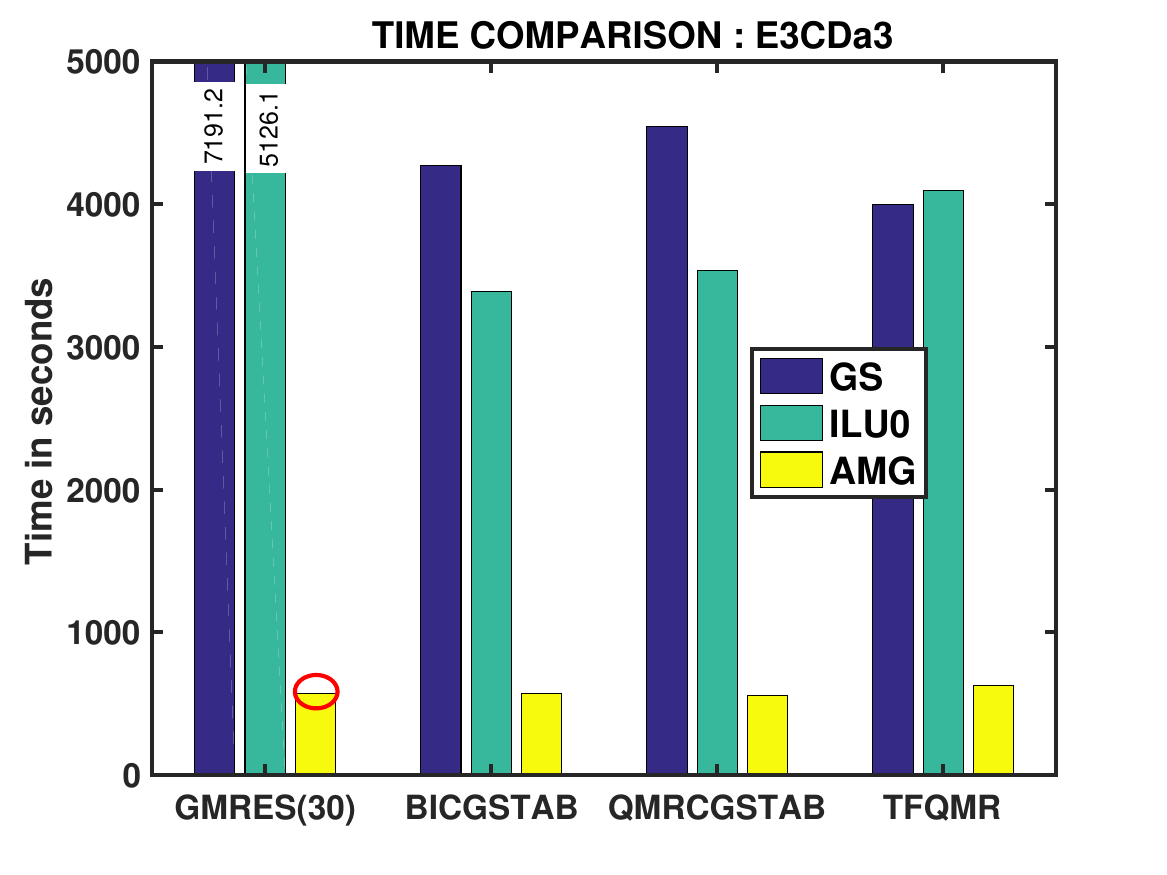}
\par\end{center}%
\end{minipage}

\begin{minipage}[t]{0.45\textwidth}%
\begin{center}
\includegraphics[width=1\textwidth]{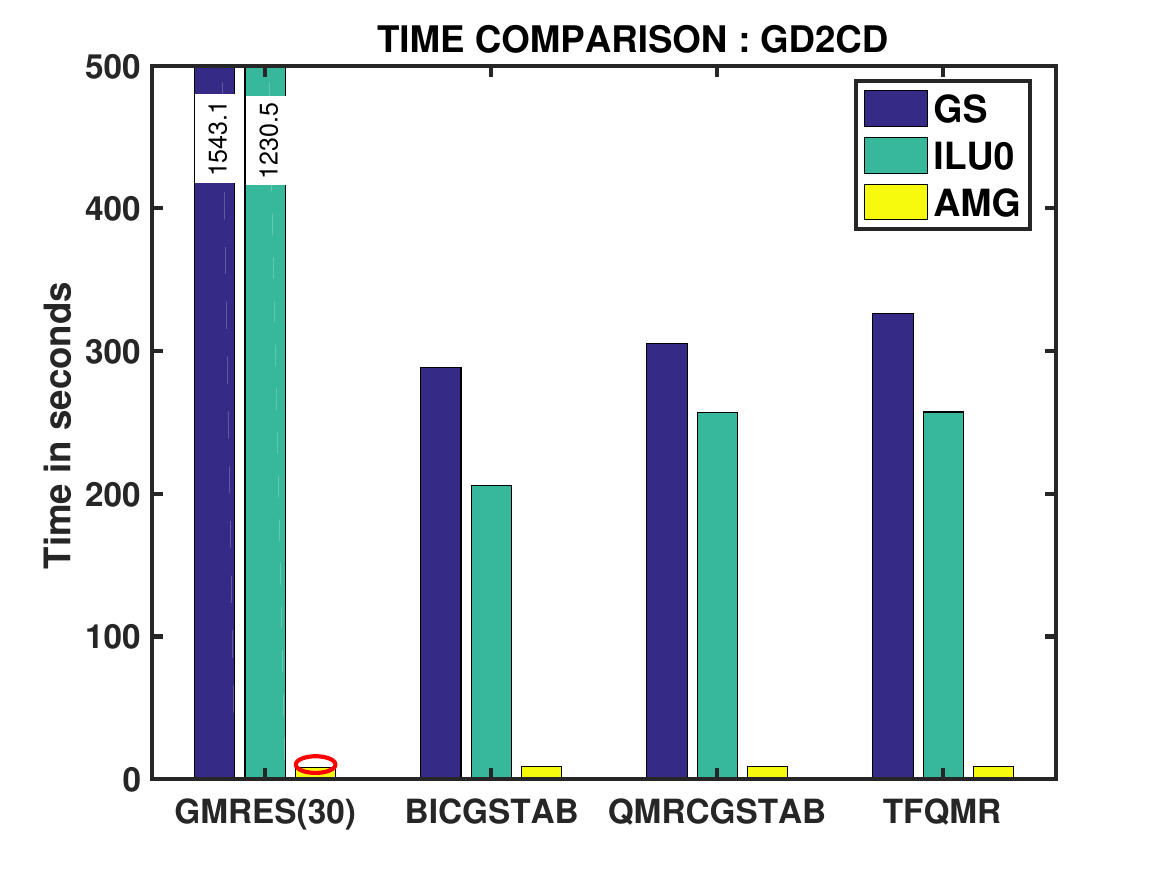}
\par\end{center}%
\end{minipage}\hfill{} %
\begin{minipage}[t]{0.45\textwidth}%
\begin{center}
\includegraphics[width=1\textwidth]{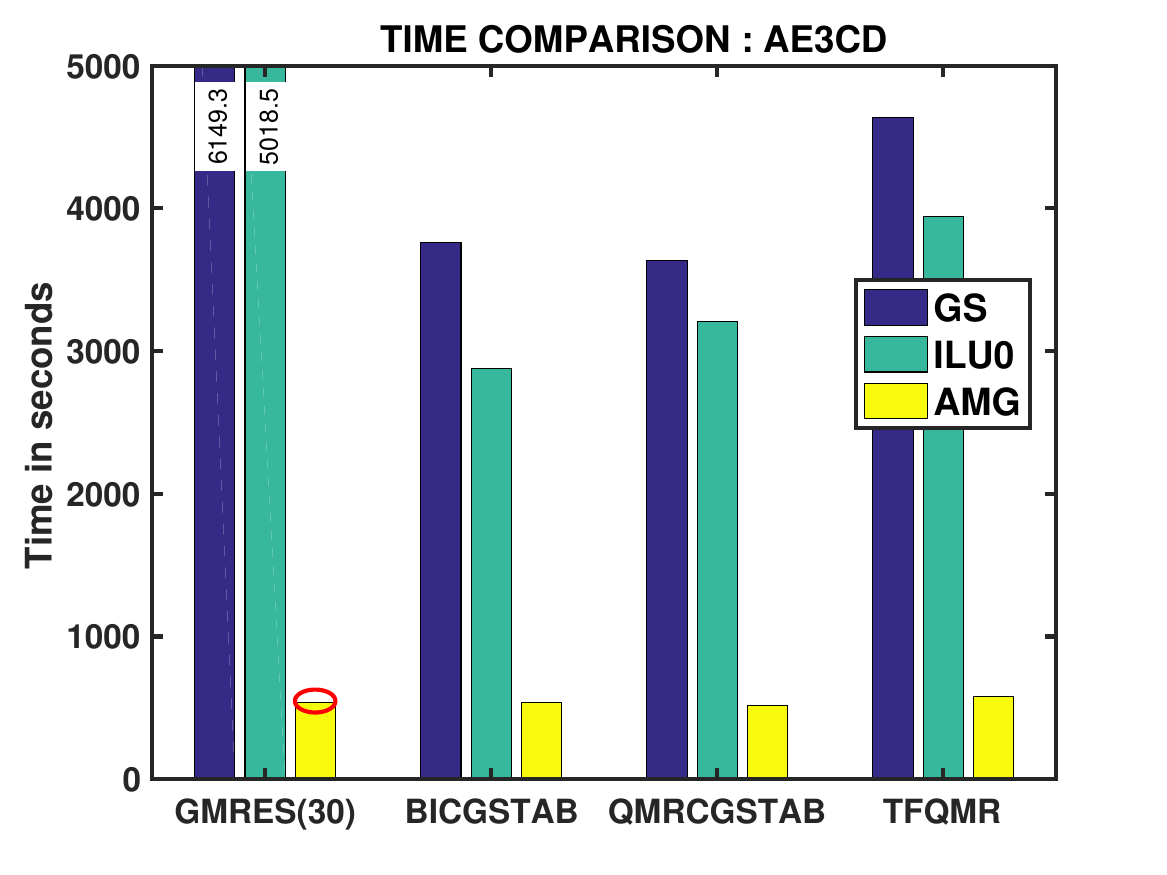}
\par\end{center}%
\end{minipage}

\begin{minipage}[t]{0.45\textwidth}%
\begin{center}
\includegraphics[width=1\textwidth]{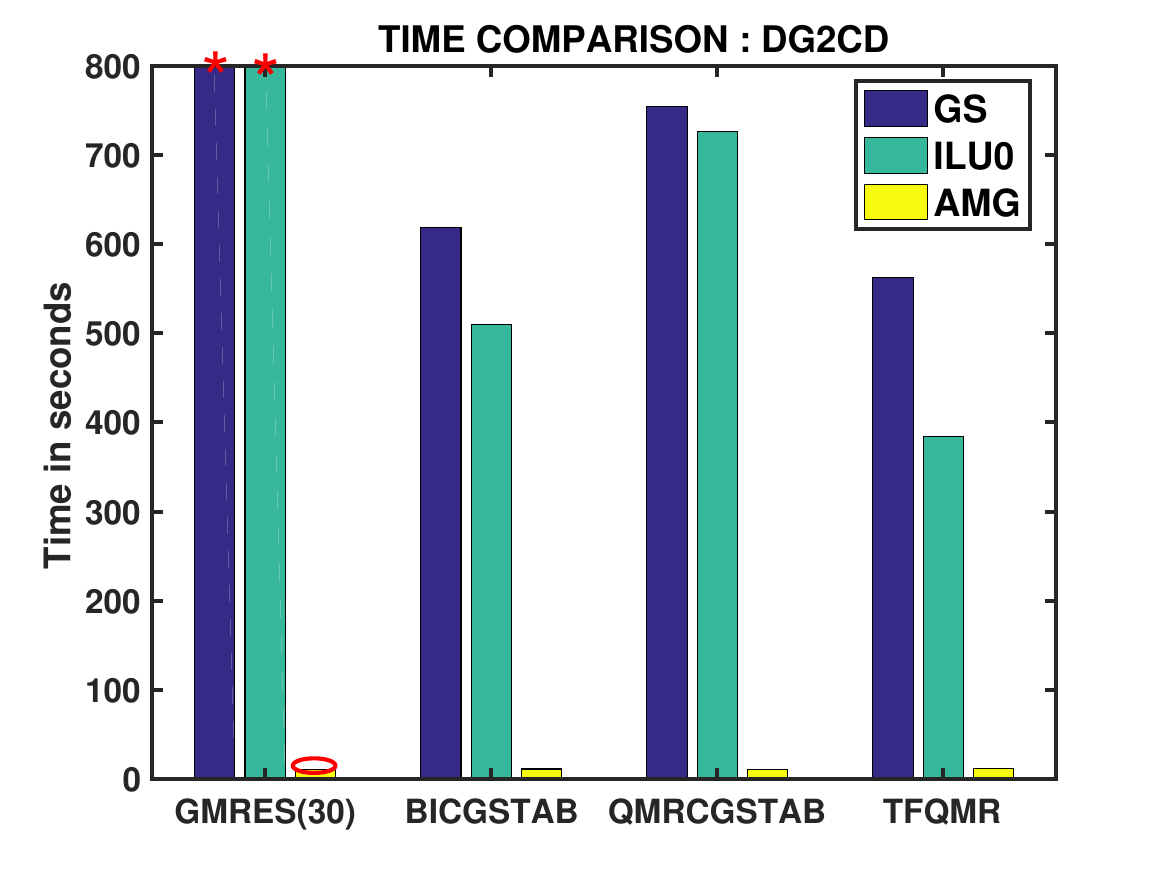}
\par\end{center}%
\end{minipage}\hfill{} %
\begin{minipage}[t]{0.45\textwidth}%
\begin{center}
\includegraphics[width=1\textwidth]{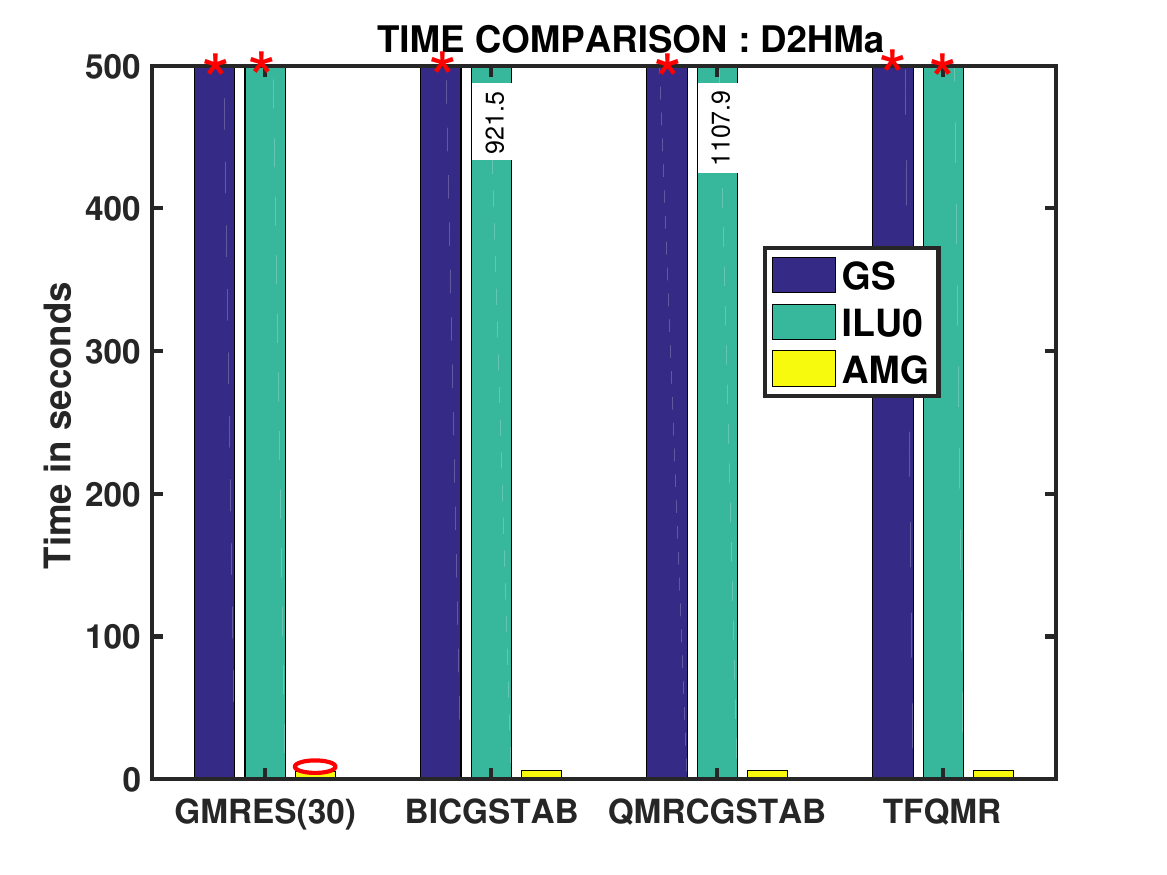}
\par\end{center}%
\end{minipage}

\begin{minipage}[t]{0.45\textwidth}%
\begin{center}
\includegraphics[width=1\textwidth]{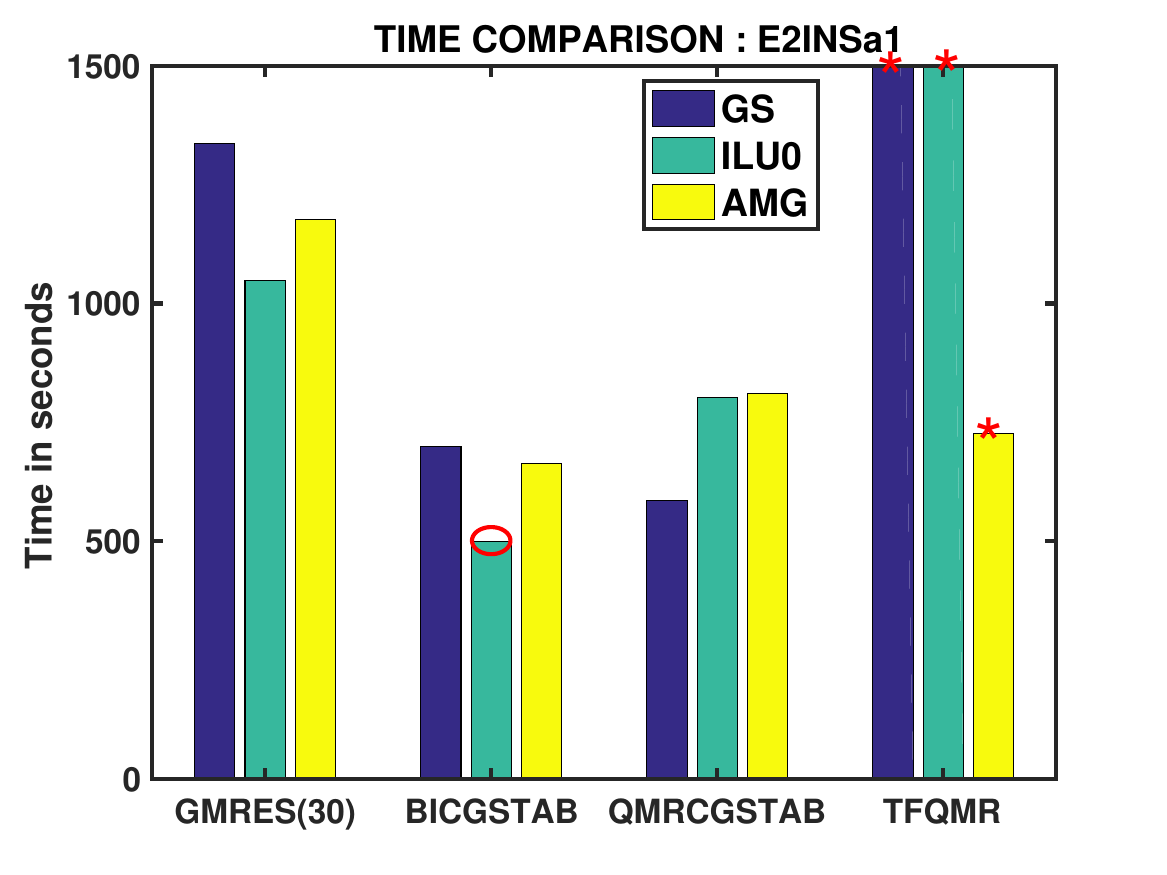}
\par\end{center}%
\end{minipage}\hfill{} %
\begin{minipage}[t]{0.45\textwidth}%
\begin{center}
\includegraphics[width=1\textwidth]{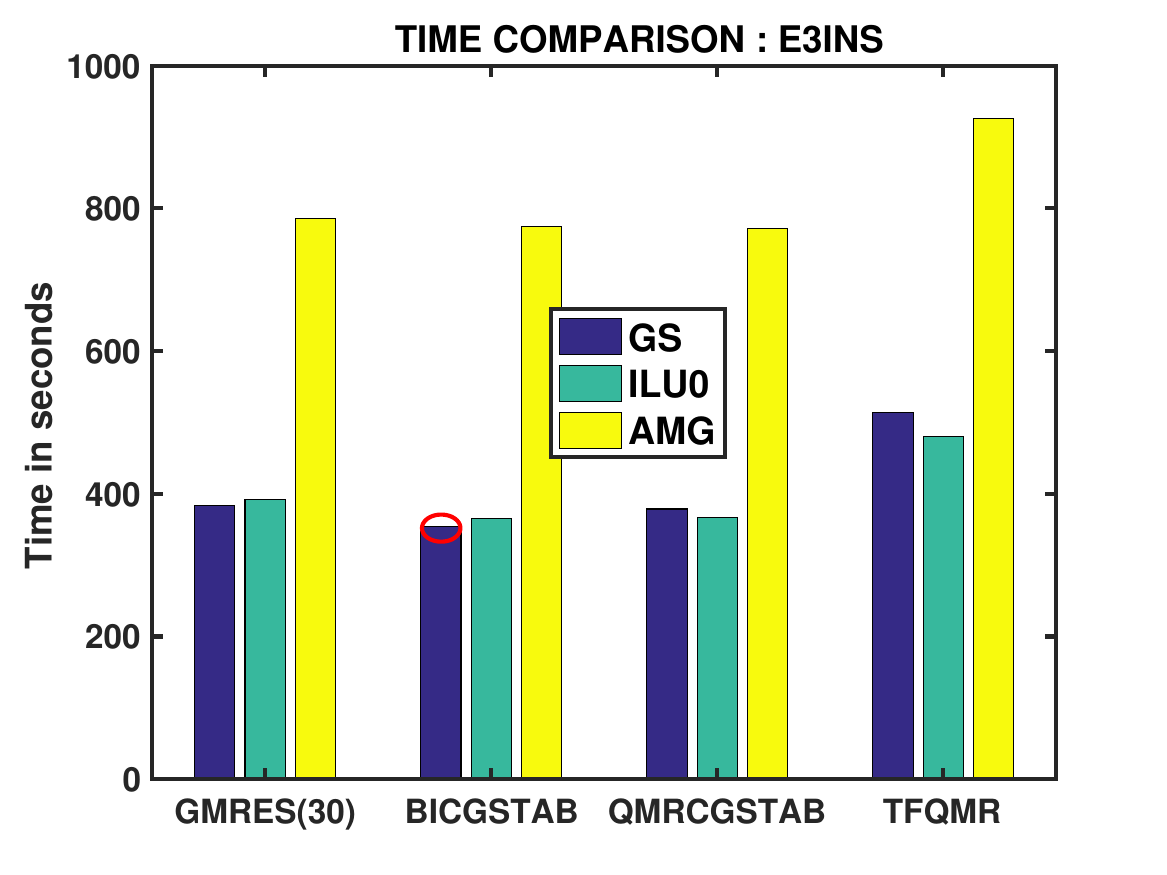}
\par\end{center}%
\end{minipage}

\caption{\label{fig:Timing}Comparison of timing results for eight representative
cases. Circled bars indicate best results. `{*}' indicates non-convergence
or stagnation after 10,000 iterations.}
\end{figure}

The convergence histories are helpful in revealing the intrinsic properties
of the KSP methods, but in practice the overall runtime is the ultimate
criterion. Figure~\ref{fig:Timing} compares the runtimes of seven
cases from those in Section~\ref{subsec:Convergence-Comparison.}
along with AE3CD. We circled out the best performances for each case.
It can be seen that for six out of eight cases, BoomerAMG accelerated
KSP significantly better than Gauss-Seidel and ILU0. GMRES with BoomerAMG
was the best in these cases, while BiCGSTAB and QMRCGSTAB were close
runners-up. With GS and ILU0, GMRES was significantly slower than
the others due to restarts. For Navier-Stokes equations, BiCGSTAB
with Gauss-Seidel and ILU0 delivered better performance than AMG.
QMRCGSTAB had similar performance as BiCGSTAB.

Overall, the timing results are consistent with the convergence results,
so the numbers of matrix-vector products are indeed good predictors
of the overall performance. Among the bi-Lanczos-based methods, BiCGSTAB
is slightly more efficient due to its lower cost per iteration, but
QMRCGSTAB is a competitive alternative, especially if smooth convergence
is desired. In addition, we observe that the continuum formulations
of the PDEs, rather than the discretization methods, tend to have
a bigger impact on the choice of preconditioners. In particular, AMG
can significantly outperform Gauss-Seidel and ILU0 for convection-diffusion
equations, independently of discretization methods, but AMG is not
very effective for incompressible Navier-Stokes equations.

\subsubsection{Asymptotic Growth with Respect to Problem Size.}

\begin{figure}[h]
\begin{minipage}[t]{0.45\textwidth}%
\begin{center}
\includegraphics[width=1\textwidth]{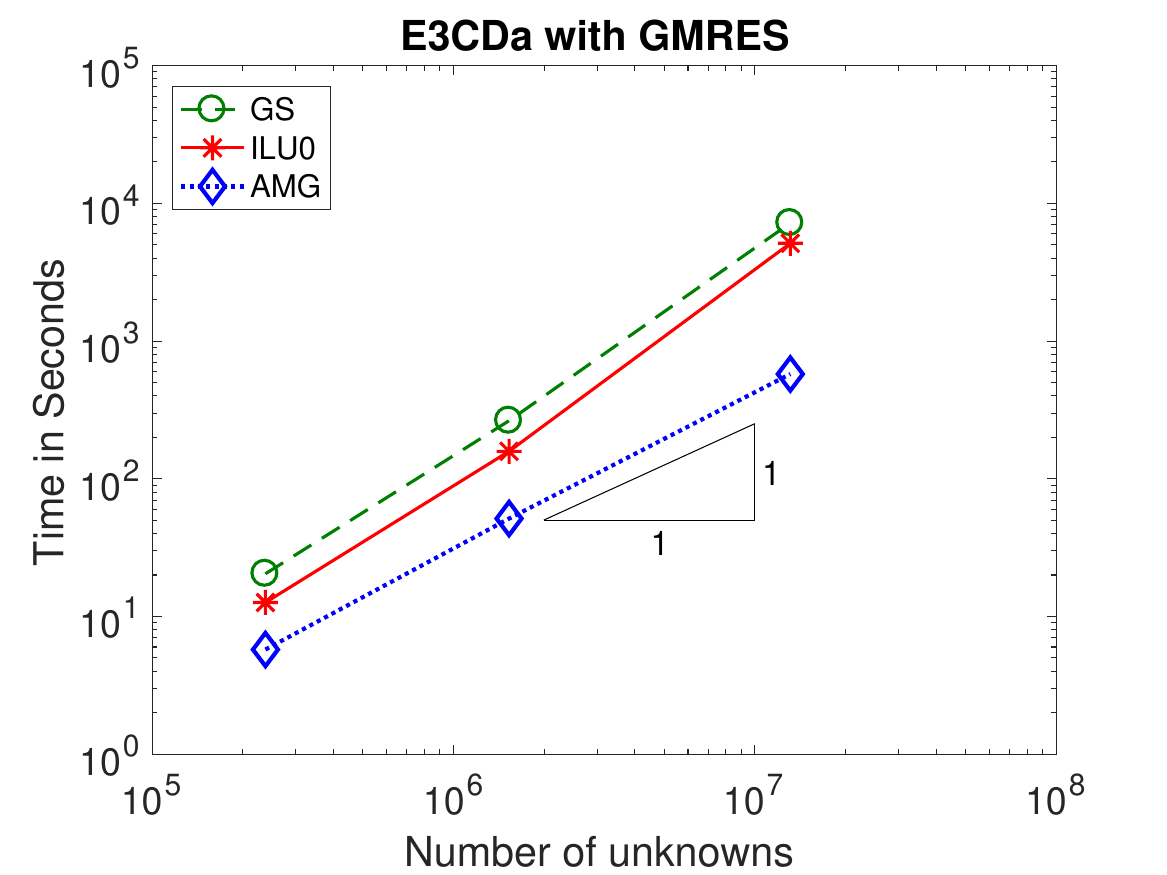}
\par\end{center}%
\end{minipage}\hfill{} %
\begin{minipage}[t]{0.45\textwidth}%
\begin{center}
\includegraphics[width=1\textwidth]{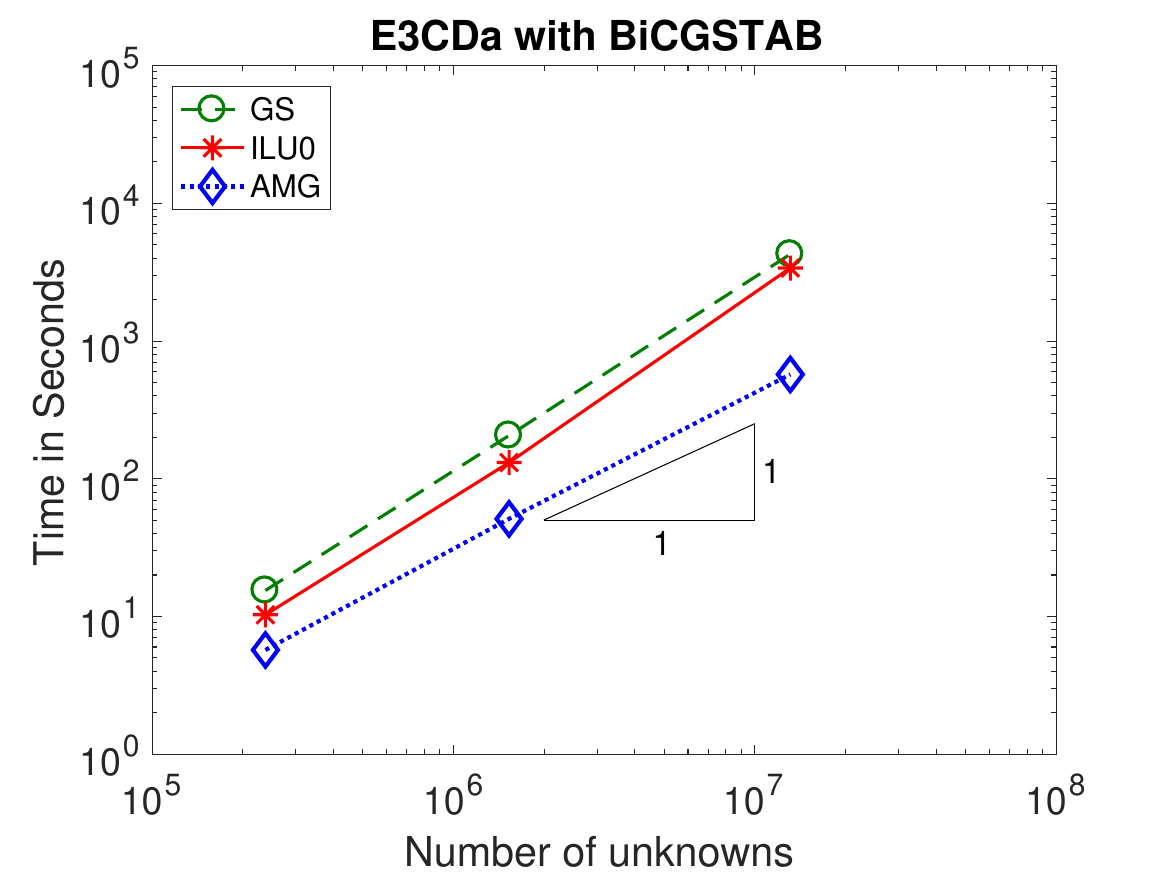}
\par\end{center}%
\end{minipage}

\begin{minipage}[t]{0.45\textwidth}%
\begin{center}
\includegraphics[width=1\textwidth]{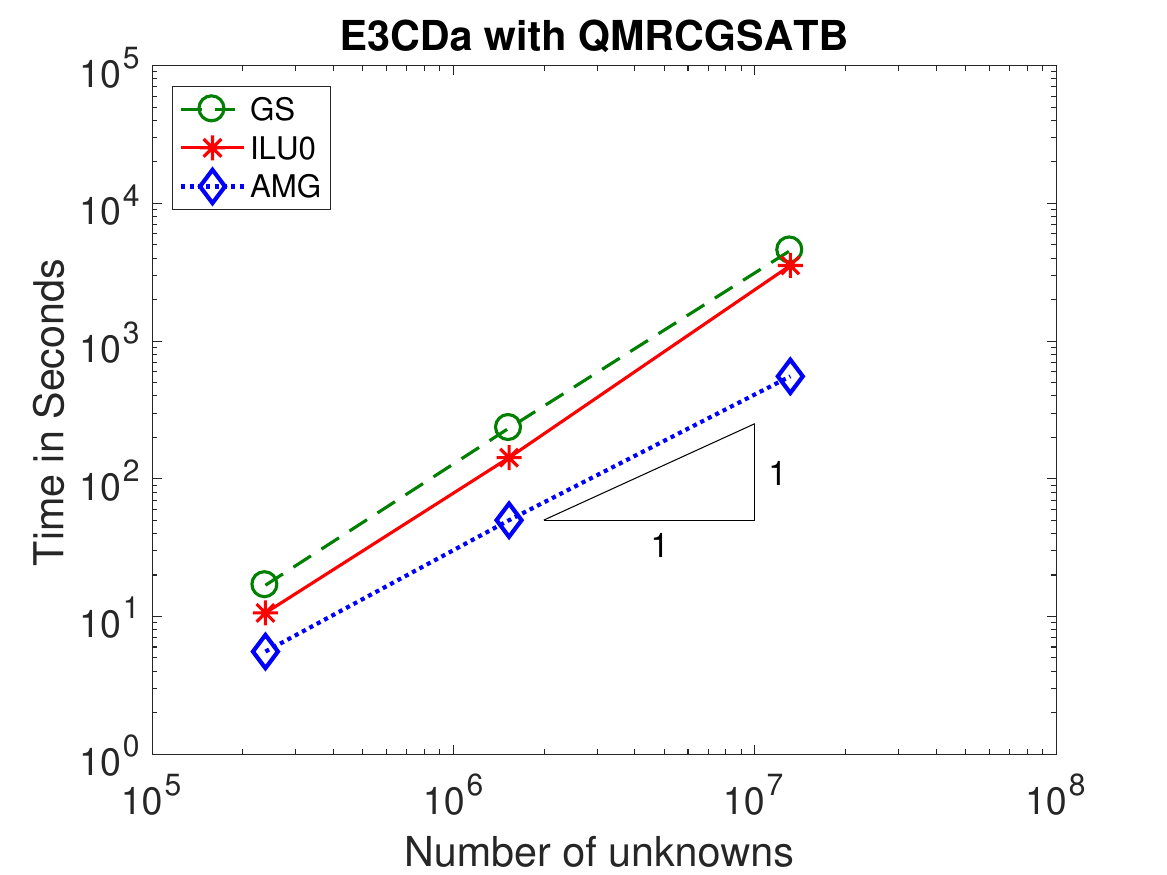}
\par\end{center}%
\end{minipage}\hfill{} %
\begin{minipage}[t]{0.45\textwidth}%
\begin{center}
\includegraphics[width=1\textwidth]{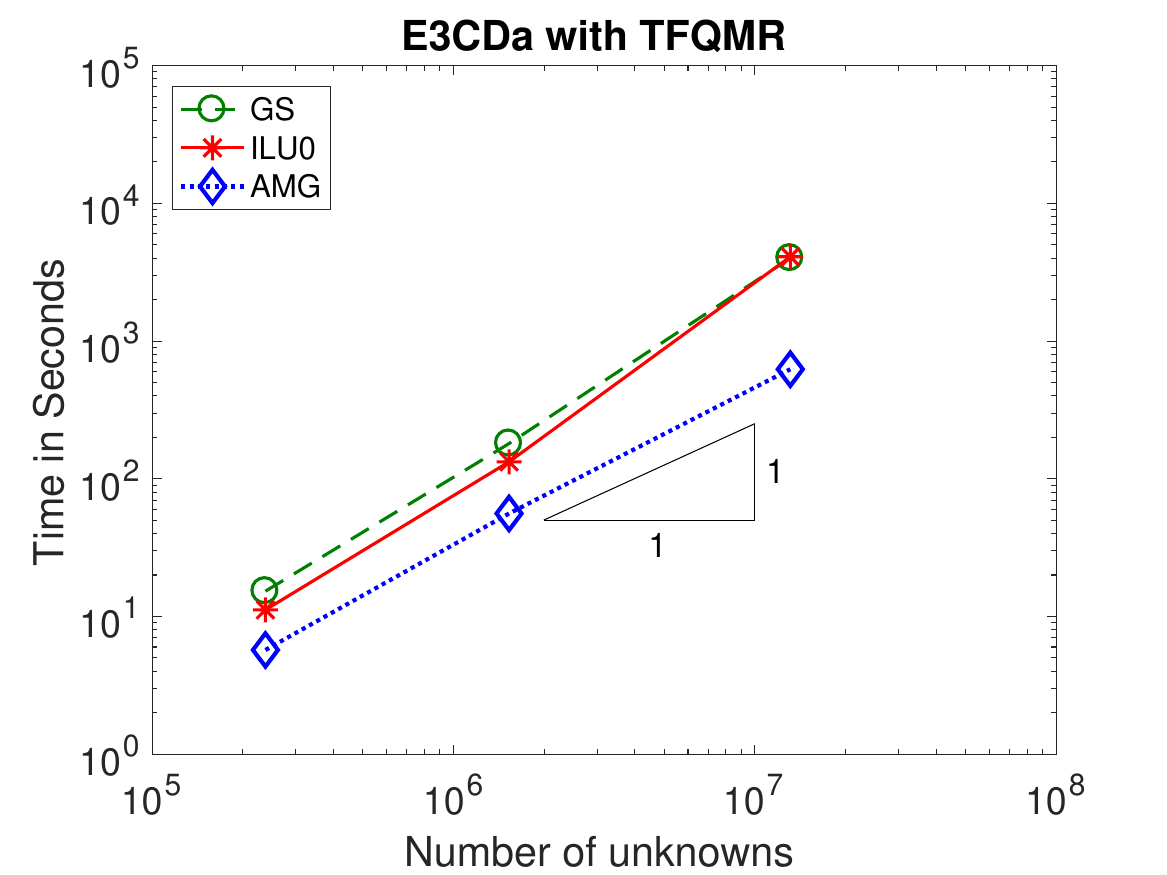}
\par\end{center}%
\end{minipage}

\caption{\label{fig:Scalability}Asymptotic growth of runtimes of the preconditioned
solvers for E3CDa1-3.}
\end{figure}

The relative performance of preconditioned KSP methods may depend
on problem sizes, because different methods may scale differently
with respect to the problem sizes. To assess the asymptotic complexity
of the methods, we consider matrices E3CDa1-3, whose numbers of unknowns
grow approximately by a factor of 8 between each adjacent pair. Figure~\ref{fig:Scalability}
shows the timing results of the four Krylov subspace methods with
Gauss-Seidel, ILU0, and BoomerAMG. The $x$-axis corresponds to the
number of unknowns, and the $y$-axis corresponds to the runtimes,
both in logarithmic scale. For a perfectly scalable method, the slope
should be 1. We observe that with AMG, the slopes for the four KSP
methods are all close to 1. The slopes for Gauss-Seidel and ILU0 were
greater than 1, so the numbers of iterations would grow as the problem
size increases. Therefore, the advantage of multigrid preconditioners
is more significant for larger systems.

\subsection{Comparison of Variants of ILU \label{subsec:ILU-Comparison} }

\begin{figure}[h]
\begin{minipage}[t]{0.45\textwidth}%
\begin{center}
\includegraphics[width=1\textwidth]{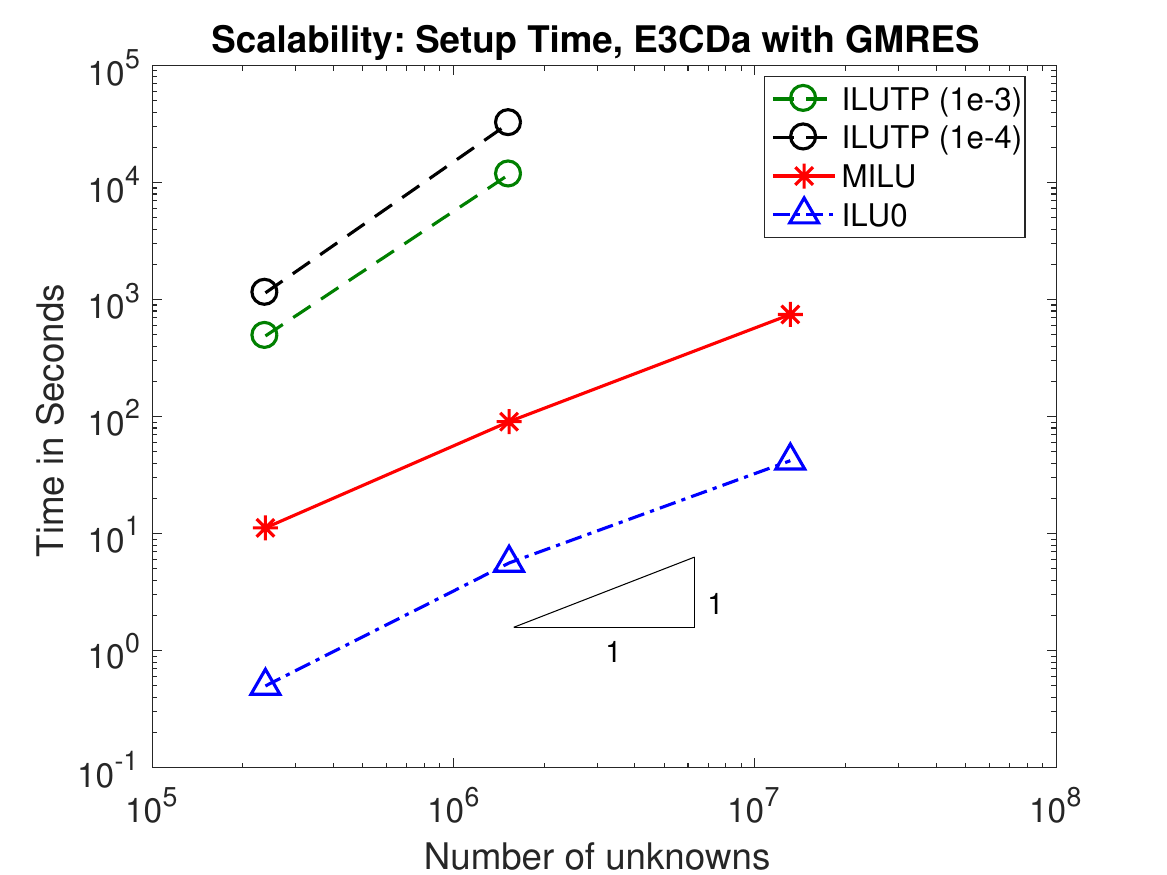}
\par\end{center}%
\end{minipage}\hfill{}%
\begin{minipage}[t]{0.45\textwidth}%
\begin{center}
\includegraphics[width=1\textwidth]{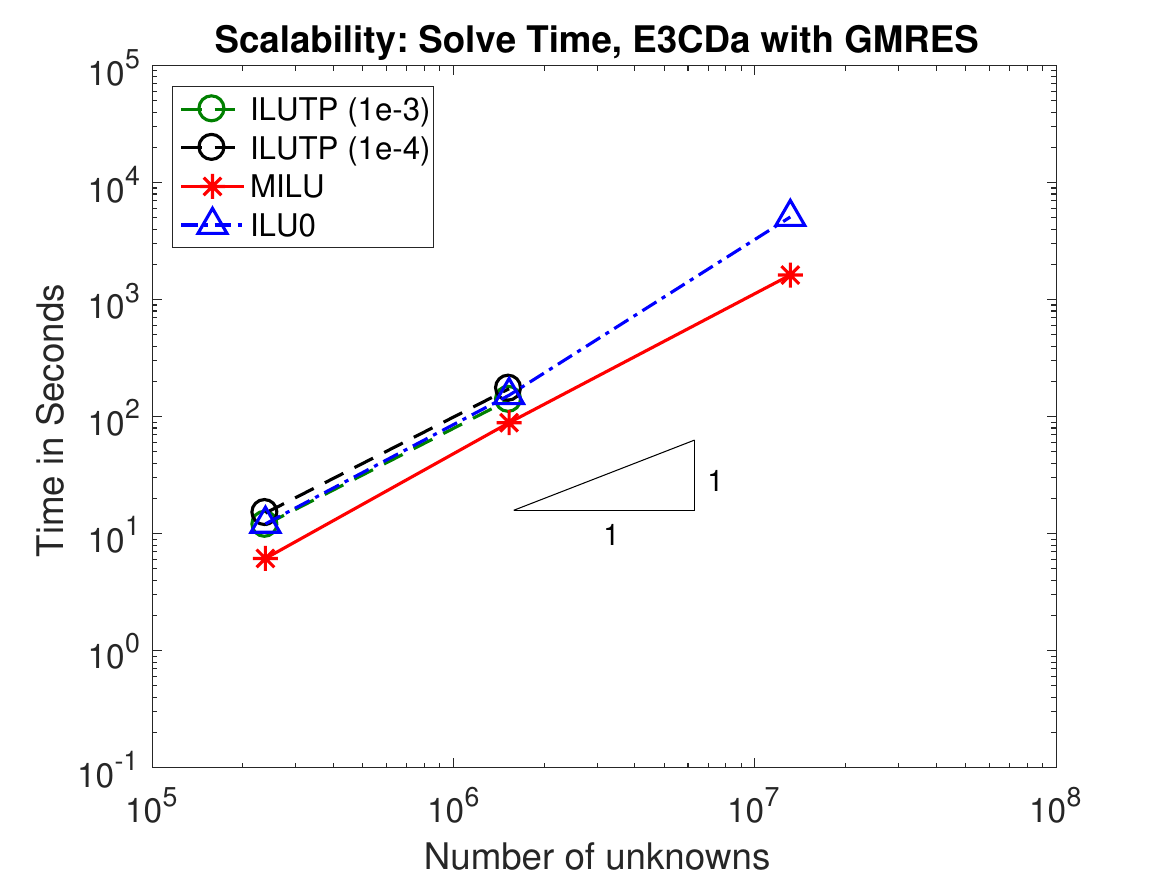}
\par\end{center}%
\end{minipage}

\caption{\label{fig:ilu_scalability}Asymptotic growth of setup (left) and
solve times of ILU0 in PETSc, supernodal ILUTP in SuperLU with drop
tolerance $10^{-4}$ (default) and $10^{-3}$, and MILU in ILUPACK
for E3CDa1-3.}
\end{figure}

We now compare the performance of GMRES with ILU0, ILUTP (in particular,
the supernodal ILUTP implementation in SuperLU v5.2.1), and MILU (in
particular, its implementation in ILUPACK v2.4). We first compare
the scalabilities of their setup and solve times with respect to the
number of unknowns in Figure~\ref{fig:ilu_scalability} for E3CDa1-3.
ILUTP failed for E3CDa3 due to malloc errors in SuperLU, and hence
we only show their results for E3CDa1-2. For ILU0, the setup times
grow linearly, but its solve time is higher than MILU with drop tolerance
$10^{-1}$. For ILUTP, the setup times grow superlinearly with either
the default drop tolerance $10^{-4}$ or a larger drop tolerance $10^{-3}$,
making it orders of magnitude slower than MILU for larger problems.
For MILU, both the setup and solve times grow nearly linearly.

\begin{figure}[h]
\begin{minipage}[b]{0.45\textwidth}%
\begin{center}
\includegraphics[height=2in]{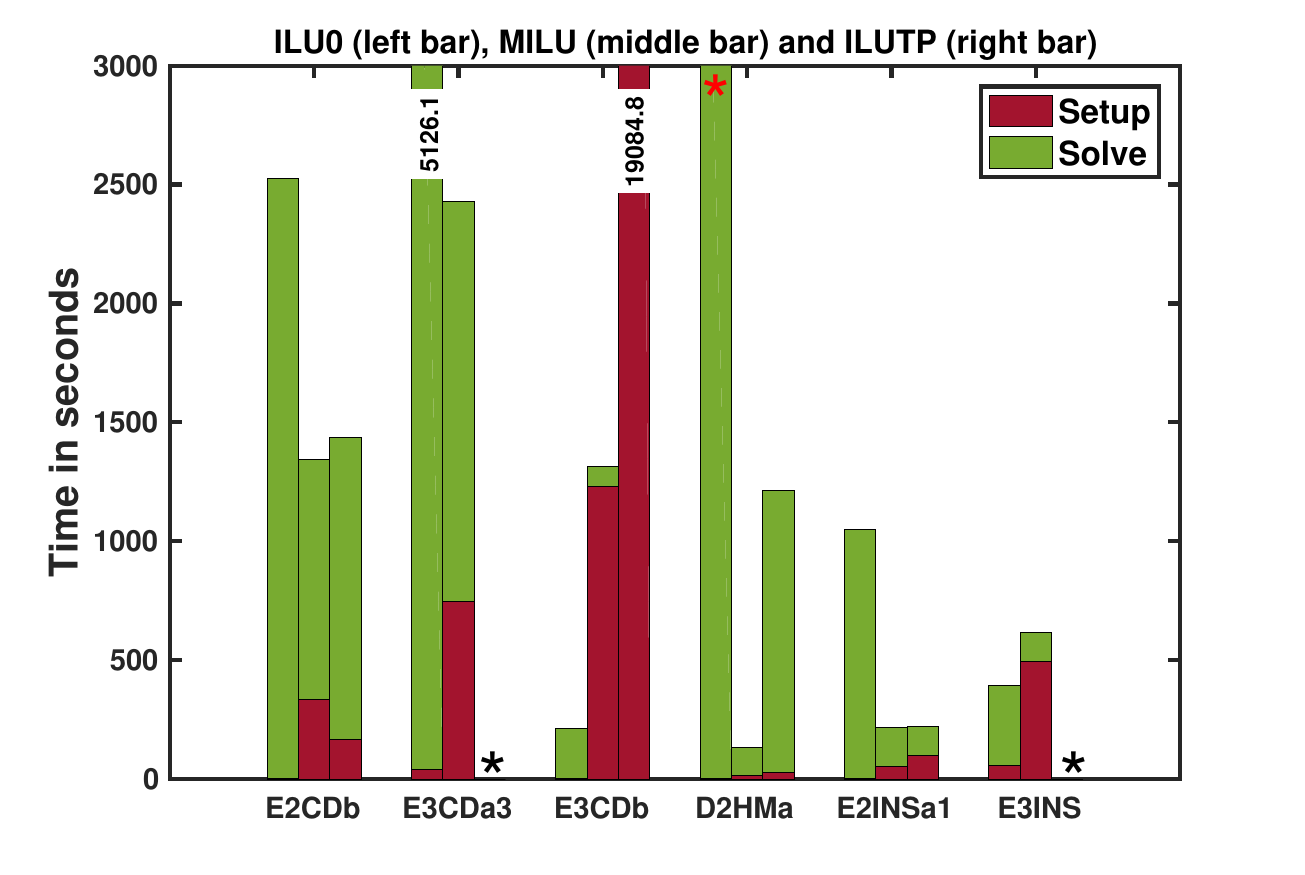}
\par\end{center}
\caption{\label{fig:ilu_timing}Comparison of runtimes of GMRES with ILU0,
MILU, and ILUTP.}
\end{minipage}\hfill{}%
\begin{minipage}[b]{0.45\textwidth}%
\begin{center}
\includegraphics[height=2in]{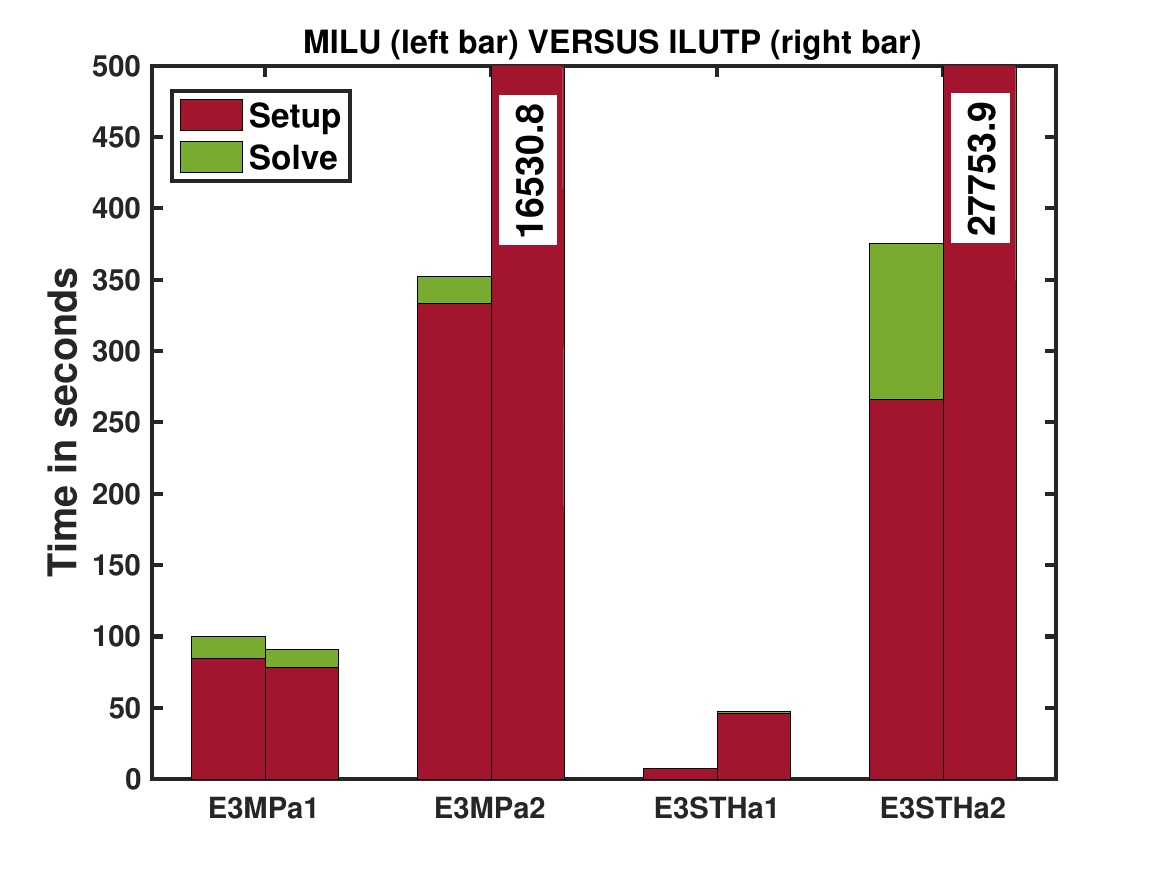}
\par\end{center}
\caption{\label{fig:ilu_saddle_point}Comparison of runtimes of MILU vs. ILUTP
for saddle-point-like problems.}
\end{minipage}
\end{figure}

To compare the overall performance of the methods, Figure\ \ref{fig:ilu_timing}
compares the runtimes of the three preconditioners for six representative
cases. For ILUTP, we used the drop tolerance $10^{-3}$; for MILU,
we used the default parameters. For E3CDa3 and E3NS, `{*}' indicates
failure of ILUTP due to malloc error in SuperLU; for D2HMa, `{*}'
indicates GMRES with ILU0 did not converge after 10,000 iterations.
ILU0 was the best in two out of six cases, due to its low setup cost.
However, it significantly underperformed in the other four cases.
Between MILU and ILUTP, MILU outperformed for all the cases. Note
that in \cite{lishao10}, the opposite conclusion was drawn, where
the test cases had thousands, instead of millions, of unknowns. However,
we note that even though MILU worked well for D2HMa and E3CDa1, it
still underperformed BoomerAMG by nearly two orders of magnitude for
these cases.

The main advantage of MILU, and also of ILUTP to some extent, is that
they are much more robust for ill-conditioned and saddle-point-like
problems. As shown in Table~\ref{tab:results_outline}, MILU was
the only preconditioner that enabled GMRES to solve V3CNS. This highlights
the robustness of MILU for ill-conditioned systems. For saddle-point-like
problems, only MILU allowed GMRES to solve the problems in a reasonable
amount of time, although GMRES with ILUTP also converged. Figure\ \ref{fig:ilu_saddle_point}
compares the setup times and runtimes of GMRES with MILU and ILUTP
for the four saddle-point-like problems. For three out of four problems,
MILU significantly outperformed ILUTP. Hence, MILU is the most robust
choice for saddle-point-like problems.

\subsection{BoomerAMG Versus ML \label{subsec:HYPRE-vs-ML}}

\begin{figure}[h]
\begin{minipage}[t]{0.45\textwidth}%
\begin{center}
\includegraphics[width=1\textwidth]{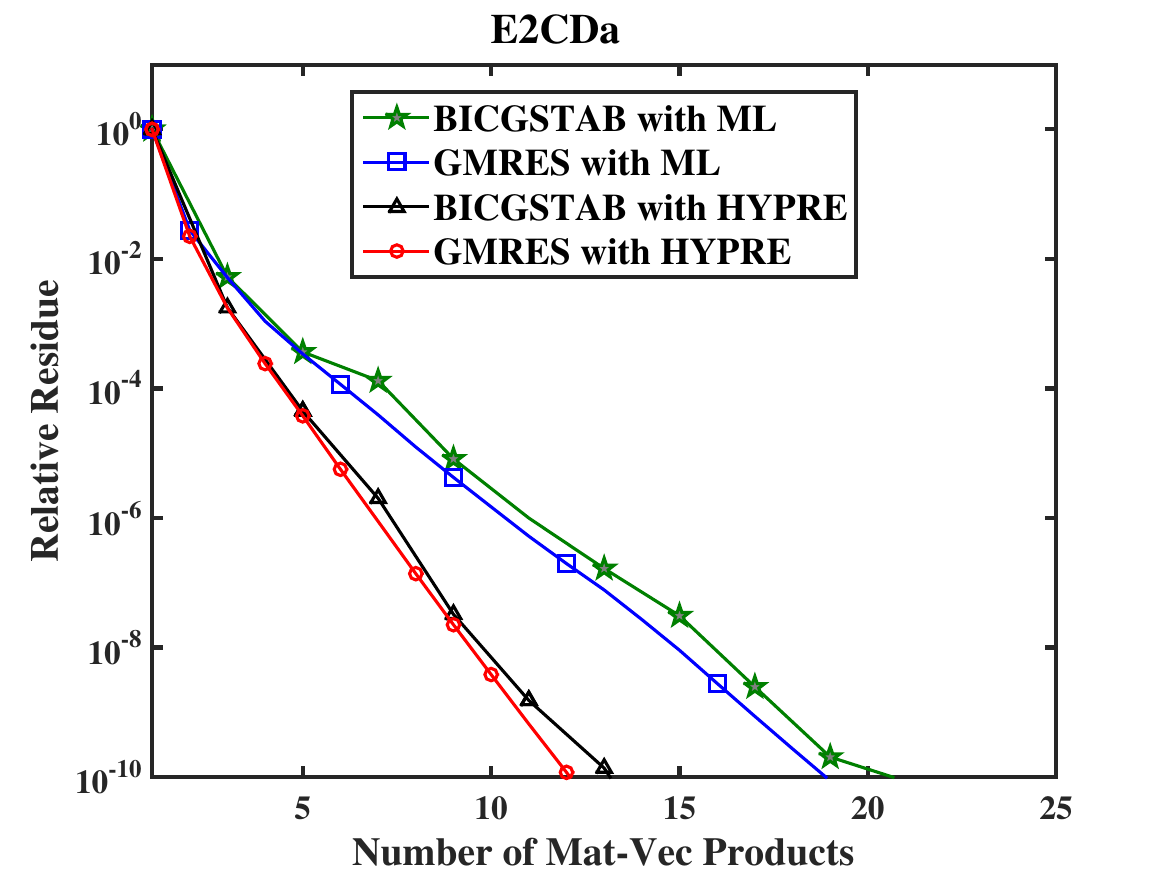}
\par\end{center}%
\end{minipage}\hfill{} %
\begin{minipage}[t]{0.45\textwidth}%
\begin{center}
\includegraphics[width=1\textwidth]{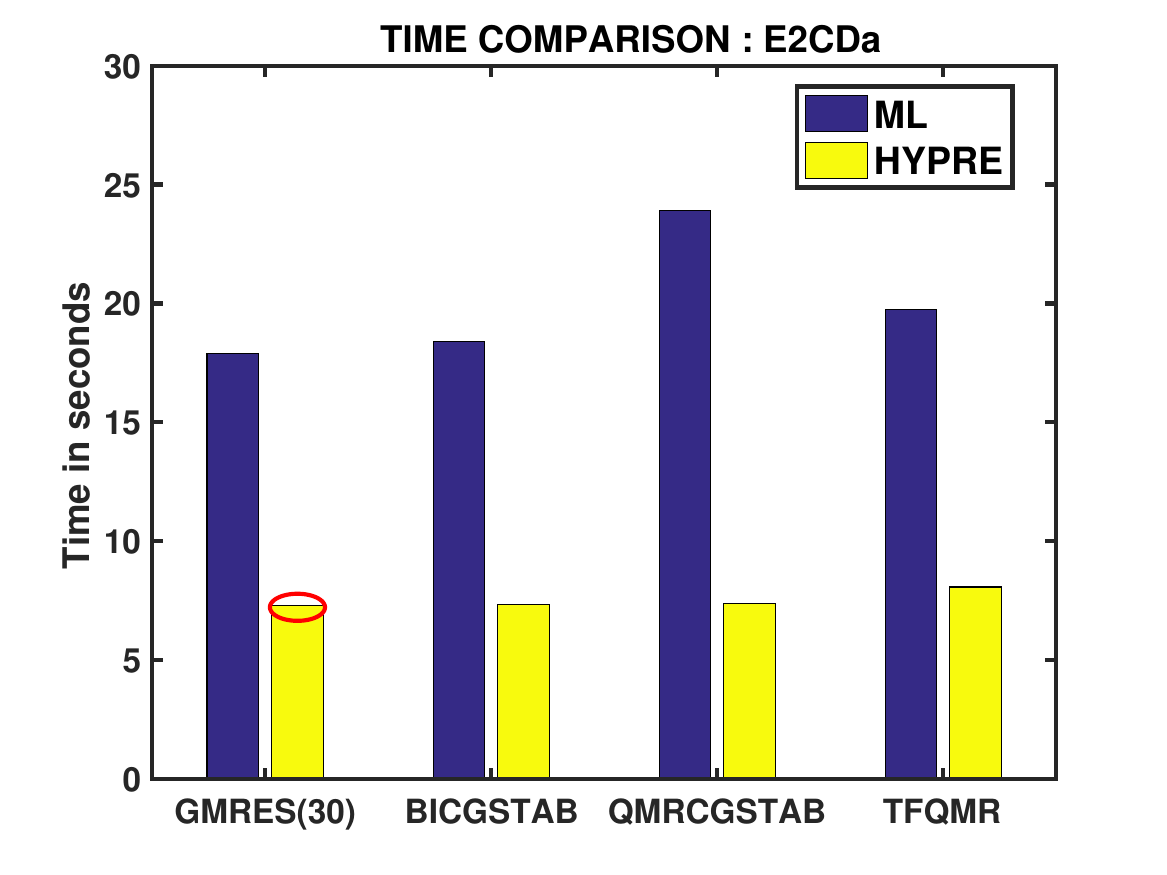}
\par\end{center}%
\end{minipage}

\begin{minipage}[t]{0.45\textwidth}%
\begin{center}
\includegraphics[width=1\textwidth]{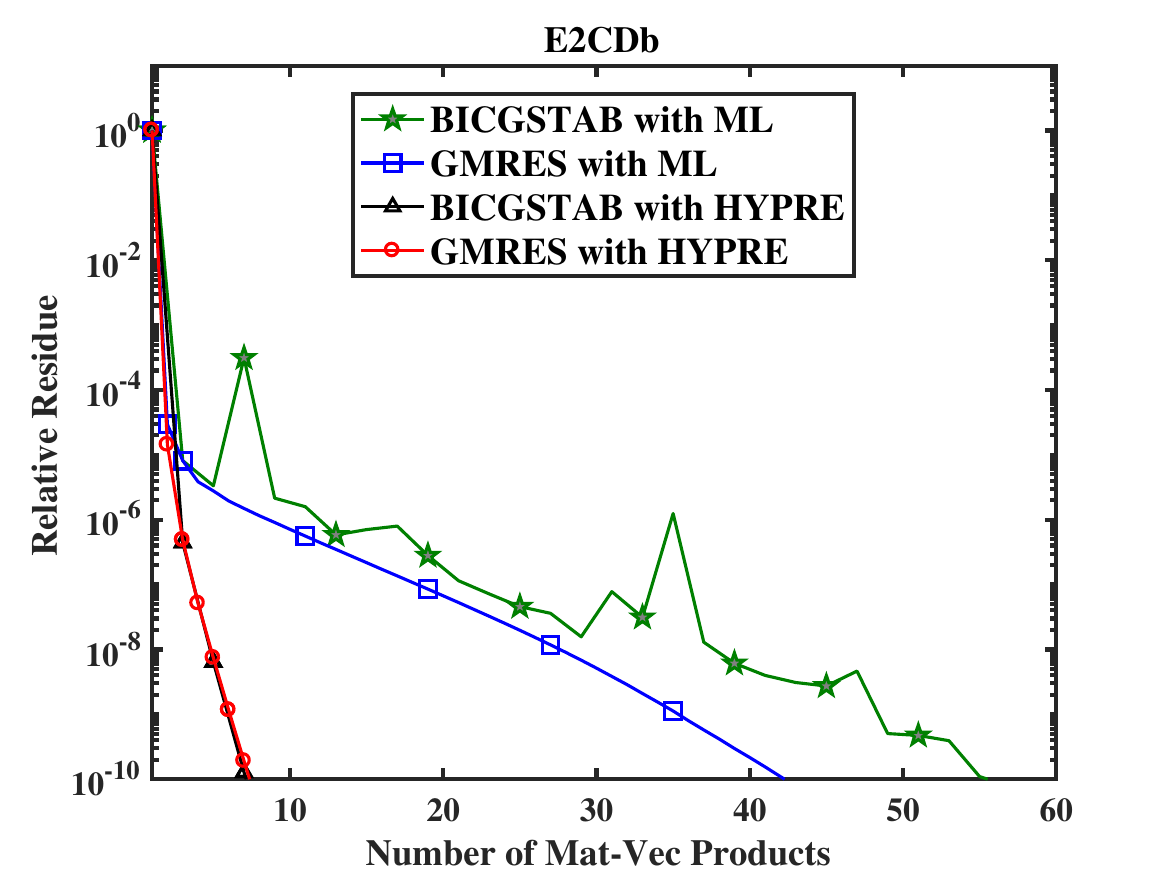}
\par\end{center}%
\end{minipage}\hfill{} %
\begin{minipage}[t]{0.45\textwidth}%
\begin{center}
\includegraphics[width=1\textwidth]{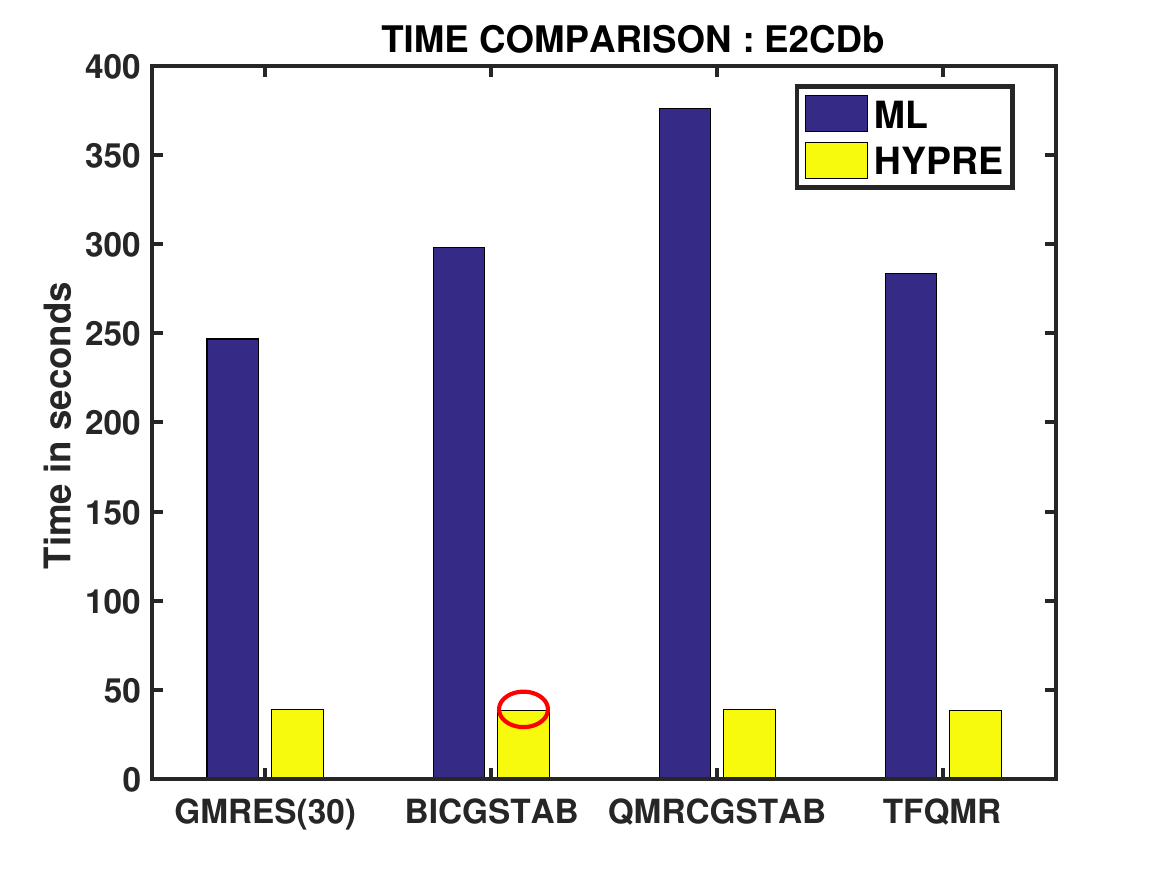}
\par\end{center}%
\end{minipage}

\begin{minipage}[t]{0.45\textwidth}%
\begin{center}
\includegraphics[width=1\textwidth]{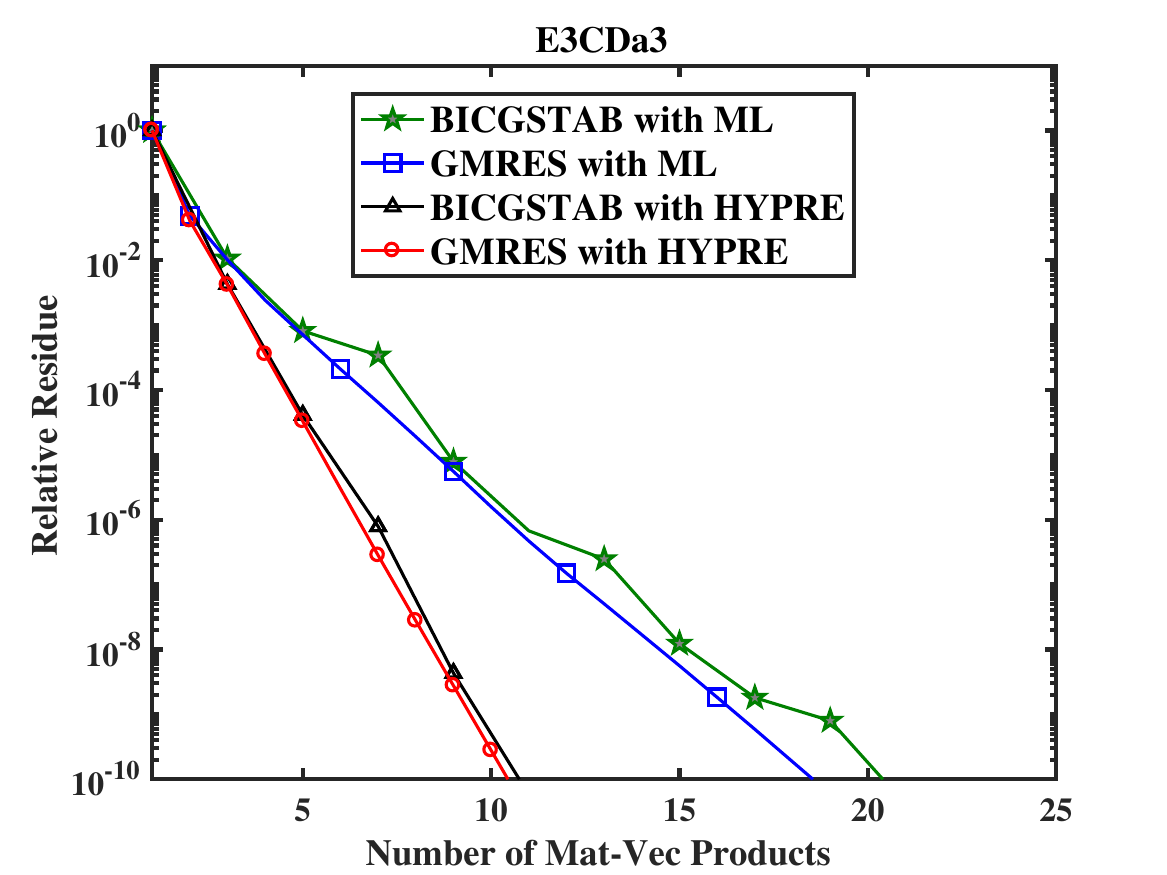}
\par\end{center}%
\end{minipage}\hfill{} %
\begin{minipage}[t]{0.45\textwidth}%
\begin{center}
\includegraphics[width=1\textwidth]{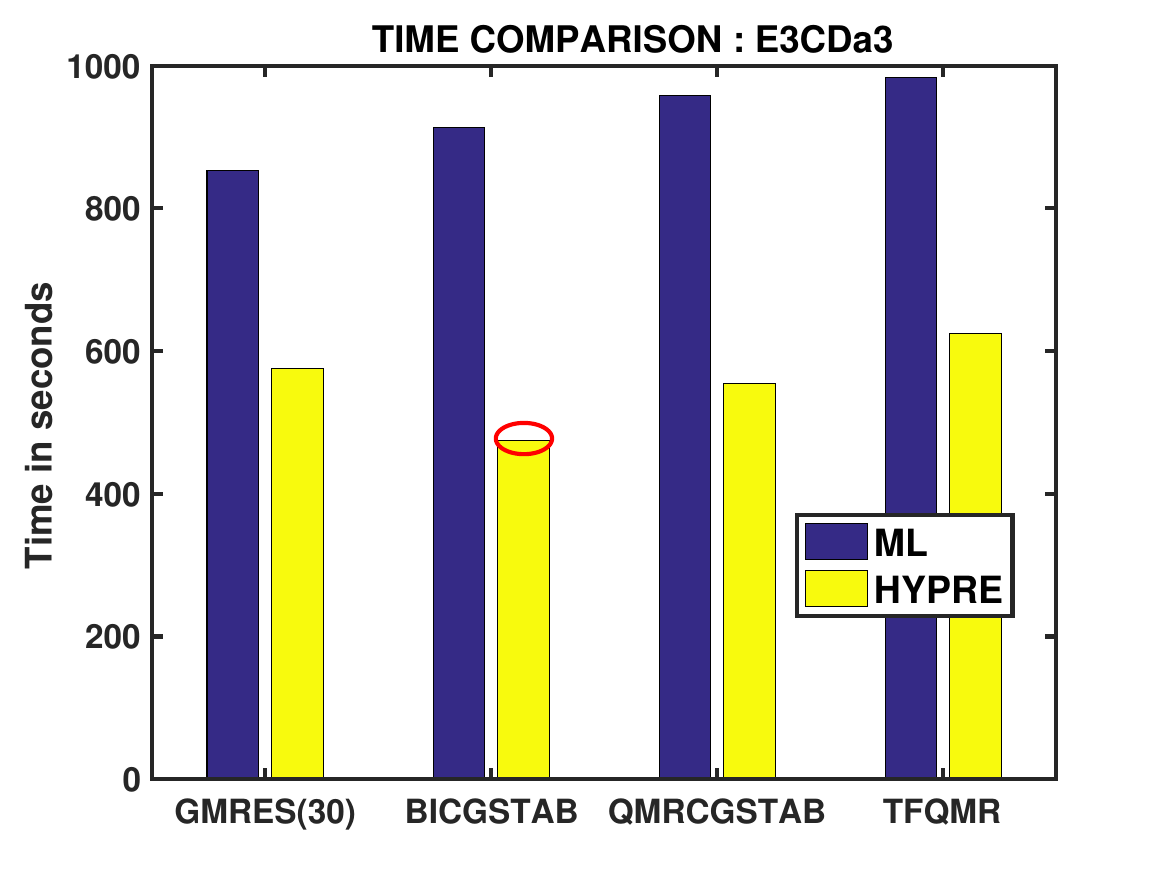}
\par\end{center}%
\end{minipage}

\begin{minipage}[t]{0.45\textwidth}%
\begin{center}
\includegraphics[width=1\textwidth]{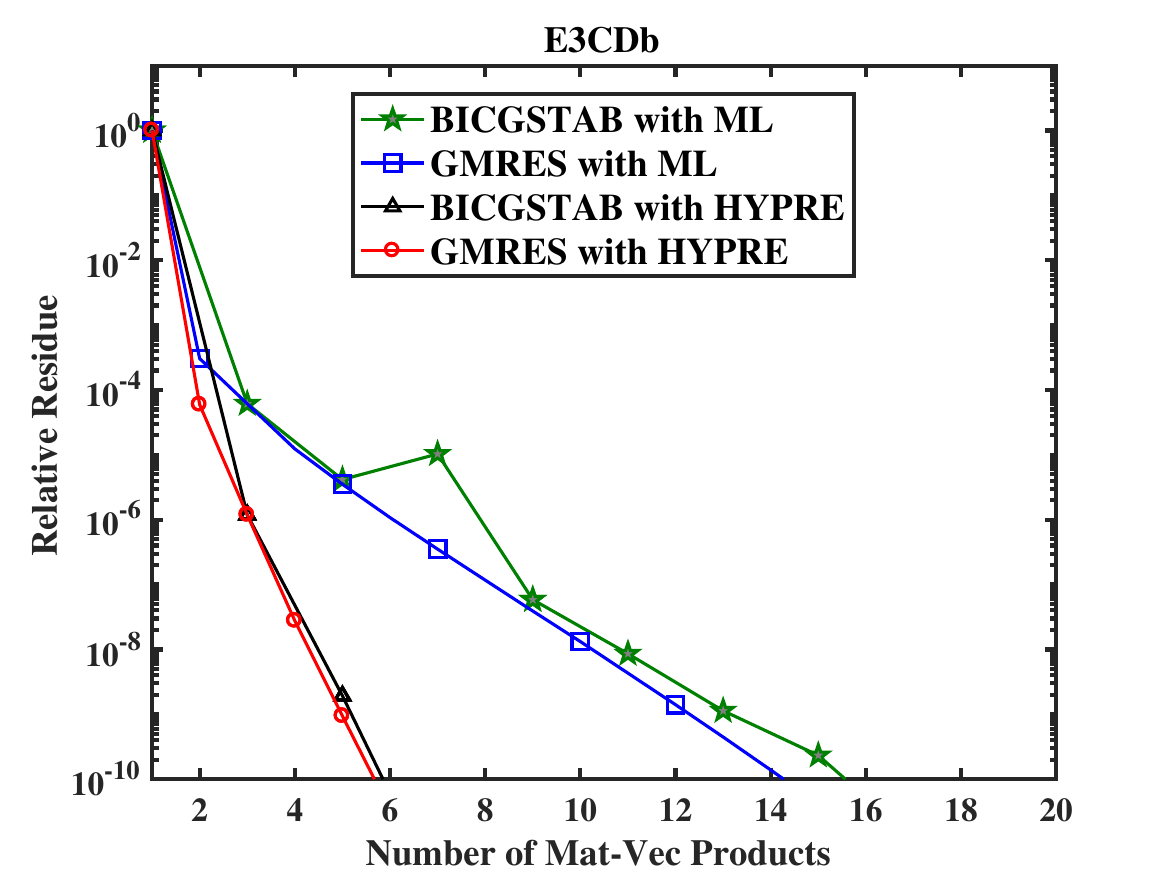}
\par\end{center}%
\end{minipage}\hfill{} %
\begin{minipage}[t]{0.45\textwidth}%
\begin{center}
\includegraphics[width=1\textwidth]{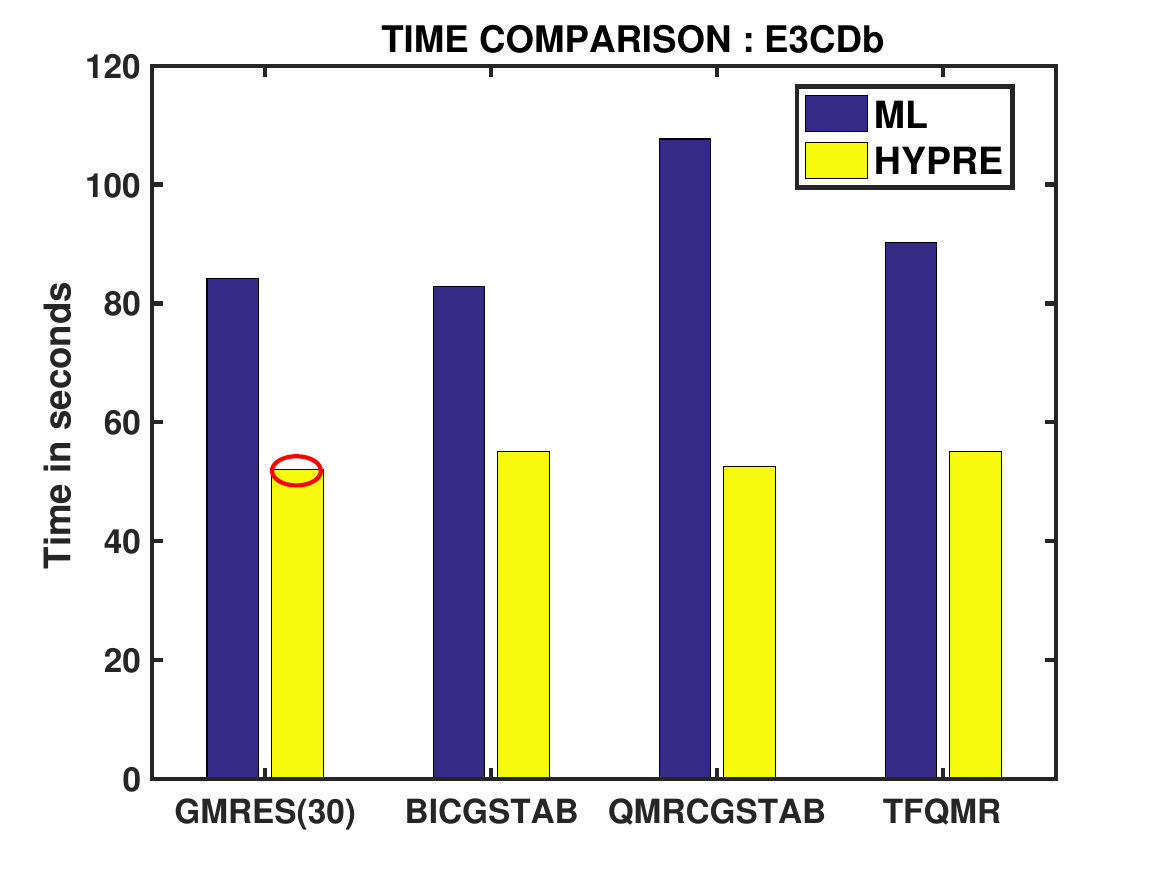}
\par\end{center}%
\end{minipage}

\caption{\label{fig:Hypre_vs_ML}Convergence history (left) and runtimes (right)
of GMRES and BiCGSTAB with BoomerAMG and ML. Circled bars indicate
best performance.}
\end{figure}

\begin{figure}[h]
\begin{minipage}[t]{0.45\textwidth}%
\begin{center}
\includegraphics[width=1\textwidth]{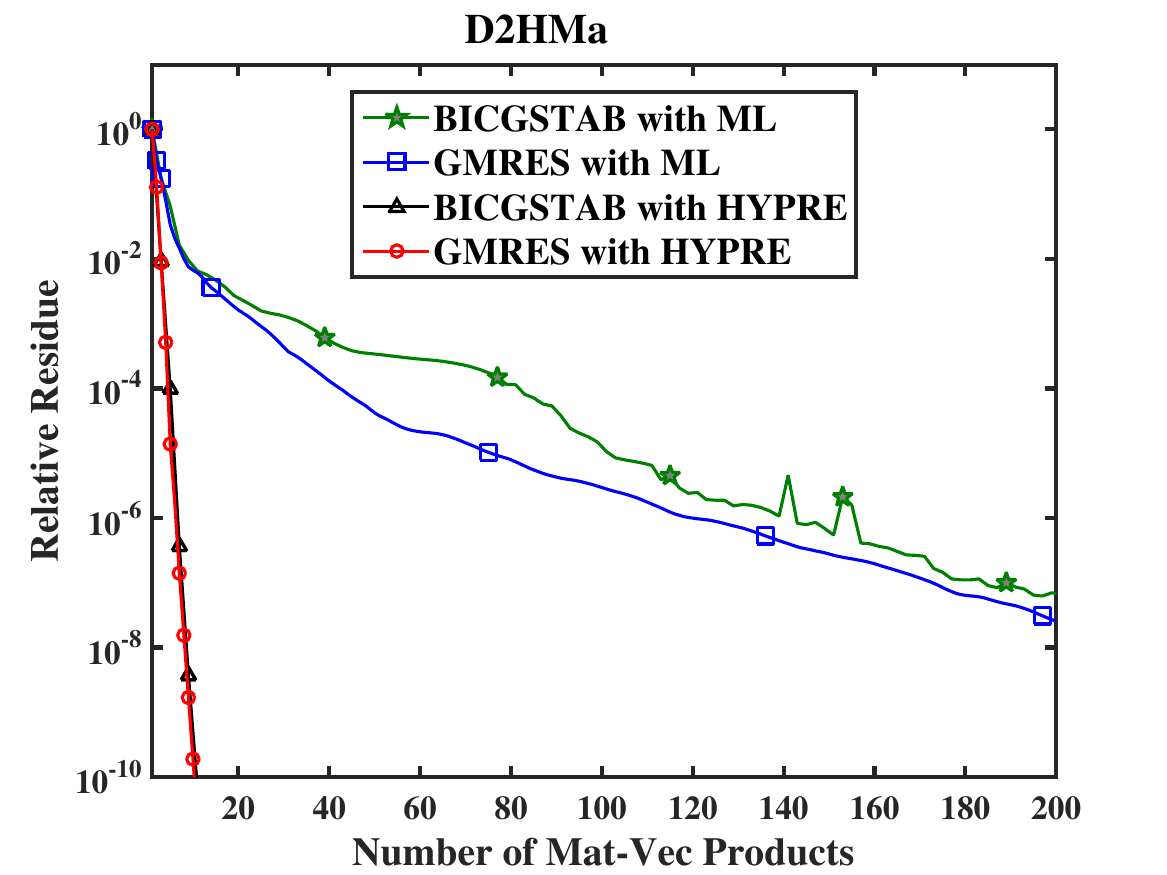}
\par\end{center}%
\end{minipage}\hfill{} %
\begin{minipage}[t]{0.45\textwidth}%
\begin{center}
\includegraphics[width=1\textwidth]{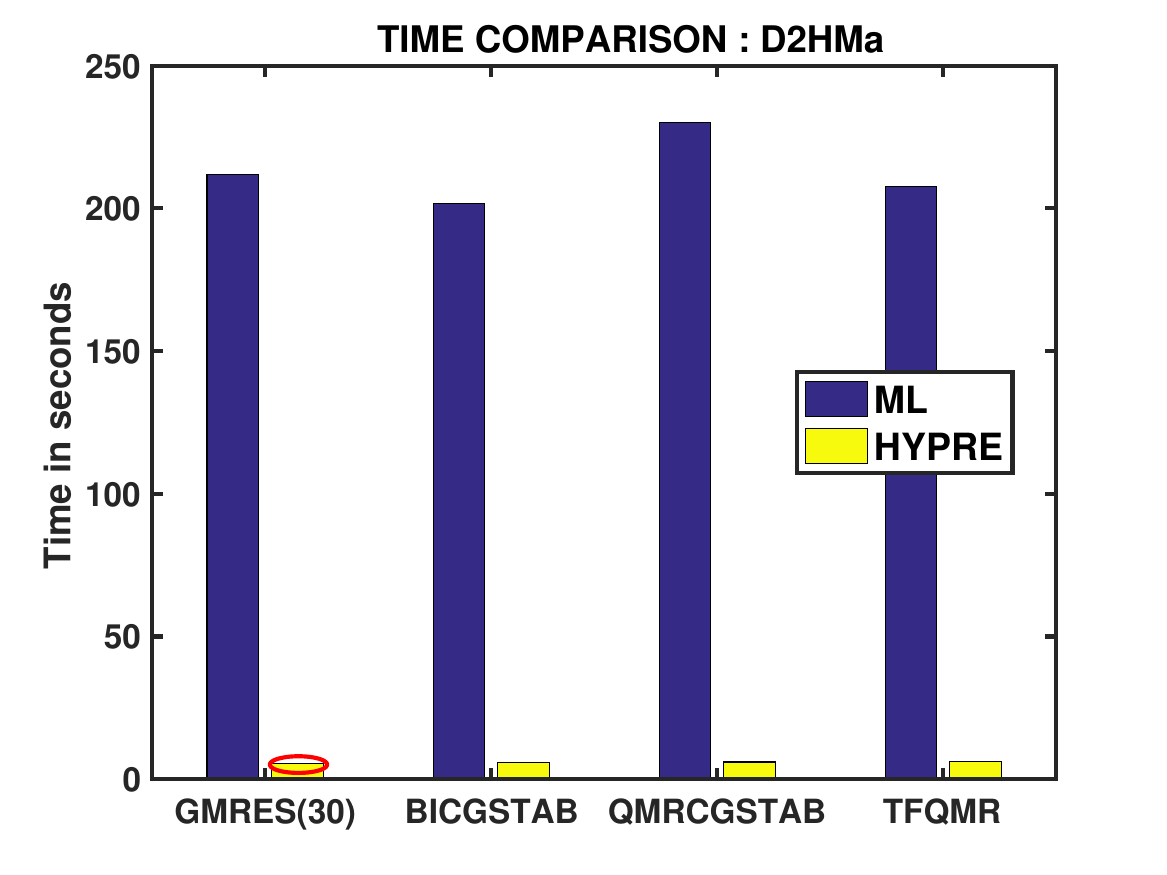}
\par\end{center}%
\end{minipage}

\begin{minipage}[t]{0.45\textwidth}%
\begin{center}
\includegraphics[width=1\textwidth]{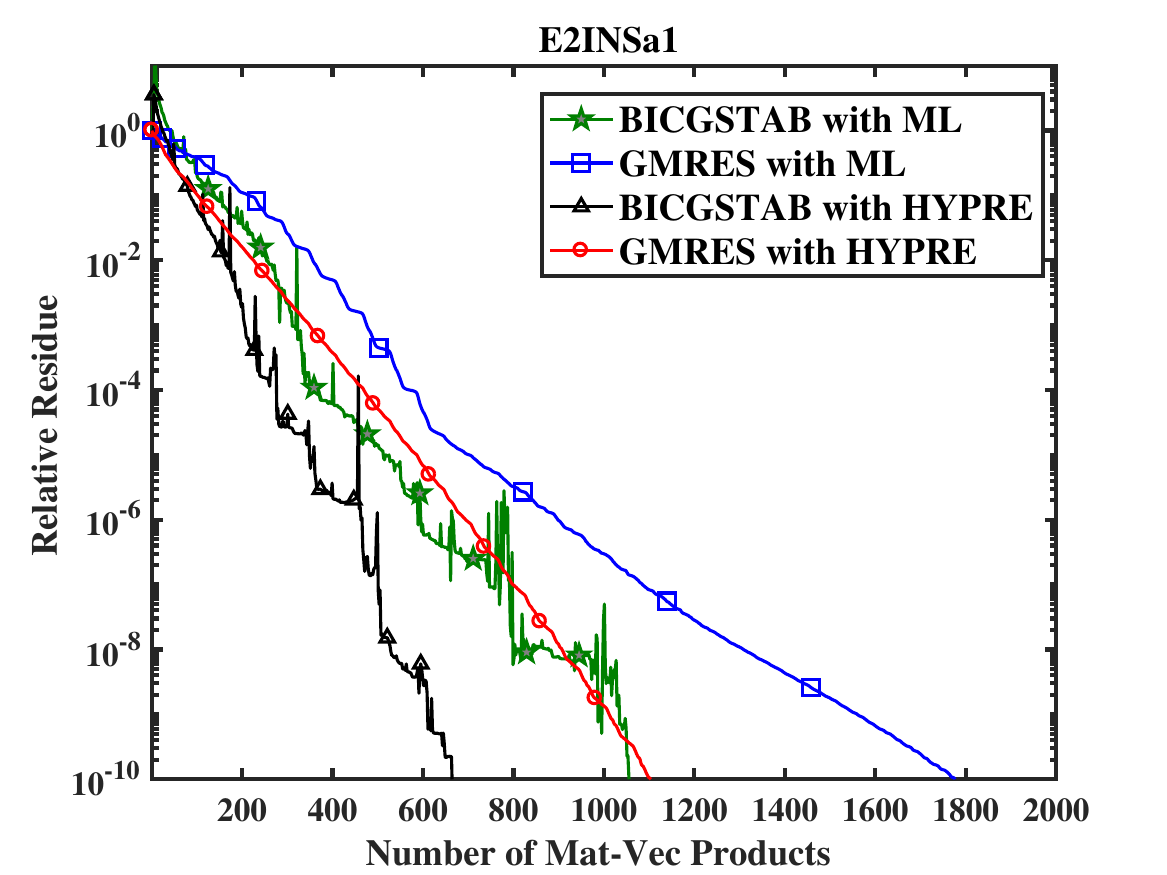}
\par\end{center}%
\end{minipage}\hfill{} %
\begin{minipage}[t]{0.45\textwidth}%
\begin{center}
\includegraphics[width=1\textwidth]{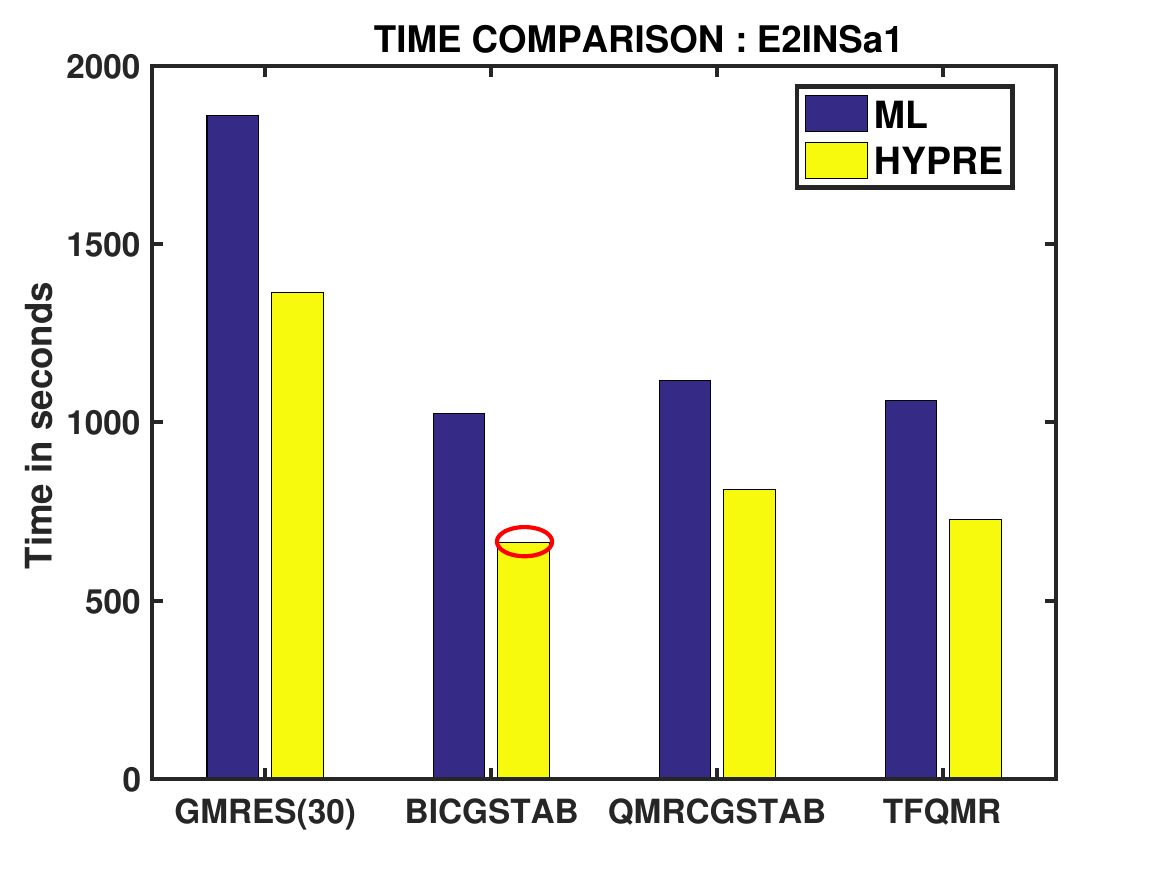}
\par\end{center}%
\end{minipage}

\begin{minipage}[t]{0.45\textwidth}%
\begin{center}
\includegraphics[width=1\textwidth]{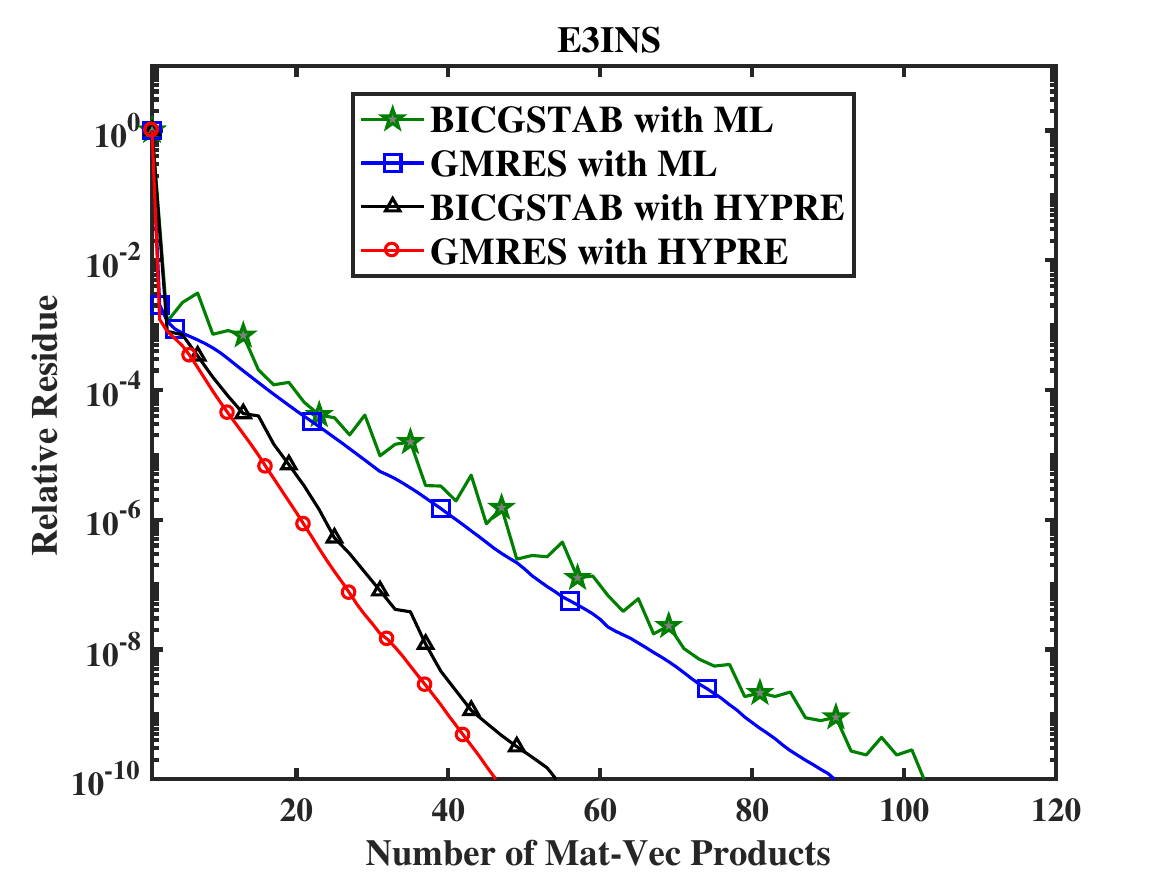}
\par\end{center}%
\end{minipage}\hfill{} %
\begin{minipage}[t]{0.45\textwidth}%
\begin{center}
\includegraphics[width=1\textwidth]{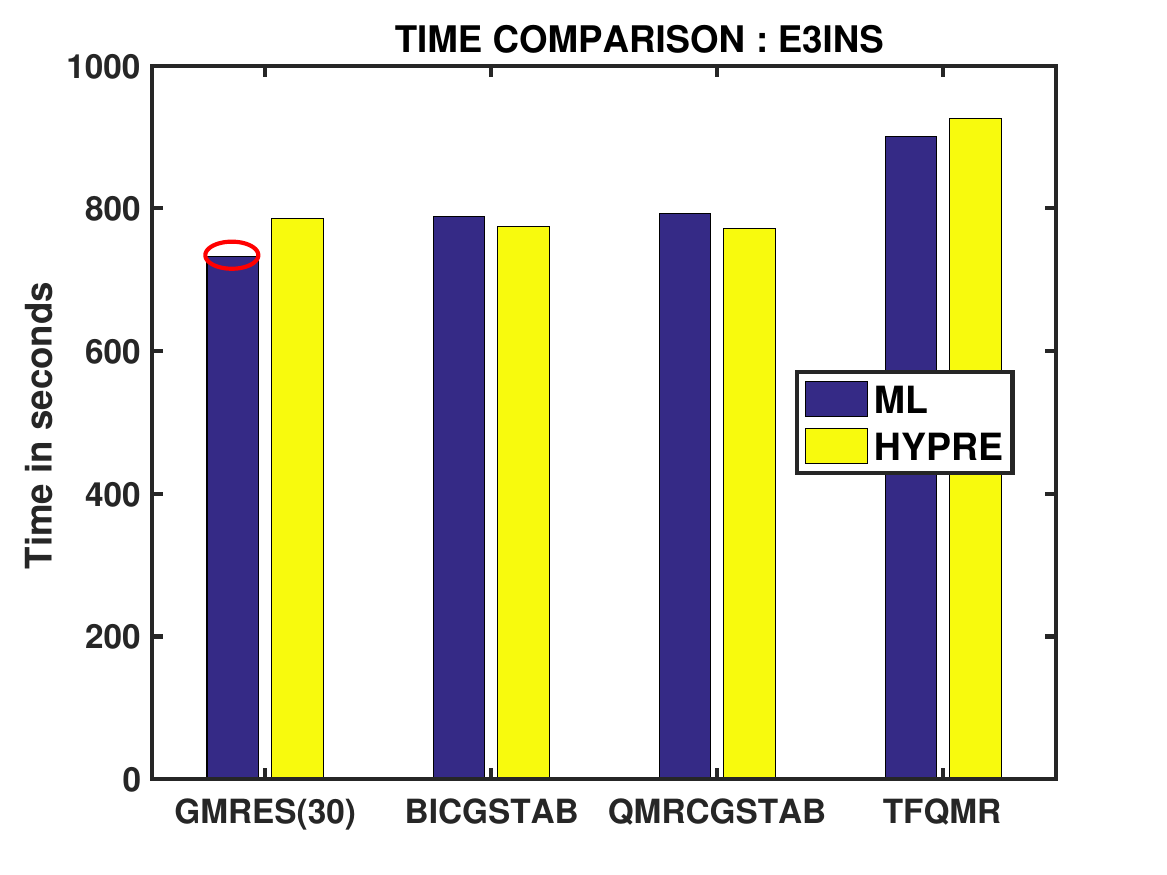}
\par\end{center}%
\end{minipage}

\caption{\label{fig:Hypre_vs_ML-1}Convergence history (left) and runtimes
(right) of GMRES and BiCGSTAB with BoomerAMG and ML. Circled bars
indicate best performance.}
\end{figure}

Finally, we compare two AMG preconditioners: BoomerAMG and ML, which
are variants of classical AMG and smoothed aggregation, respectively.
Both of these methods scale nearly linearly in terms of the number
of unknowns. However, their actual performance may differ drastically.
Figure~\ref{fig:Hypre_vs_ML} and \ref{fig:Hypre_vs_ML-1} compare
the convergence and runtimes of GMRES and BiCGSTAB with BoomerAMG
and ML for seven representative results. We used HMIS coarsening with
FF1 interpolation for BoomerAMG and the default parameters for ML.
For six out of seven cases, BoomerAMG outperformed ML. The only case
that ML outperformed BoomerAMG was E3INS, for which Gauss-Seidel also
outperforms ML. For D2HMa, which is ill-conditioned, BoomerAMG was
more than 20 times faster than ML. This is probably because for smoothed
aggregation, too many aggregates can lead to the growth of complexity
and irregularly shaped aggregates, which can cause poor convergence
\cite{yang2006parallel}. We note that we have tried many different
options in ML, and the results were qualitatively the same.

\begin{figure}[h]
\begin{minipage}[t]{0.45\textwidth}%
\begin{center}
\includegraphics[width=1\textwidth]{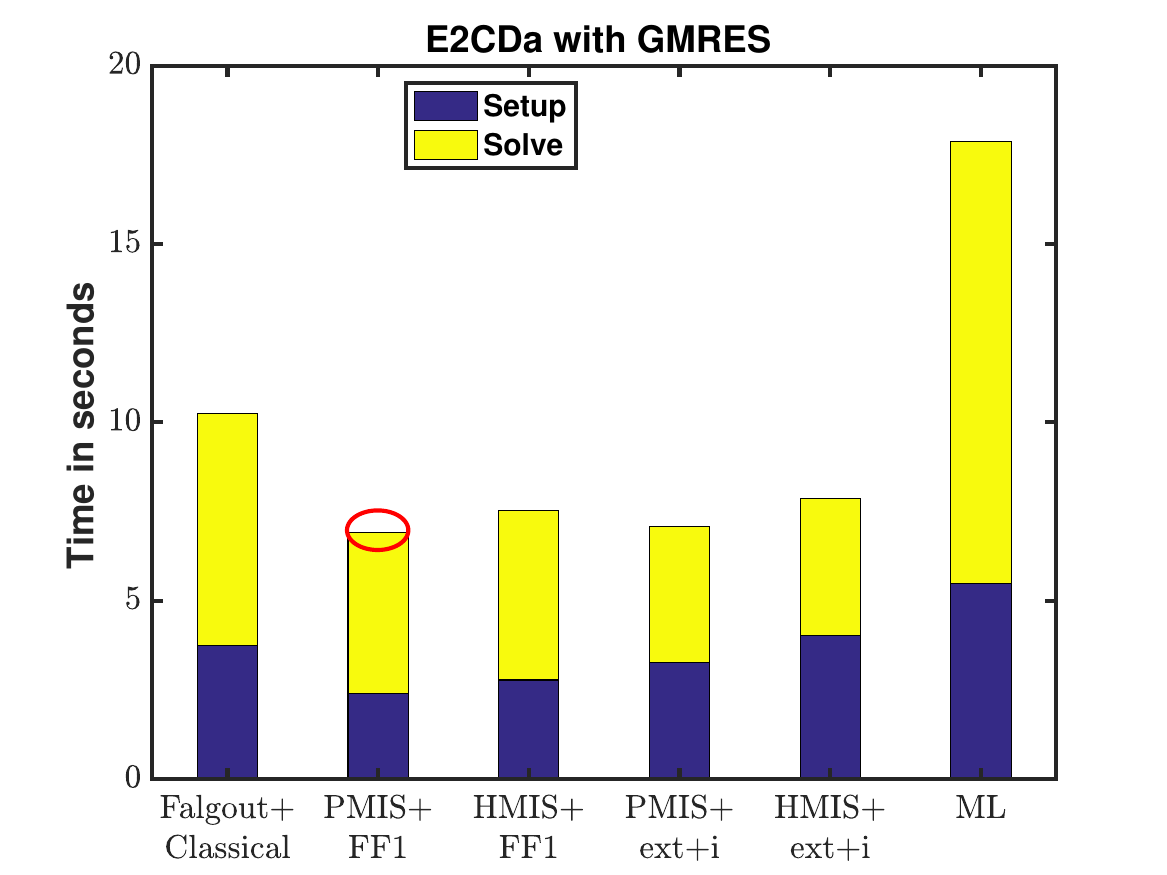}
\par\end{center}%
\end{minipage}\hfill{} %
\begin{minipage}[t]{0.45\textwidth}%
\begin{center}
\includegraphics[width=1\textwidth]{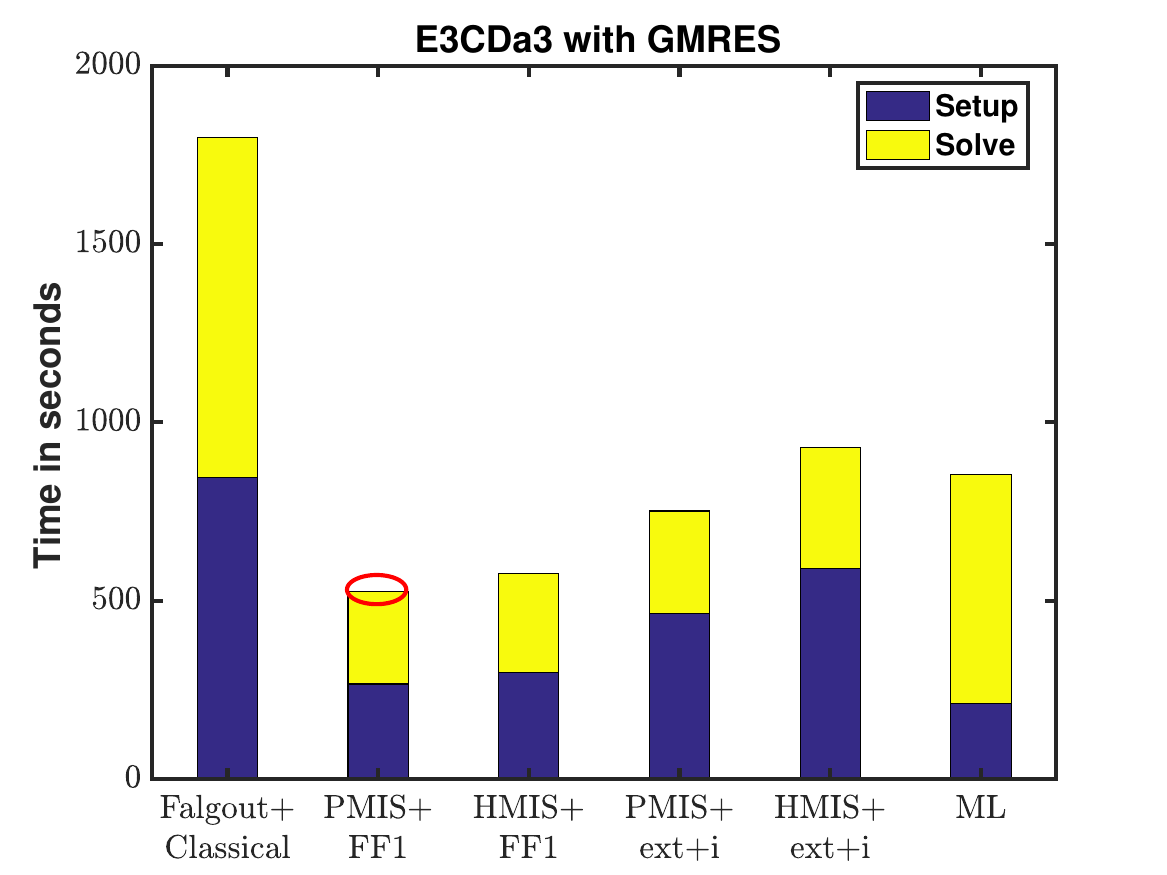}
\par\end{center}%
\end{minipage}

\begin{minipage}[t]{0.45\textwidth}%
\begin{center}
\includegraphics[width=1\textwidth]{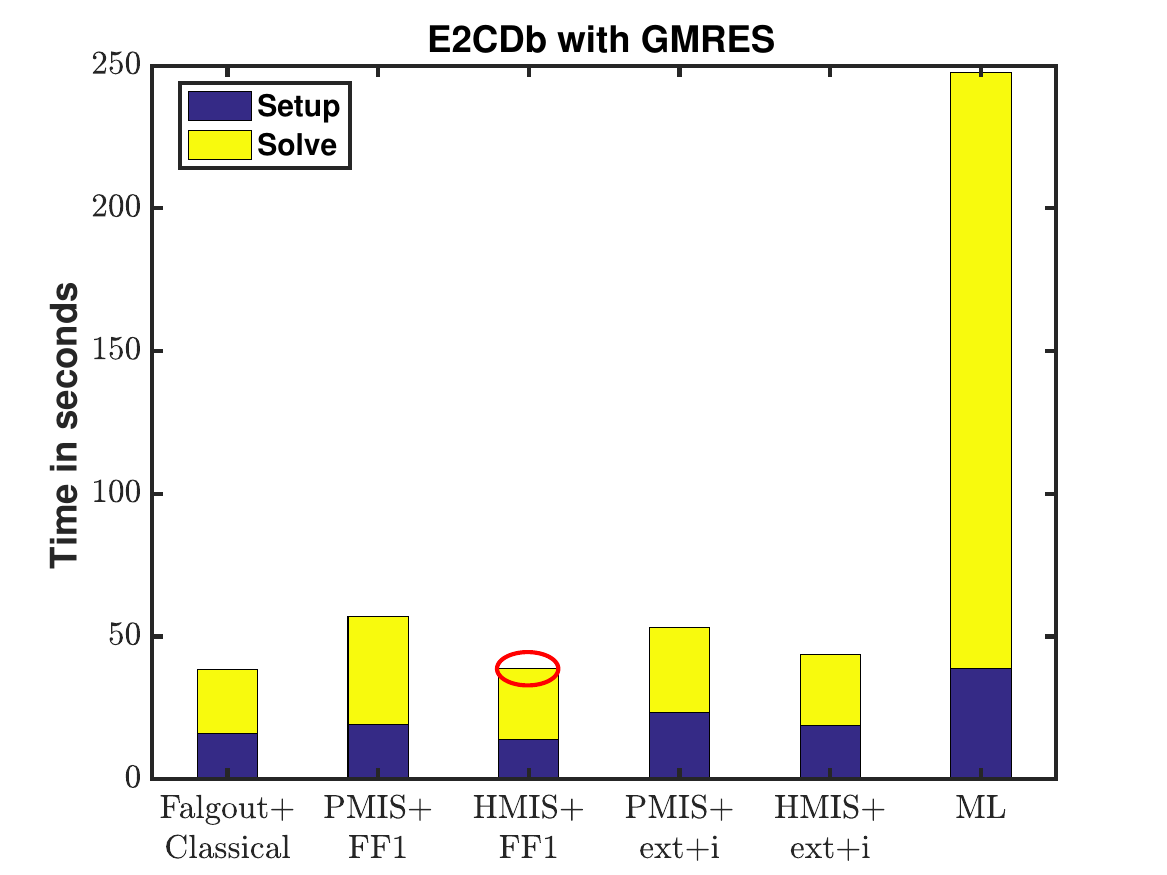}
\par\end{center}%
\end{minipage}\hfill{} %
\begin{minipage}[t]{0.45\textwidth}%
\begin{center}
\includegraphics[width=1\textwidth]{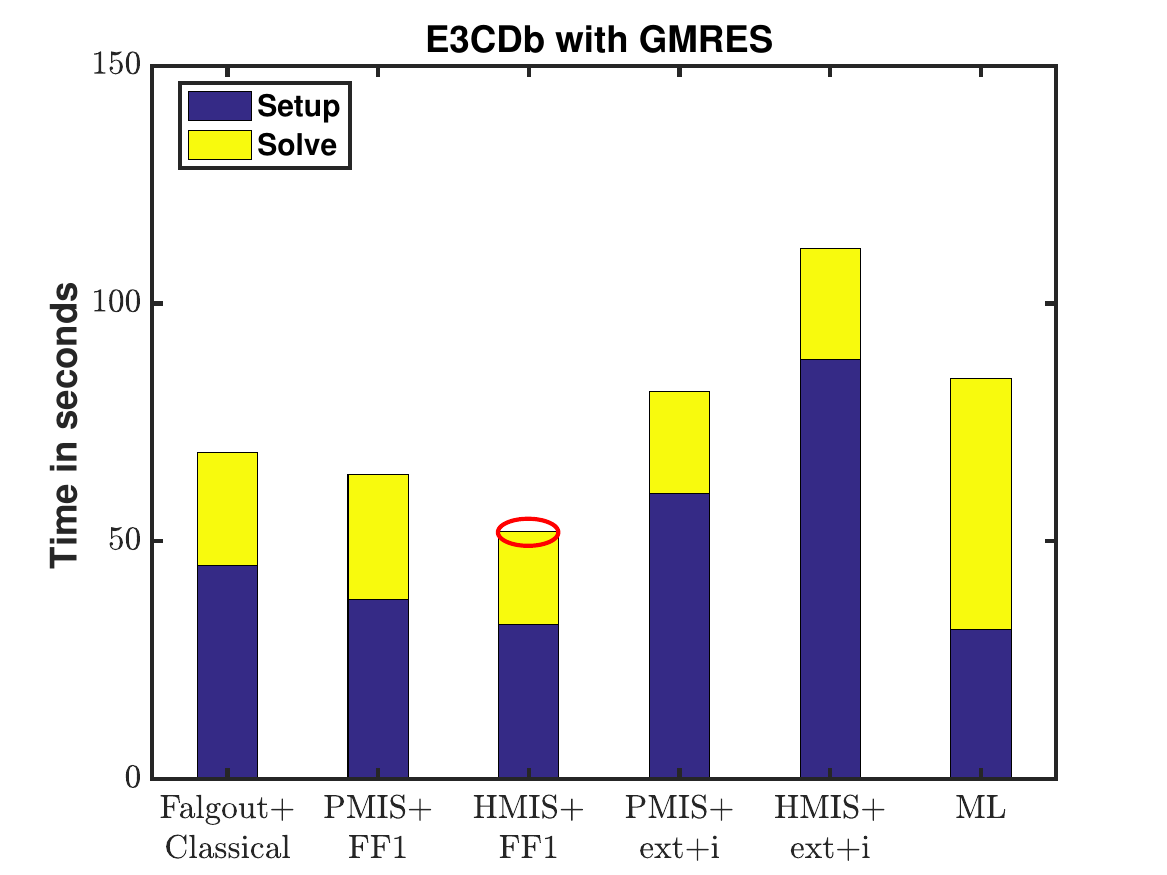}
\par\end{center}%
\end{minipage}

\begin{minipage}[t]{0.45\textwidth}%
\begin{center}
\includegraphics[width=1\textwidth]{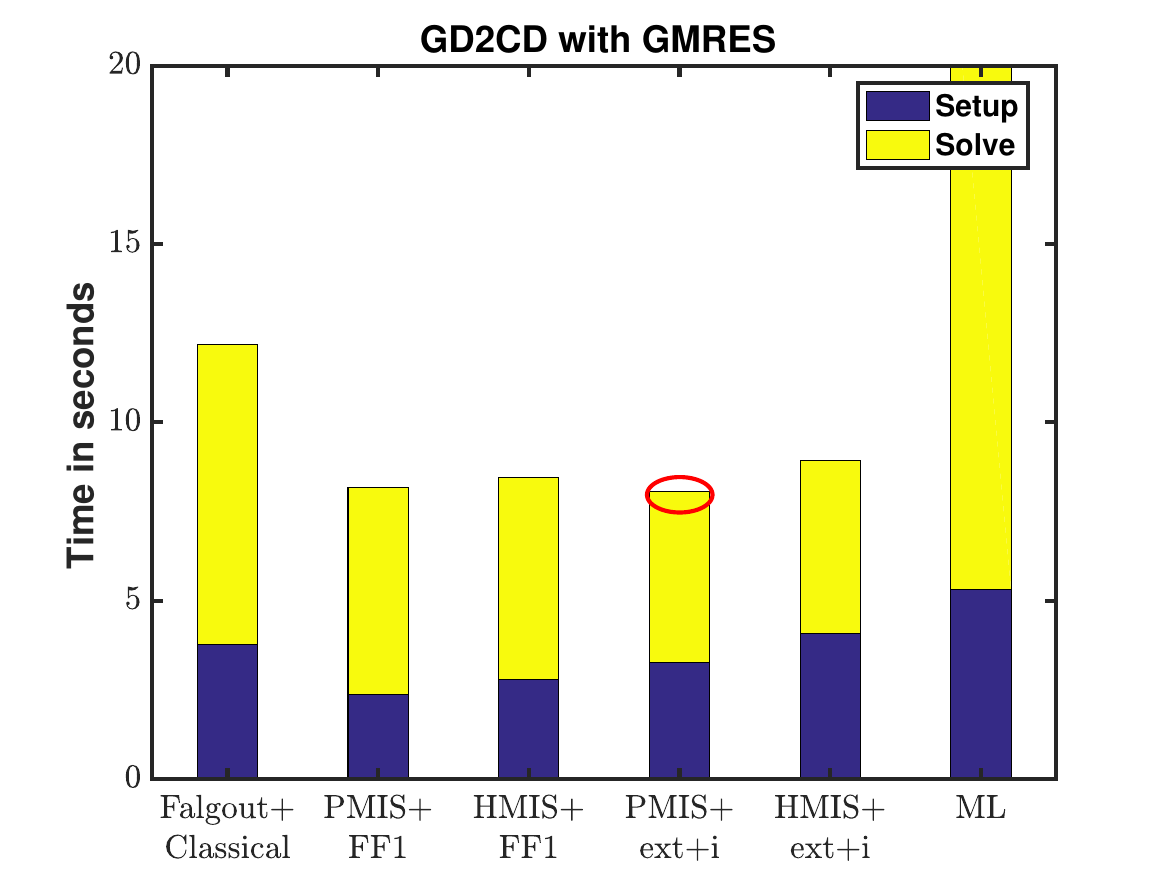}
\par\end{center}%
\end{minipage}\hfill{} %
\begin{minipage}[t]{0.45\textwidth}%
\begin{center}
\includegraphics[width=1\textwidth]{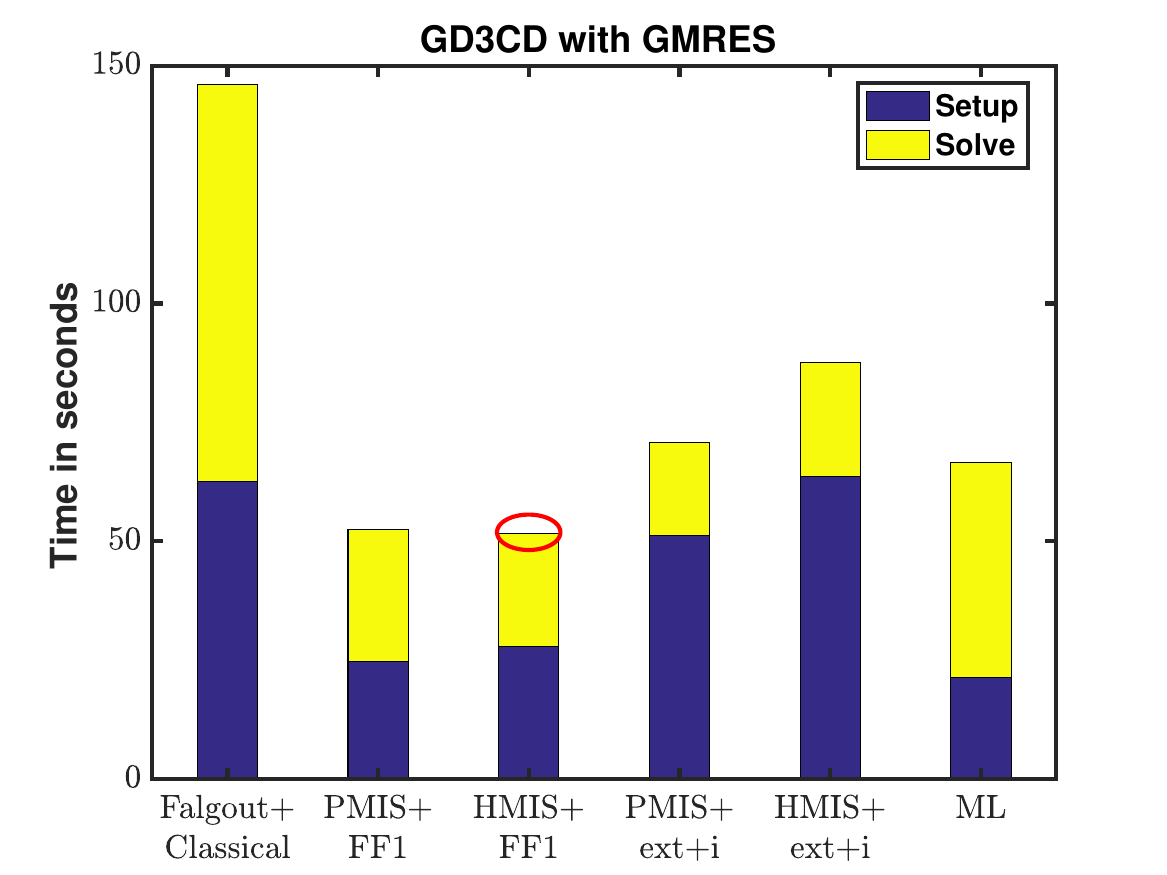}
\par\end{center}%
\end{minipage}

\begin{minipage}[t]{0.45\textwidth}%
\begin{center}
\includegraphics[width=1\textwidth]{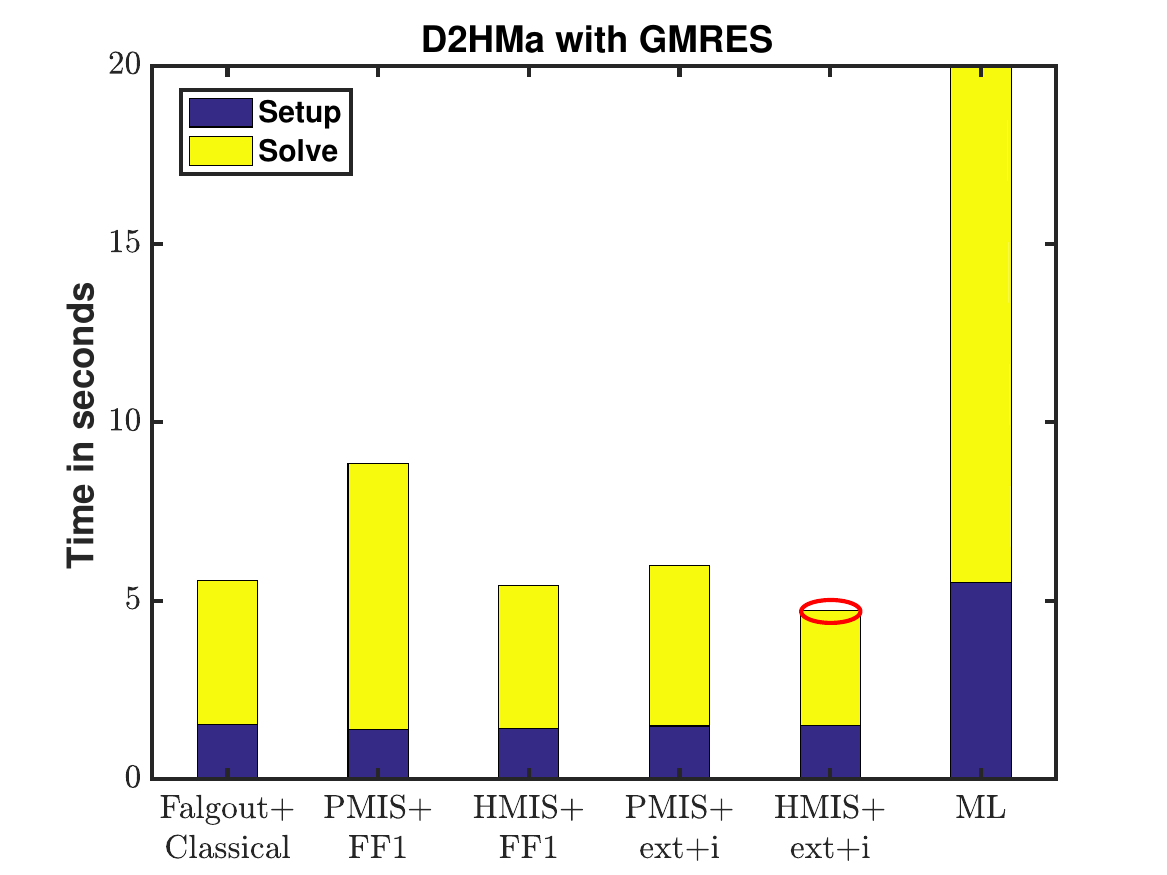}
\par\end{center}%
\end{minipage}\hfill{} %
\begin{minipage}[t]{0.45\textwidth}%
\begin{center}
\includegraphics[width=1\textwidth]{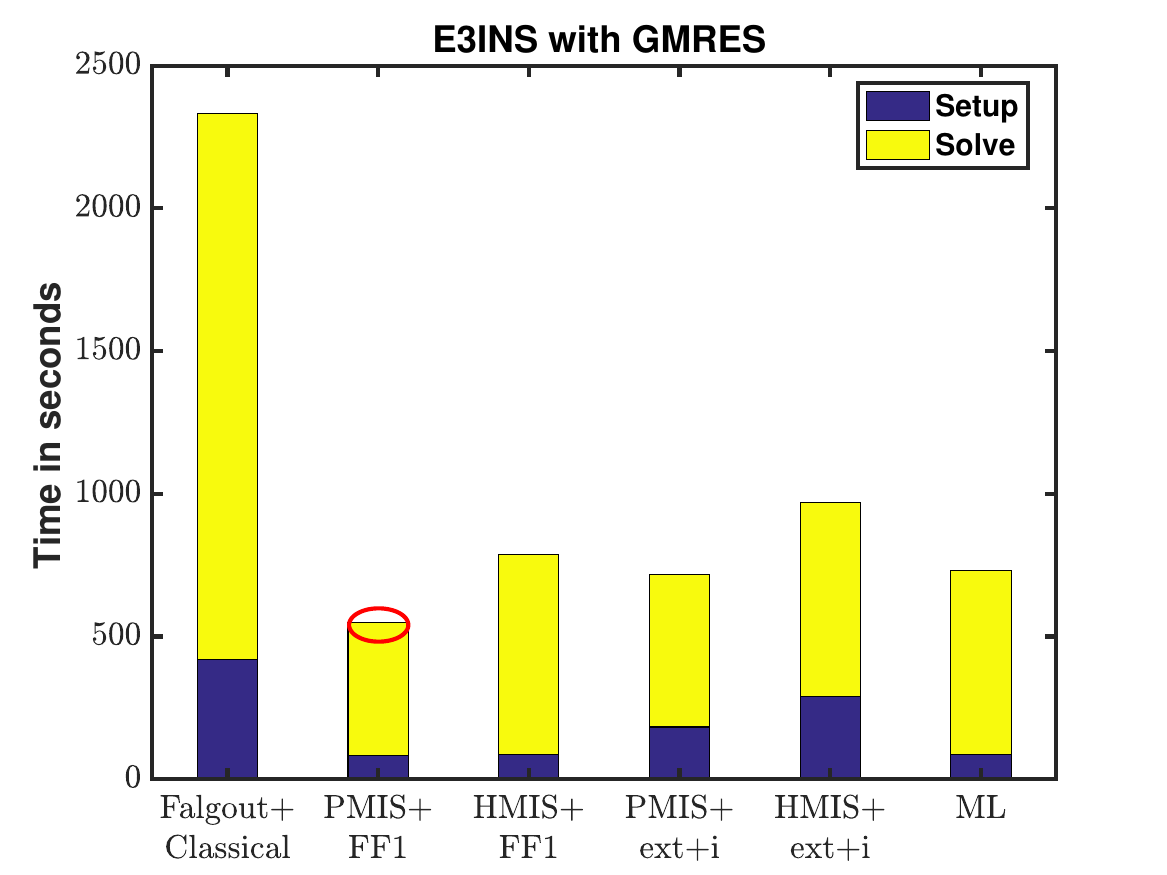}
\par\end{center}%
\end{minipage}

\caption{\label{fig:Hypre-options-ML}Comparisons of setup and solve times
of BoomerAMG with five coarsening+interpolation strategies and ML.
Circled bars indicate best performance.}
\end{figure}

In the above tests, BoomerAMG was clearly the winner. However, different
coarsening and interpolation techniques in BoomerAMG may lead to drastically
different performance, especially for 3D problems. Figure~\ref{fig:Hypre-options-ML}
presents a more in-depth comparison of the three coarsening strategies,
namely Falgout, PMIS, and HMIS, along with ML. For Falgout coarsening,
the ``classical'' interpolation is the most effective. For PMIS
and HMIS, we consider both the extended+i interpolation and the modified
FF interpolation (FF1). In all the cases, the interpolation matrices
were truncated to 4 elements per row, as recommended by hypre manual
\cite{hypre-user}. We used the default values for the other parameters.
It can be seen that the Falgout coarsening, delivered good performance
in 2D, but it underperformed ML for 3D problems. However, PMIS+FF1
and HMIS+FF1 delivered significantly better performance than ML for
3D problems. Between PMIS+FF1 and HMIS+FF1, their performances were
comparable, but HMIS outperformed PMIS for ill-conditioned problems.
In addition, FF1 interpolation outperformed the extended+i interpolation
for seven out of eight cases, probably because the truncation for
extended+i is not very effective. Note that for E3INS, although ML
outperformed HMIS+FF1, BoomerAMG with PMIS+FF1 still outperformed
ML. Therefore, we consider BoomerAMG with HMIS+FF1 or PMIS+FF1 as
the overall winner.

\textbf{}

We note that in the literature, it is sometimes reported that smoothed
aggregation works better than classical AMG for 3D elasticity problems,
which are generalizations of Poisson equations to vector-valued functions.
We did not consider elasticity problems in this work, because their
linear systems are symmetric; we refer readers to \cite{baker2010improving}
for more discussions on AMG for elasticity. However, we comment that
the aforementioned comparison for elasticity problems refers to the
classical Ruge-Stüben coarsening, which, like Falgout coarsening in
BoomerAMG, does not work well for 3D problems or for vector-valued
PDEs. An advantage of smoothed aggregation is that it can naturally
take into consideration the coupling of the variables in vector-valued
PDEs. However, in Table~\ref{tab:results_outline}, the results for
incompressible Navier-Stokes equations indicate that HMIS coarsening
in BoomerAMG compares favorably to smoothed aggregation in ML for
3D vector-valued PDEs.

It is also worth noting that HMIS coarsening is by no means foolproof.
In particular, GMRES failed to converge with HMIS coarsening and FF1
(or extended+i) interpolation, although it converged rapidly with
PMIS coarsening and FF1 interpolation. In addition, GMRES with BoomerAMG
could not solve any of the saddle-point problems, because the smoothers
could not handle large zero diagonal blocks. Hence, there is still
room for improvements for BoomerAMG in terms of both coarsening techniques
and smoothers. 

\section{Conclusions and Discussions\label{sec:Conclusions-and-Future}}

In this paper, we presented a systematic comparison of a few preconditioned
Krylov subspace methods, including GMRES, TFQMR, BiCGSTAB, and QMRCGSTAB,
with Gauss-Seidel, several variants of ILU, and AMG as right preconditioners.
We compared these methods at a theoretical level in terms of mathematical
formulations and operation counts. More importantly, we reported empirical
comparisons in terms of convergence, runtimes, and asymptotic complexity
with respect to problem sizes. To facilitate this comparative study,
we generated a number of large benchmark problems from various numerical
methods for a range of PDEs. Overall, our results show that GMRES
with multigrid as right preconditioner tends to be the most effective
for problems that are well-conditioned and are not saddle-point-like.
This is because right-preconditioned GMRES without restart minimizes
the 2-norm of the residuals within the Krylov subspace, and with an
effective preconditioner such as AMG, its cost of orthogonalization
is minimal when the iteration count is low.

Among AMG preconditioners, we observe that BoomerAMG in hypre, whose
algorithms are extensions of classical AMG with more sophisticated
coarsening and interpolation, tends to converge faster than ML, whose
algorithms are variants of smoothed aggregation.

Based on these results, we make the following primary recommendation:
\begin{quotation}
\emph{For large, moderately conditioned, non-saddle-point problems,
use GMRES with BoomerAMG as right preconditioner, with HMIS (or PMIS)
coarsening and FF1 interpolation.}
\end{quotation}
We emphasize the importance of coarsening and interpolation techniques,
especially for 3D problems. For BoomerAMG, we recommend HMIS coarsening
with FF1 interpolation. It delivers similar performance as PMIS+FF1
for well-conditioned systems, but it may outperform the latter significantly
for ill-conditioned systems. However, HMIS+FF1 may fail sometimes,
and in those cases one can try PMIS+FF1. The default in hypre before
version 2.11.2 was Falgout coarsening with classical interpolation,
which works well only for 2D problems; the new default in version
2.11.2 is HMIS coarsening with extended+i interpolation, which underperformed
HMIS+FF1 in our comparisons.

The easiest way to leverage the above recommendation is to use existing
software packages. PETSc \cite{petsc-user-ref} is an excellent choice,
since it supports both left and right preconditioning for GMRES and
BiCGSTAB, and it supports BoomerAMG with various options. Note that
PETSc uses left preconditioning by default, so we recommend explicitly
setting the option to use right preconditioning to avoid premature
or delayed termination for large systems. Note that with an effective
multigrid preconditioner, BiCGSTAB often converges almost as smoothly
as GMRES. Therefore, assuming a multigrid precondition is available,
BiCGSTAB can be used in place of GMRES, especially if right-preconditioned
GMRES is unavailable but right-preconditioned BiCGSTAB is (such as
in MATLAB).

Our comparative study also draws attention to QMRCGSTAB. Like BiCGSTAB,
QMRCGSTAB enjoys a three-term recurrence, so it may outperform restarted
GMRES if many iterations are needed. Furthermore, QMRCGSTAB converges
much more smoothly than BiCGSTAB, and its extra cost is negligible.
However, QMRCGSTAB is not available in PETSc or MATLAB. We plan to
make our implementations publicly available in the future.

A key component of multigrid preconditions is the smoother, which
is typically based on some variant of stationary iterative methods
or ILU0. These smoothing techniques are not robust for ill-conditioned
problems, and they fail for saddle-point-like problems. In the absence
of robust smoothers for multigrid preconditioners, we make the following
complementary recommendation to practitioners:
\begin{verse}
\emph{For ill-conditioned or saddle-point-like problems, use GMRES
with multilevel ILU as a right preconditioner.}
\end{verse}
We do not recommend the use of BiCGSTAB with multilevel ILU, but a
right-preconditioned QMRCGSTAB is advisable. If MILU is unavailable,
ILUTP may be used for moderate-sized systems; however, parameter tuning
for ILUTP is problematic for large-scale problems.

In this paper, we did not consider geometric multigrid (GMG) methods.
GMG often suffices as a standard-alone solver, and it typically outperforms
preconditioned Krylov-subspace methods significantly if applicable;
see e.g. \cite{LJM14HYGA} for comparisons of GMG, AMG, and a hybrid
multigrid method. If GMG is not robust enough as a standalone solver,
our recommendations regarding AMG preconditioners also holds to GMG
as right preconditioners. For ill-conditioned systems, GMG, as a solver
or preconditioner, may have significant advantages over AMG; for example,
for the very ill-conditioned Helmholtz equations, our team have observed
GMG outperforming AMG by orders of magnitude, which we plan to report
elsewhere. For incompressible Navier-Stokes equations, for which AMG
did not outperform Gauss-Seidel or ILU0 in our test, GMG can still
significantly outperform them; we refer readers to \cite{wesseling2001geometric}
for GMG for fluid-dynamics applications.

In terms of future work for research on preconditioners, one of the
critical areas is the robust smoothers for saddle-point-like problems.
This applies to both AMG and GMG. Although some customized multigrid
preconditioners have been developed (see e.g. \cite{Adams04AMG}),
they are not very general. Multilevel ILU is the most robust approach
in the state of the art, but it is not yet widely available, especially
in terms of parallel implementations. More importantly, MILU significantly
underperforms multigrid methods for large-scale systems. Hence, we
pose this open problem: \emph{Develop preconditioners that are as
robust as multilevel ILU but are as efficient and scalable as BoomerAMG
and GMG}. This likely will require some hybrid approaches.

One limitation of this work is that we did not compare parallel performance
and the scalability of the iterative methods with respect to the number
of cores. This omission was necessary to make the scope of this study
manageable, and also due to a lack of efficient parallel implementations
of ILU. For engineering applications that require only a small number
of cores, our primary recommendations are still relevant, in that
the MPI-based parallel implementation of right-preconditioned GMRES
and BiCGSTAB are available in PETSc, and both BoomerAMG and ML support
MPI. Hence, PETSc is an excellent choice for solving non-saddle-point
problems on distributed-memory machines. Parallel implementation of
ILU is not available in PETSc as of v3.9. A recent algorithm of iterative
computation of ILU seems to be promising \cite{chow2015fine}, but
parallel multilevel ILU is still an open problem. Another limitation
of PETSc is that it does not support multithreading. Some OpenMP and
CUDA-based implementations are available, such as the commercial version
of Paralution \cite{Paralution}, which seems to support only left
preconditioning as of v1.1. Eigen 3 \cite{eigenweb} is an open-source
software, and its BiCGSTAB supports OpenMP and right-preconditioning.
Finally, we note that the conclusions in this study cannot, and should
not, be extrapolated to extreme-scale applications, such as leading-edge
scientific or defense applications with billions of unknowns on hundreds
of thousands of cores. These applications would require more sophisticated
and more customized preconditioners, such as the hybrid multigrid
method in \cite{rudi2015extreme}. 

\section*{Acknowledgements}

The authors acknowledge supports by DoD-ARO under contract \#W911NF0910306
and by Argonne National Laboratory under Contract DE-AC02-06CH11357
for the SciDAC program funded by the Office of Science, Advanced Scientific
Computing Research of the U.S. Department of Energy. Results were
obtained using the LI-RED computer system at the Institute for Advanced
Computational Science of Stony Brook University, funded by the Empire
State Development grant NYS \#28451.

\bibliographystyle{abbrv}
\bibliography{refs/refs,refs/multigrid,refs/compkrylov_refs}

\end{document}